\documentclass{amsart}
\usepackage{amsmath} 
\usepackage{amssymb}
\usepackage{mathrsfs}

\newtheorem{theorem}{Theorem}[section] 
\newtheorem{claim}[theorem]{Claim}
\newtheorem{lemma}[theorem]{Lemma} 
 
\newtheorem{conclusion}[theorem]{Conclusion}
\newtheorem{observation}[theorem]{Observation}

\theoremstyle{definition}
\newtheorem{definition}[theorem]{Definition}
\newtheorem{convention}[theorem]{Convention}
\newtheorem{example}[theorem]{Example}
\newtheorem{fact}[theorem]{Fact}
\newtheorem{se}[theorem]{Subexample}

\newtheorem{discussion}[theorem]{Discussion}

\theoremstyle{remark}
\newtheorem{remark}[theorem]{Remark}

\newtheorem{notation}[theorem]{Notation}

\newcommand{\tp}{{\rm tp}}
\newcommand{\tr}{{\rm tr}}
\newcommand{\Av}{{\rm Av}}
\newcommand{\inv}{{\rm inv}}
\newcommand{\INV}{{\rm INV}}

\newcommand{\qf}{{\rm qf}}
\newcommand{\lx}{{\rm lx}}
\newcommand{\oor}{{\rm or}}
\newcommand{\org}{{\rm org}}
\newcommand{\otp}{{\rm otp}}
\newcommand{\proj}{{\rm proj}}
\newcommand{\stat}{{\rm stat}}

\newcommand{\Reg}{{\rm Reg}}

\newcommand{\OTOP}{{\rm OTOP}}
\newcommand{\Card}{{\rm Card}}

\newcommand{\BA}{{\rm BA}}
\newcommand{\CH}{{\rm CH}}

\newcommand{\sub}{{\rm sub}}

\newcommand{\acc}{{\rm acc}}

\newcommand{\id}{{\rm id}}

\newcommand{\PC}{{\rm PC}}
\newcommand{\EM}{{\rm EM}}
\newcommand{\EF}{{\rm EF}}

\newcommand{\Min}{{\rm Min}}
\newcommand{\Dom}{{\rm Dom}}
\newcommand{\Rang}{{\rm Rang}}
\newcommand{\rang}{{\rm rang}}
\newcommand{\isoI}{{\dot{\bbI}}}

\newcommand{\rest}{{\restriction}}

\newcommand{\wilog}{{\rm without loss of generality}}
\newcommand{\Wilog}{{\rm Without loss of generality}}

\newcommand{\then}{{\underline{then}}}
\newcommand{\when}{{\underline{when}}}
\newcommand{\Then}{{\underline{Then}}}

\newcommand{\Iff}{{\underline{iff}}}
\newcommand{\mn}{{\medskip\noindent}}
\newcommand{\sn}{{\smallskip\noindent}}

\newcommand{\bbI}{{\mathbb I}}

\newcommand{\cC}{{\mathscr C}}

\newcommand{\cD}{{\mathscr D}}

\newcommand{\cH}{{\mathscr H}}

\newcommand{\varp}{{\varepsilon}}

\newcommand{\bbL}{{\mathbb L}}
\newcommand{\bbG}{{\mathbb G}}

\newcommand{\cM}{{\mathscr M}}

\newcommand{\bbP}{{\mathbb P}}
\newcommand{\cP}{{\mathscr P}}

\newcommand{\bbS}{{\mathbb S}}

\newcommand{\cW}{{\mathscr W}}

\newcommand{\cL}{{\mathscr L}}
 
\newcommand{\cU}{{\mathscr U}}

\newcommand{\cf}{{\rm cf}}

\newcount\skewfactor
\def\mathunderaccent#1#2 {\let\theaccent#1\skewfactor#2
\mathpalette\putaccentunder}
\def\putaccentunder#1#2{\oalign{$#1#2$\crcr\hidewidth
\vbox to.2ex{\hbox{$#1\skew\skewfactor\theaccent{}$}\vss}\hidewidth}}

\newenvironment{PROOF}[2][\proofname.]
   {\begin{proof}[#1]\renewcommand{\qedsymbol}{\textsquare$_{\rm #2}$}}
   {\end{proof}}

\begin{document}

\title {General non-structure theory}
\author {Saharon Shelah}
\address{Einstein Institute of Mathematics\\
Edmond J. Safra Campus, Givat Ram\\
The Hebrew University of Jerusalem\\
Jerusalem, 91904, Israel\\
 and \\
 Department of Mathematics\\
 Hill Center - Busch Campus \\ 
 Rutgers, The State University of New Jersey \\
 110 Frelinghuysen Road \\
 Piscataway, NJ 08854-8019 USA}
\email{shelah@math.huji.ac.il}
\urladdr{http://shelah.logic.at}
\thanks{The author thanks Alice Leonhardt for the beautiful typing.
This is a revised version
of \cite[Ch.III,\S1-\S3]{Sh:300}, has existed (and somewhat revised) for
many years.
% 10.05.30 Hopefully will be Ch.III of the authors book on
% 10.05.30 non-structure theory, but it is by itself
% 10.05.30 \cite{Sh:E59}.
Was mostly ready in the early nineties, and
public to some extent.
For the sake of \cite{LwSh:687}  we add the part of \S1 from
\ref{1.22new}.  For the sake of \cite{ShUs:928} we add in the end of \S2.
% 10.05.31  bdoq   We \underline{may sometime} add some word after
Recently this work was used and continued in Farah-Shelah \cite{FaSh:954}.
This was written as Chapter III of the book \cite{Sh:e},
which hopefully will materialize some day, but in meanwhile 
it is \cite{Sh:E59}.
% 10.05.30 end of addition
The intentions were: \cite{Sh:E58} (revising \cite{Sh:229}) for Ch.I,
\cite{Sh:421} for Ch.II, \cite{Sh:E59} for Ch.III, \cite{Sh:309} for Ch.IV,
\cite{Sh:363} for Ch.V, \cite{Sh:331} for Ch.VI, \cite{Sh:511} for
Ch.VII, \cite{Sh:E60}, a revision of \cite{Sh:128} for Ch.VIII,
\cite{Sh:E62} for the appendix, and probably \cite{Sh:757}, \cite{Sh:384},
\cite{Sh:482}, and \cite{Sh:800}.
References like
\cite[3.7=Lc2]{Sh:E62} means that c2 is the label of 3.7 in
\cite{Sh:E62}, will only help the author if changes in the paper
\cite{Sh:E62} will change the number.}

% 10.05.30 end footnote, second } - for the title



% Previous revision: January 8, 2016

\date {January 12, 2016}

\begin{abstract}
The theme of the first two sections, is to prepare the
framework of how from a ``complicated'' family of index
models $I \in K_1$ we build many and/or complicated structures in
a class $K_2$. The index models are characteristically linear orders,
trees with $\kappa+1$ levels (possibly with linear order on the
set of successors of a member) and linearly ordered graph,
for this we phrase relevant complicatedness properties
(called bigness).
% 10.05.30

We say when $M \in K_2$ is represented in $I \in K_1$.
We give sufficient conditions when $\{M_I:I \in K^1_\lambda\}$ is 
complicated where for each $I\in K^1_\lambda$ we build $M_I \in K^2$ 
(usually $\in  
% 10.08.02
K^2_\lambda$) represented in it and
reflecting to some degree its structure (e.g. for
$I$ a linear order we can build a model of an unstable
first order class reflecting the order). If we understand
enough we can even build e.g. rigid members of $K^2_\lambda$.
% 10.05.30

Note that we mention ``stable'', ``superstable'', but in a self
contained way, using an equivalent definition which is useful
here and explicitly given. We also frame the use of
generalizations 
% 10.05.30
of Ramsey and Erd\" os-Rado theorems to get models in which any $I$
from the relevant $K_1$ is reflected.
We give in some detail how this may apply to the class
of separable reduced Abelian $\dot p$-group and how we get relevant
models for ordered graphs (via forcing).

In the third section we show stronger results
concerning linear orders. If for each linear order $I$
of % 10.05.30
 cardinality $\lambda > \aleph_0$ we can attach a model
$M_I \in K_\lambda$ in which the linear order can be embedded
such that for enough cuts of $I$, their being omitted is reflected
in $M_I$, then there are $2^\lambda$ non-isomorphic cases. 

But in the end of the second section we show how the
results on trees with $\omega+1$ levels
% 10.05.30 
(on which concentrate \cite{Sh:331}
gives results on linear ordered (not covered by \S3),
on trees with $\omega + 1$ levels see \cite{Sh:331}. 
% 10.05.30
To get more we prove explicitly more on such trees.
Those will be enough for results in model theory of
Banach space of Shelah-Usvyatsov \cite{ShUs:928}.
\end{abstract}

\maketitle
\numberwithin{equation}{section}
\setcounter{section}{-1}
\newpage

\section {Introduction}

The main result presented in this paper is (in earlier proofs we have it
only in ``most" cases):

\begin{theorem}
\label{0.1new}
If $\psi \in \bbL_{\chi^+,\omega},\varphi(\bar x,\bar y) \in 
\bbL_{\chi^+,\omega},\ell g(\bar x) = \ell g(\bar y)=\partial$ and 
$\psi$ has the $\varphi(\bar x,\bar y)$--order property (see 
Definition \ref{1.2}(5)) \then \ $\isoI 
% 10.05.30
(\lambda,\psi)=2^\lambda$ provided that for example:
$\lambda \ge \chi+ \aleph_1,\partial<\aleph_0$ or $\lambda=\lambda^\partial+
\chi + \partial^++ \aleph_1$ or $\lambda > \chi + \partial^+$ 
or $\lambda^{\partial^+} < 2^\lambda,\lambda \ge \chi$.
\end{theorem}

\begin{PROOF}{\ref{0.1new}}
By \ref{3.14}(2), 
% 10.08.02
clause $(b)$ of \ref{3.14}(2) holds. When $\lambda \ge \chi
+ \aleph_1,\partial < \aleph_0$, by Theorem \ref{3.10}(3), 
$\dot{\bbI}(\lambda,\psi)= 2^\lambda$. 

So we can assume that $\lambda \ge \chi$ and $\partial \ge \aleph_0$.
When $\lambda^\partial=\lambda$ or $\lambda^{(\partial^+)} <
2^\lambda$ the conclusion holds by \ref{3.11}(a), \ref{3.11}(b),
respectively, using $\kappa=\partial^+$ and the existence of 
such models follows from \ref{1.11} they are as required by \ref{3.4}(4).
When $\lambda>\chi+\partial^+$ the conclusion holds by
\ref{3.14}(1).
So we are done. 
\end{PROOF}

Note that although some notions connected to stability appear, they are not
used in any way which require knowing them: we define what we use and at
most quote some results. In fact, the proof covered problems with no
(previous) connection to stability. For understanding and/or checking, the
reader does not need to know the works quoted below: they only help to see
the  background. Note that 

For later chapters (please give specific numbers)
\S2 is essential to some of % 10.08.02 i
the later parts of non-structure (see \cite{Sh:309},
\cite{Sh:331} \cite{Sh:511}) 
% 10.05.30 
them but not \S1 or \S3  still 
% 10.08.02 
but better read \ref{1.1}-\ref{1.7}.

Generally the construction of many models (up to isomorphism in this paper)
in $K_\lambda$ ($=:\{M\in K:\|M\|=\lambda\})$ goes as follows. 
We are given a class $K$ of models (with fix vocabulary), and we are 
trying to prove that $K$ has many complicated members. 
To help us, we have a class $K^1$ of ``index models''
(this just indicates their role; supposedly they are well understood; they
usually are linear orders or a class of trees). By the ``non-structure
property of $K$'', for some formulas $\varphi_\ell$ (see below), for every
$I \in K^1_\lambda$ there is $M_I\in K_\lambda$ and $\bar a_t\in M_I$ for $t
\in I$, which satisfies (in $M_I$) some instances of $\pm\varphi_\ell$.

We may demand on $M_I$:
\mn
\begin{enumerate}
\item[$(0)$]  nothing more (except the restriction on the
  cardinality),
\sn
\item[$(1)$]  $\langle\bar a_t:t\in I\rangle$ behaves nicely: like a
skeleton (see \ref{3.1}(1)), \underline{or} even
\sn
\item[$(2)$]  $M_I$ is ``embedded'' in a model built from $I$ in a simple way
($\Delta$--represented; see Definition \ref{2.2}(c)), \underline{or}
\sn
\item[$(3)$]  $M_I$ is built from $I$ in a simple way, an 
%10.05.30  
the extreme case being $EM_\tau$
% 10.05.30
$(I,\Phi)$; see Definition \ref{1.6} where $\tau = 
% 10.08.02
\tau(M_I)$ of course. 
% 10.05.30
\end{enumerate}
\mn
Now even for (0) we can have meaningful theorems
(see % REF[IV]{4.2}
\cite[1.1]{Sh:309} and
% REF[IV]{4.4}
\cite[1.3]{Sh:309}); but we cannot have all we would naturally
like to have --- see 
% REF[IV]{4.6}
\cite[1.8]{Sh:309} (i.e., we cannot prove much better results in
this direction, as shown by a consistency proof).

Though it looks obvious by our formulation, 
% 10.08.02
experience shows that we must stress that the
formulas $\varphi_\ell$ need not be first order, they just have to have the
right vocabulary (but in results on ``no $M_i$ embeddable in $M_j$'' this
usually means embedding preserving $\pm \varphi_\ell$ (but see the proof of
% REF[VI]{7.15}
\cite[3.22(2)]{Sh:331}. 
% 10.08.02  (2)).
So they are just properties of sequences in the
structures we are considering preserved by the morphism we have in mind.

Another point is that though it would be nice to prove

\[
[I \not\cong J\quad\Rightarrow\quad M_I\not\cong M_J];
\]

\mn
this does not seem realistic. What we do is to construct a family

\[
\{I_\alpha:\alpha<2^\lambda\}\subseteq K^1_\lambda
\]

\mn
such that for $\alpha \neq \beta$, in a strong sense $I_\alpha$ is not
isomorphic to (or not embeddable into) $I_\beta$ (see \ref{2.3}, \ref{3.4},
% REF[VI]{7.1}
     \cite[1.1]{Sh:331}, % REF[VI]{7.3}
     \cite[1.4]{Sh:331}), such that now we have $M_{I_\alpha},
M_{I_\beta}$ not isomorphic for $\alpha \ne \beta$. We are thus led to the
task of constructing such $I_\alpha$'s, which, probably unfortunately,
splits to cases according to properties of the cardinals involved. Sometimes
we just prove $\{\alpha :M_\alpha\cong M_\beta\}$ is small for each $\beta$.

A point central to
\cite{Sh:E58}, \cite{Sh:421}, \cite{Sh:511},\cite{Sh:384} and
\cite{Sh:482}
% 10.05.30 smallbox{please check these numbers}
           % chapters I, II, VII, X, XI
but incidental here, is the construction of a model 
which is for example rigid or has few endomorphisms,
etc. In particular in \cite{Sh:511} we could use linear order for ``the
gluing".

The methods here can be combined with \cite{Sh:220} or \cite{Sh:188} to get
non-isomorphic ${\bbL}_{\infty,\lambda}$--equivalent models of cardinality
$\lambda$; Instead ``${\bbL}_{\infty,\lambda}$-equivalent
non-isomorphic model of $T$'' we can consider equivalence by stronger
games, e.g. $\EF_{\alpha,\lambda}$-equivalence started in
Hyttinen-Tuuri \cite{HyTu91}, and then Hyttinen-Shelah \cite{HySh:474},
\cite{HySh:529}, \cite{HySh:602}; See
V\"a\"an\"anen % 10.05.30 Vaananen
\cite{Va95} or such games.

In the next few paragraphs we survey the results of this paper. In this
survey we omit some parameters for 
% 10.05.30 
at various defined notions. These parameters
are essential for an accurate statement of the theorems. We suppress them
here trying to make it easier reading while still communicating essential
points.

In \S1 we mainly represent {E.M.} models. This is how in a natural way we
construct a model from an ``index model''. The proof of existence many times
rely on partition theorems. We give definition, deal with the framework,
quote important cases, and present general theorems for getting the {E.M.}
models, i.e., the templates; we then, as an example, deal with random graphs
for theories in ${\bbL}_{\kappa^+,\omega}$.

In \S2 we discuss a method of ``representability'' (from
\cite{Sh:136}). This is a natural way to get for ``a model gotten from an
index model $I$" that ``$I$ is complicated" implies ``$M$ is
complicated". We discuss applications (to separable reduced Abelian
$\dot p$-groups and Boolean algebras), but the aim is to explain; 
full proofs of full results will appear later (see
% 10.05.30  VI, \S3, VII,
\cite[\S3]{Sh:331}, \cite{Sh:511} 
% 10.07.23
respectively). We
introduce two strongly contradictory notions, the $\Delta$--representability
of a structure $M$ in the ``free algebraq" (i.e., ``polynomial algebra") of
an index model (Definition \ref{2.2}) and the $\varphi(\bar x, \bar
y)$--unembeddability of one index model in another. Now, to show that a
class $K$ has many models it suffices if for some formula $\varphi$, one
first shows that (a) an index class $K_1$ has many pairwise
$\varphi$--unembeddable structures, second that (b) for each 
$I \in K_1$, 
%10.08.02
there is a model $M_I$ which is $\Delta$--representable in the free algebra on
$I$, and finally that (c) if $M_I \cong M_J$ and $M_J$ is
$\Delta$--represented in the free algebras on $J$ then $I$ is
$\varphi$--embeddable in $J$.

However, for building for example a rigid model of cardinality $\lambda$, it is
advisable to use $\langle I_\alpha:\alpha<\lambda\rangle$ such that
$I_\alpha$ is $\varphi$--unembeddable into $\sum\limits_{\beta\ne\alpha}
I_\beta$. (See \ref{2.6}, \ref{2.7}, more in \cite{Sh:511}).
Generally having suitable sequence of $I \in K_1$ 
% 10.08.02  % 10.05.30
is expressed by ``$K_1$  
% 10.08.02
has a suitable bigness property". 
% 10.05.30

Now, \S3 does not depend on \S2. The point is that in this section our non-
isomorphisms proofs are so strong that they do not need
``representability", we use a much weaker property. In \S3 we extend and
simplify the argument showing that an unstable first order theory $T$ has
$2^\lambda$ models of cardinality $\lambda$ if $\lambda \ge |T|
+\aleph_1$. Rather
than constructing Ehrenfeucht--Mostowski models we consider a weaker notion
--- that a linear order $J$ indexes a weak $(\kappa,\varphi)$-skeleton like
sequence in a model $M$. In this section, $K_1$ is the class of linear
orders. The formula $\varphi(\bar x,\bar y)$ need not be first order and
after \ref{3.10} may have infinitely many arguments. Most significantly we
make no requirement on the means of definition of the class $K$ of models
(for example first order, $\bbL_{\infty,\infty}$, etc). We require only that
for each linear order $J$ there are an $M_J \in K$ and a sequence $\langle
\bar a_s:s\in J\rangle$ which is weakly $(\kappa,\varphi)$--skeleton like in
$M_J$.

% 10.05.30 for the book: if you get lost somewhere in {\S}3, you can
% 10.05.30 jump to the next chapter.  \smallbox{Please reformulate}
Note
that having bigness properties for $K^\kappa_\tr$ implies the ones for
$K_{\oor}$ see \ref{2.21}, Ehrenfeucht and Mostowski \cite{EhMo56}
built what are here $EM_\tau 
% 10.05.30 
(I,\Phi)$
for $I$ a linear order and first order $T$
where $\tau = \tau_T$. 
% 10.05.30
Ehrenfeucht \cite{Eh57}, \cite{Eh58} (and Hodges in \cite{Ho73}
improved the set theoretic  assumption) proved that if
$T$ has the property $(E)$ then it has at least two non-isomorphic models
(this property is a precursor of being unstable). Recall that
the property $(E)$  % 10.05.30 means
says that:
some a formula $R(x_1,\ldots x_n)$ is asymmetric on some
infinite subset of some model of $T$; note that $(E)$
is not equivalent to being unstable as
 the theory of random graphs fail it. Morley
\cite{Mo65} prove that for well ordered $I$, the model is stable in
appropriate cardinalities, to prove that non-totally transcendental countable
theories are not categorical in any $\lambda>\aleph_0$. See more in
\cite[VII,VIII]{Sh:c}; by it if
$T\subseteq T_1$ are unstable, complete first order and $\lambda\geq
|T_1|+\aleph_1$ then $T_1$ has $2^\lambda$ models of
cardinality  
% 10.05.30
$\lambda$ 
% 10.08.02
with pairwise non-isomorphic reducts 
% 10.08.02
to $\tau_T$. On the cases for $\bbL_{\chi^+,\omega},\lambda > \chi$,  
% 10.08.02
see Grossberg-Shelah \cite{GrSh:222}, \cite{GrSh:259} which continue
\cite{Sh:11}.

This paper is a revised version of sections \S1,\S2,\S3
of chapter III of \cite{Sh:300}.
\newpage

\section {Models from Indiscernibles}

We survey here \cite[Ch.VIII,\S3]{Sh:a}, which was the starting point for the
other works appearing or surveyed in this paper and \cite{Sh:309},
\cite{Sh:363}.  So we concentrate on
building many models for first order theories, using {E.M.} models, i.e., in
all respects taking the easy pass. Our aim there was

\begin{theorem}
\label{1.1}
If $T$ is a complete first order theory, unstable and $\lambda\geq |T|+
\aleph_1$, \then \, $\dot{\bbI}(\lambda,T)=2^\lambda$,
\end{theorem}

\noindent
where
\begin{definition}
\label{1.2}
$T$ is unstable \when \,
% 10.08.02 \underline{when}
 % 10.08.02  % 10.05.30 if
for some first order formula $\varphi(\bar x,\bar y)$
($n = \ell g(\bar x) = \ell g(\bar y)$)
in the vocabulary $\tau_T$ of $T$ of course, 
% 10.05.30
for every $\lambda$ there is a model $M$ of
$T$ and $\bar a_i \in {}^n M$ for $i<\lambda$ such that

\[
M \models\varphi[\bar a_i,\bar a_j] \text{ Iff } i<j\ (<\lambda).
\]
\end{definition}

\begin{definition}
\label{1.3new}
For a theory $T$ and vocabulary $\tau \subseteq \tau_T$, 
% 10.08.02

$\isoI(\lambda,T) = \text{ the number of models of } T$ of cardinality
$\lambda$, up to isomorphism,

$\isoI_\tau(\lambda,T) = \text{the number of } \tau \text{-reducts of
models of } T \text{ of cardinality } \lambda$,\\

\hskip25pt up to isomorphism.
\end{definition}

\begin{definition}
\label{1.4new}
1)  For a class $K$ of models and set $\Delta$ of formulas:
\mn
\begin{enumerate}
\item[${{}}$]  $\isoI(\lambda,K) =$ the number of models in $K$ 
of cardinality $\lambda$ up to isomorphism,
\sn
\item[${{}}$]  $\isoI(K) =$ the number of models in $K$ up to
  isomorphism,
\sn
\item[${{}}$]  $\dot I \dot E_\Delta(\lambda,K) = \sup\{\mu$: there
  are $M_i \in K_\lambda$, for $i < \mu$,  such that for $i \ne j$ 
there is no $\Delta$-embedding of $M_i$ to $M_j\}$.
\end{enumerate}
\mn
see part (2); and we may write $\tau$ instead
$\Delta = \bbL(\tau_K)$, may omit $\Delta$ when it is $\bbL(\tau_M)$.

\noindent
2) $f:M \longrightarrow N$ is a $\Delta$-embedding (of $M$ into $N$)
\Iff \, ($f$ is a function from $|M|$ into $|N|$ and) for every
$\varphi(\bar x)\in\Delta$ and $\bar a \in {}^{\ell g(\bar a)}|M|$,
we have: 

\[
M \models \varphi [\bar a] \Rightarrow N \models\varphi[f(\bar a)].
\]

\mn
(so if $(x \ne y) \in \Delta$ then $f$ is one to one).
\end{definition}

\begin{definition}
\label{1.5new}
1)  A sentence $\psi\in \bbL_{\chi^+,\omega }$ is $\partial$-unstable
\Iff \, there are $\alpha<\partial$ and a formula $\varphi(\bar x,\bar
y)$ from $\bbL_{\chi^+,\omega}$ with $\ell g(\bar x) = \ell g(\bar y)
=\alpha$ such that $\psi$ has the $\varphi$-order property, i.e., 
for every $\lambda$ there is a model
$M_\lambda$ of $\psi$ and a sequence $\bar a_\zeta$ 
of length $\alpha$ from $M_\lambda$ such that for $\zeta,\xi<\lambda$
we have

\[
M_\lambda\models\varphi[\bar a_\zeta,\bar a_\xi] \Leftrightarrow
\zeta < \xi.
\]

\mn
If $\partial = \aleph_0$ we may omit it.

\noindent
2)  For $\kappa$ regular and $T$ first order, we say $\kappa<\kappa(T)$
\Iff \, there are first order formulas $\varphi_i(\bar x,\bar y_i)
\in \bbL(\tau_T)$ for $i<\kappa$ and for every
$\lambda$ there is a model $M_\lambda$ of $T$ and for $i \le \kappa,
\eta \in {}^i \lambda$ a sequence $\bar a_\eta$ from $M_\lambda$, with

\[
i < \kappa \Rightarrow \ell g(\bar a_\eta) = \ell g(\bar y_i)
\]

\[
i = \kappa \Rightarrow \ell g(\bar a_\eta) = \ell g(\bar x)
\]

\mn
such that:
if $\nu \in {}^i\lambda,\ \eta \in {}^\kappa\lambda,\ \nu
\vartriangleleft \eta$ then $M_\lambda \models \varphi_{i+1}
[\bar a_\eta,\bar a_{\nu \char 94 \langle \alpha\rangle}] \Leftrightarrow
\eta(i)=\alpha$. [We shall not use this except in 
\ref{1.8} below.]

\noindent
3) $T$, a first order theory, is unsuperstable if 
$\aleph_0 < \kappa(T)$ [but we shall use it only in \ref{1.8}].
\end{definition}
\bigskip

\centerline{$* \qquad * \qquad *$}
\bigskip

\begin{definition}
\label{1.5}
1) $\langle\bar a_t:t\in I\rangle$ is $\Delta$-indiscernible (in $M)$
\Iff \:
\mn
\begin{enumerate}
\item[$(a)$]  $I$ is an index model (usually linear order or tree),
  i.e., it can be any model but its role will be as an index set,
\sn
\item[$(b)$]  $\Delta$ is a set of formulas in the vocabulary of $M$
(i.e. in $\cL_{\tau(M)}$ for some logic $\cL$)
\sn
\item[$(c)$]  the $\Delta$-type in $M$ of $\bar a_{t_0} \char 94 \ldots
\char 94 \bar a_{t_{n-1}}$ for any $n<\omega$ and 
$t_0,\ldots t_{n-1}\in I$) depends only on the quantifier free type of
$\langle t_0,\ldots,t_{n-1}\rangle$ in $I$.
\end{enumerate}
\mn
Recall that the $\Delta$-type of $\bar a$ in $M$ is $\{\varphi(\bar x)\in
\Delta: M \models \varphi(\bar a)\}$, where $\bar a,\bar x$ are indexed by
the same set.  So the length of $\bar{a_t}$ depend just on 
the quantifier free type which $\ell g(\bar{a_t})$ realizes in $I$. 

If we allow $\varphi(\bar x) \in \Delta,\kappa > \alpha = \ell g(\bar x)\ge
\omega$ and we allow $\langle t_i:i<\alpha\rangle$ above, \then\ we 
say $(\Delta,\kappa)$-indiscernible.

\noindent
2) For a logic $\cL, ``\cL$-indiscernible" will mean
$\Delta$-indiscernible for the set of $\cL$-formulas in the
vocabulary of $M$.  If $\Delta,\cL$ are not mentioned we mean first
order logic.

\noindent
3) Notation: Remember that if $\bar t=\langle t_i:i<\alpha\rangle$ then
$\bar a_{\bar t}=\bar a_{t_0} \char 94 \bar a_{t_1} \char 94 \ldots$.
\end{definition}

Many of the following definitions are appropriate for counting the number of
models in a pseudo elementary class. Thus, we work with a pair of
vocabularies, $\tau \subseteq \tau_1$. Often $\tau_1$ will contain Skolem
functions for a theory $T$ which is $\subseteq \cL(\tau)$.

\begin{convention}
\label{1.5f}
For the rest of this section all predicates and function 
symbols have finite number of places (and similarly 
$\varphi(\bar x)$ means $\ell g(\bar x)<\omega)$.
\end{convention}

\begin{definition}
\label{1.6}
1)  $M = \EM(I,\Phi)$ \Iff \, for some vocabulary $\tau=\tau_\Phi=\tau(\Phi)$
(called $L^\Phi_1$ in \cite[Ch.VII]{Sh:a}) and sequences $\bar a_t(t\in I)$
we have:
\mn
\begin{enumerate}
\item[$(i)$]  $M$ is a $\tau_\Phi$-structure and is generated by 
$\{\bar a_t:t \in I\}$,
\sn
\item[$(ii)$]  $\langle\bar a_t:t \in I\rangle$ is quantifier 
free indiscernible in $M$, 
\sn
\item[$(iii)$]  $\Phi$ is a function, taking (for $n<\omega$) the quantifier
free type of $\bar t = \langle t_0,\ldots,t_{n-1}\rangle$ 
in $I$ to the quantifier free type of $\bar a_{\bar t} =
\bar a_{t_1} \char 94 \ldots \char 94 \bar a_{t_n}$ in $M$
(so $\Phi$ determines $\tau_\Phi$ uniquely).
\end{enumerate}
\mn
2) A function $\Phi$ as above is called a template and we say it is
proper for $I$ if there is $M$ such that $M =\EM(I,\Phi)$. 
We say $\Phi$ is proper for $K$ if $\Phi$ is proper for every 
$I\in K$, and lastly $\Phi$ is proper for $(K_1,K_2)$ if it is 
proper for $K_1$ and $\EM^1(I,\Phi)\in K_2$ for $I\in K_1$.

\noindent
3) For a logic $\cL$, or even a set $\cL$ of formulas in the
vocabulary of $M$, we say that $\Phi$ is almost $\cL$-nice (for $K$)
\Iff \, it is proper for $K$ and:
\mn
\begin{enumerate}
\item[$(*)$]  for every $I \in K,\langle \bar a_t:t\in I\rangle$ is
$\cL$-indiscernible in $EM(I,\Phi)$.
\end{enumerate}
\mn
4) In part (3), $\Phi$ is $\cL$-nice \Iff \, it is almost
$\cL$-nice and 
\mn
\begin{enumerate}
\item[$(**)$]  for $J\subseteq I$ from $K$ we have $\EM(J,\Phi)
\prec_{\cL} \EM(I,\Phi)$.
\end{enumerate}
\mn
5) In part (3) we say that $\Phi$ is $(\cL,\tau)$-nice \when \,
$\tau \subseteq \tau_\Phi$, it is almost $\bbL$-nice and (see
\ref{1.7}(1))
\mn
\begin{enumerate}
\item[$(***)$]  for $I \subseteq J$ from $K$ we have $\EM_\tau(J,\Phi) 
\prec_{\cL} \EM_\tau(I,\Phi)$.
\end{enumerate}
\end{definition}

\noindent
In the book \cite{Sh:a}, always $\bbL_{\omega,\omega}(\tau_\Phi)$-nice
$\Phi$ were used and $\EM(I,\Phi),\EM_\tau(I,\Phi)$ here
are $\EM^1(I,\Phi),\EM(I,\Phi)$ there. 
\begin{definition}
\label{1.7}
1) $\EM_\tau(I,\Phi) =\EM(I,\Phi) \rest \tau$, i.e., $\tau$-reduct
of $\EM(I,\Phi)$,  where $\tau\subseteq\tau_\Phi$.
We may omit $\tau$ when clear from the context and write $\EM(I,\Phi)$.
Saying ``an $\EM$-model will mean ``a model of the 
form $\EM_\tau(I,\Phi)$" where $\Phi,I,\tau$
are understood from the context. 

\noindent
2)  We identify $I \subseteq {}^{\kappa\ge} \lambda$ which is closed under
initial segments, with the model 
$(I,P_\alpha,\cap,<_\lx,\vartriangleleft)_{\alpha\leq\kappa}$, where

$P_\alpha=I \cap {}^\alpha\lambda$,

$\rho=\eta \cap \nu$ if $\rho = \eta \rest \alpha$ for the maximal
$\alpha$ such that $\eta \rest \alpha = \nu \rest \alpha$,

$\vartriangleleft=$ being initial segment of (including equality),

$<_{\lx}=$ the lexicographic order.

\noindent
3)  Similarly to (2), for any linear order $J$, every $I\subseteq
{}^{\kappa\ge}J$ which is closed under initial segments is identified with
$(I,P_\alpha,\cap,<_{\lx},\vartriangleleft)_{\alpha\le\kappa}$
($\le_{\lx}$ is still well defined).

\noindent
4)  $K^\kappa_{\tr}$ is the class of such models, i.e., models isomorphic
to such $I$, i.e., to $(I,P_\alpha,\cap,<_{\lx},\vartriangleleft)_{\alpha
\le \kappa}$ for some $I\subseteq {}^{\kappa\ge}J$
which is closed under initial segments, $J$ a linear order
($\tr$ stands for tree). We call $I$ standard if $J$ is an ordinal or at
least well ordered.

\noindent
5)  $K_{\oor}$ is the class of linear orders.
\end{definition}

\begin{remark}
\label{1.7A}
The main case here is $\kappa=\aleph_0$. We need such trees for $\kappa>
\aleph_0$, for example if we would like to build many $\kappa$-saturated models
of $T$, $\kappa(T)>\kappa$, $\kappa$ regular. If $\kappa(T) \le \kappa$ there
may be few $\kappa$--saturated models of $T$.
\end{remark}

\noindent
In \cite[Ch.VIII]{Sh:a} we have also proved:
\begin{lemma}
\label{1.8}
1)  If $T \subseteq T_1$ are complete first order theories, $T$ is
unstable as exemplified by $\varphi=\varphi(\bar x,\bar y)$, say $n =
\ell g(\bar x) = \ell g(\bar y)$, \then \, for some template $\Phi$ 
proper for the class of linear orders and nice for first order logic, 
$|\tau_\Phi|=|T_1|+ \aleph_0$ and for any linear order $I$ and 
$s,t \in I$ we have

\[
\EM(I,\Phi) \vDash \varphi[\bar a_s,\bar a_t] \text{ iff } I \vDash s<t.
\]

\mn
2) If $T \subseteq T_1$ are complete first order theories and $T$ is
unsuperstable, \then \, there are first order $\varphi_n(\bar x,
\bar y_n)\in \bbL(\tau_T)$ and a template $\Phi$ proper for every
$I \subseteq {}^{\omega \ge}\lambda$ such that for any such $I$ we
have:
\mn
\begin{enumerate}
\item[$(a)$]  $\eta \in {}^\omega\lambda,\nu\in {}^n\lambda$ implies
$\EM(I,\Phi) \models \varphi_n[\bar a_\eta,\bar a_\nu]$ \Iff \,$\eta
\rest n = \nu$
\sn
\item[$(b)$]  $\EM(I,\Phi)\models T_1$ and $\Phi$ is $\bbL_{\omega,\omega}
(\tau_\Phi)$-nice, $|\tau_\Phi|=|T_1|+\aleph_0$ (note that for $\eta_1,
\eta_2$ of the same length, $\eta_1 \ne \eta_2 \Rightarrow \bar a_{\eta_1}
\ne \bar a_{\eta_2})$\footnote{In fact $EM^1(I,\Phi)$ is 
well defined for $I \in K^\omega_{\tr}$.}.
\end{enumerate}
\mn
3) If $T \subseteq T_{1}$ are complete first order theories and $\kappa=
\cf(\kappa)<\kappa(T)$ \then
\mn
\begin{enumerate}
\item[$(a)$]   there is a sequence of first order formulas
$\varphi_i(\bar x,\bar y_i)$ (for $i<\kappa$) witnessing $\kappa<\kappa(T)$
i.e. there are a model $M$ of $T$ and sequences $\bar{a}_\eta$ for $\eta
\in {}^{\kappa\le}\lambda$ such that for $\eta\in {}^{\kappa}\lambda,
\nu \in {}^i\lambda,i<\kappa,\alpha<\lambda$ we have $M\models
\varphi_i [\bar{a}_\eta, \bar{a}_{\nu \char 94 \langle \alpha\rangle}]$
\Iff \, $\alpha=\eta(i)$
\sn
\item[$(b)$] for any $\langle \varphi_i(\bar{x},\bar{y}):i<\kappa\rangle$
as in (a) there is a nice template $\Phi$ proper for 
$K^\kappa_{\tr}$ such that for any $\lambda$:
\sn
\begin{enumerate}
\item[$(\alpha)$]  if $\eta\in {}^\kappa \lambda$, $\nu\in {}^i\lambda,i<
\kappa,\alpha<\lambda$ then
\[
\EM(^{\kappa\ge}
\lambda,\Phi)\models\varphi_i[\bar{a}_\eta,\bar{a}_{\nu \char 94
\langle\alpha\rangle}] \text{ iff } \alpha=\eta(i);
\]
\sn
\item[$(\beta)$]  $\EM(I,\Phi)\models T_1$,
\sn
\item[$(\gamma)$]  $\Phi$ is $\bbL_{\omega,\omega}(\tau_\Phi)$-nice,
\sn
\item[$(\delta)$]  $|\tau_\Phi|=|T_1|+\aleph_0$.
\end{enumerate}
\end{enumerate}
\end{lemma}

\begin{PROOF}{\ref{1.8}}
See \cite[Ch.VII,\S3]{Sh:a}, but here we can consider the conclusion
as the  definition of unstable or unsuperstable and of $\kappa<\kappa(T)$,
respectively.
\end{PROOF}

\begin{remark}
On $K^\omega_{\tr}$ for $\bbL_{\lambda^+,\omega}$ we
need the Ramsey property defined below, see \ref{1.12} (and \ref{1.13}+
\ref{1.14}).
\end{remark}

\noindent
In \cite[Ch.VIII,\S2]{Sh:a} we actually proved:
\begin{theorem}
\label{1.9}
1) If $\lambda>|\tau_\Phi|$, and $\Phi,\tau_\Phi,\langle\varphi_n:n<
\omega\rangle$ are as in Lemma \ref{1.8}(2) (and $\Phi$ is almost
$\bbL_{\omega,\omega}$-nice) \then \,: we can find $I_\alpha
\subseteq {}^{\omega\ge}\lambda$ (for $\alpha<2^\lambda$), 
$|I_\alpha|=\lambda$ such that for
$\alpha \ne \beta$ there is no one-to-one function from 
$\EM(I_\alpha,\Phi)$ onto $\EM(I_\beta,\Phi)$ preserving the 
$\pm\varphi_n$ for $n<\omega$.

\noindent
2)  If $\lambda$ is regular, also for $\alpha \ne\beta$ there is no
one-to-one function from $\EM(I_\alpha,\Phi)$ into $\EM(I_\beta,\Phi)$
preserving the $\pm\varphi_n$ for $n<\omega$.

\noindent
3) The $\varphi_n$'s do not need to be first order, just their vocabularies
should be $\subseteq\tau_\Phi$. But instead ``$\Phi$ is almost
$\bbL_{\omega,\omega}(\tau_\Phi)$-nice" we need just ``$\Phi$ is almost
$\{\pm\varphi_n(\ldots,\sigma_\ell(\bar
x_\ell),\ldots)_{\ell<\ell(n)}:
n<\omega,\sigma_\ell \text{ terms of } \tau_\Phi\}$-nice" 
and we should still demand (as in all this section)
\mn
\begin{enumerate}
\item[$(*)$]  the $\bar a_\eta$ are finite (and we are assuming that the
functions are finitary).
\end{enumerate}
\mn
4) So if as in Lemma \ref{1.8}, $\varphi_n\in {\mathscr L}(\tau)$
\then \, $\{M_\alpha \rest \tau:\alpha<2^\lambda\}$ are $2^\lambda$
non-isomorphic models of $T$ of cardinality $\lambda$.
\end{theorem}

\begin{PROOF}{\ref{1.9}}
This is proved in \cite[\S2 of Ch.VIII]{Sh:a} (though it is not
explicitly claimed, it was used elsewhere and there is no need to change the
proofs). Also we shall later (in \cite[3.1]{Sh:331} we prove better theorems,
mainly getting \ref{1.9}(2) also for singular $\lambda$.
\end{PROOF}

\begin{remark}
\label{1.9A}
1) Applying \ref{1.9}, we usually look at the $\tau$-reducts of the
models $\EM^1(I,\Phi)$ as the objects we are interested in, where the
$\varphi_n$'s are in the vocabulary $\tau$. E.g., for $T\subseteq T_1$ first
order, $T$ unsuperstable, we use $\varphi_n \in \bbL(T)$.

\noindent
2)  The case $\lambda=|\tau_\Phi|$ is harder. In
\cite[Ch.VIII,\S2,\S3]{Sh:a}, the existence of many models in $\lambda$ is
proved for $T$ unstable, $\lambda=|\tau_\Phi|+\aleph_1$ and there 
(in some cases) ``$T_1,T$ first order" is used.
\end{remark}
\bigskip

\centerline {$* \qquad * \qquad *$}
\bigskip

How do we find templates $\Phi$ as required in \ref{1.8} and parallel
situations?

Quite often in model theory, partition theorems (from finite or infinite
combinatorics) together with a compactness argument (or a substitute)
are used to build models. Here we phrase this generally. Note that the size
of the vocabulary ($\mu$ in the ``$(\mu,\lambda)$-large")) is a variant of
the number of colours, whereas $\lambda$ is usually $\mu$; it becomes larger
if our logic is complicated.

\begin{definition}
\label{1.10}
Fix a class $K$ (of index models) and a logic (or logic fragment)
$\cL$.

\noindent
1)  An index model $I\in K$ is called $(\mu,\lambda,\chi)$-Ramsey
for $\cL i$ \when \,:
\mn
\begin{enumerate}
\item[$(a)$]  the cardinality of $I$ is $\le \chi$ and every $\qf$ 
(= quantifier free) type $p$ (in $\tau(K)$) which is realized in some
$J \in K$ is realized in $I$,
\sn
\item[$(b)$]  for every vocabulary $\tau_1$ of cardinality $\le \mu$, 
a $\tau_1$-model $M_1$ and an indexed set $\langle \bar b_t:t \in I
\rangle$ of finite sequences from $|M_1|$ with $\ell g(\bar{b}_t)$ 
determined by the quantifier free type which $t$ realizes in $I$ 
\underline{there is} a template $\Phi$, which is proper for $K$, with 
$|\tau_\Phi|\le \lambda$ such that ($\tau_1 \subseteq \tau_\Phi$ and):
\sn
\begin{enumerate}
\item[$(*)$]  for any $\tau(K)$-quantifier free type $p,I_1 \in K$
  and $s_0,\ldots,s_{n-1} \in I_1$ for which $\langle s_0,\ldots,
s_{n-1} \rangle$ realizes $p$ in $I_1$ and for any formula
\[
\varphi = \varphi(x_0,\ldots,x_{m-1}) \in \cL(\tau_1)
\]
and $\tau_1$-terms $\sigma_\ell(\bar y_0,\ldots,\bar y_{n-1})$ for $\ell=0,
\ldots,m-1$ we have
\sn
\item[$(**)$]  if for every $t_0,\ldots,t_{n-1} \in I$ such that $\langle
t_0, \ldots,t_{n-1}\rangle$ realizes $p$ in $I$ we have
$M_1 \models \varphi[\sigma_0(\bar b_{t_0},\ldots,\bar{b}_{t_{n-1}}),
\sigma_1(\bar b_{t_0},\ldots,\bar b_{t_{n-1}}),\ldots,
\sigma_{m-1} (\bar b_{t_0},\ldots,\bar b_{t_{n-1}})]$ \then \,
$\EM(I_1,\Phi)\models\varphi[\sigma_0(\bar a_{s_0},\ldots,
\bar a_{s_{n-1}}),\sigma_1(\bar a_{s_0},\ldots,\bar a_{s_{n-1}}),
\ldots,\sigma_{m-1}(\bar a_{s_0},\ldots,\bar a_{s_{n-1}})]$.
\end{enumerate}
\end{enumerate}
\mn
2) The class $K$ of index models is called explicitly
$(\mu,\lambda,\chi)$-Ramsey for $\cL$ \Iff \, some $I\in K$
of cardinality $\le \chi$ is $(\mu,\lambda)$-Ramsey for $\cL$. 
A class $K'\subseteq K$ of index models is called
$(\mu,\lambda,i,\chi)$-Ramsey (inside $K$, which is usually
understood from context), \Iff
\mn
\begin{enumerate}
\item[$(a)$]   every member of $K'$ has cardinality $\le \chi$ and
every quantifier free type $p$ in $\tau(K')$ realized in some $J\in K$ is
realized in some $I\in K'$,
\sn
\item[$(b)$]  for every vocabulary $\tau_1$ of cardinality $\leq\mu$ and
$\tau_1$-models $M_I$ for $I\in K'$, and $\bar b_{I,t} \in {}^{k(I,t)}(M_I)$,
where $k(I,t)< \omega$ depends just on $\tp_{\qf}(\langle t \rangle,
\emptyset,I)$ \underline{there is} a template $\Phi$ proper for $K$
with $|\tau_\Phi|\le \lambda$ such that $\tau^1 \subseteq \tau_\Phi$ 
we have $(*)$ only in $(**)$ we should also say ``every $I\in K'$".
Let ``$(\mu,\chi)$-Ramsey" mean ``$(\mu,\mu,\chi)$-Ramsey".
Let ``$\mu$-Ramsey" mean ``$(\mu,\chi)$-Ramsey for some $\chi$".
\end{enumerate}
\mn
3) In all parts of \ref{1.10}, \ref{13Anew}, \ref{1.10A}, if $\cL$
is first order logic, we may omit it.

\noindent
4) For $f:\Card \longrightarrow \Card,K$ is $f$-Ramsey \Iff \, it is
$(\mu,f(\mu))$-Ramsey for $\cL$ for every (infinite) cardinal $\mu$.
We say $K$ is Ramsey for $\cL$ if it is $(\mu,\mu)$-Ramsey for $\cL$ 
for every $\mu$.

\noindent
5)  We say $K$ is $*$-Ramsey for $\cL$ if it is $f$-Ramsey for
$\cL$ for some $f:\Card \longrightarrow \Card$.
\end{definition}

\begin{definition}
\label{13Anew}
Let $K$ be a class of (index) models and ${\mathscr L}$ a logic.

\noindent
1)  We say $I\in K$ is (almost) $\cL$-nicely
$(\mu,\lambda,\chi)$-Ramsey for $K$ \Iff \, \ref{1.10}(1) holds, and 
$\Phi$ is (almost) $\cL$-nice.  
Similarly replacing $I$ by a set $K'\subseteq K$.

\noindent
2) The class $K$ is called explicitly (almost) $\cL$-nice 
$(\mu,\lambda,\chi)$-Ramsey \Iff \, some $I\in K$ is (almost)
$\cL$-nicely $(\mu,\lambda,\chi)$-Ramsey. 

\noindent
3) For $f:\Card\longrightarrow\Card$, we say $K$ is (almost)
$\cL$-nicely $f$-Ramsey \Iff \, for every $\mu$ we have: $K$ is (almost)
$\cL$-nicely $(\mu,f(\mu))$-Ramsey for every (infinite) cardinal $\mu$. We
omit $f$ for the identity function.

\noindent
4)  We say $K$ is (almost) $\cL$-nicely $*$-Ramsey \Iff \,
for some $f$, it is (almost) $\cL$-nicely $f$-Ramsey.
\end{definition}

\begin{definition}
\label{1.10A}
In \ref{1.10}, \ref{13Anew} we add ``strongly'' if we strengthen
\ref{1.10}(1) by asking in $(*)$ in addition that for any
$\tau(K)$-quantifier free type $p$ and $s_0,\ldots,s_{n-1} \in I_1$ 
such that $\langle s_0,\ldots,s_{n-1}\rangle$ realizes $p$ in 
$I_1$) we can find some $t_0,\ldots,t_{n-1}$ suitable for all
$\varphi,\sigma_0,\ldots$ simultaneously (this helps for omitting types).
\end{definition}

\begin{theorem}
\label{1.11}
1)  For $\bbL_{\omega,\omega}$, the class of linear orders is
nicely Ramsey, moreover every infinite order is $(\mu,\lambda)$-Ramsey for
any $\mu \le \lambda$.

\noindent
2)  For $\bbL_{\omega_1,\omega}$ the class of linear orders is nicely
$*$-Ramsey. In fact nicely $f$-Ramsey for the functions 
$f(\mu)=\beth_{(2^\mu)^+}$.

\noindent
3)  For any fragment of $\bbL_{\lambda^+,\omega}$ or of
$\Delta(\bbL_{\lambda^+,\omega})$ (see, e.g. \cite{Mw85})
of cardinality $\lambda$, the class of linear orders is nicely
$f$-Ramsey when $f(\mu)=\beth_{(2^\mu)^+}$, even strongly;
moreover is strongly nicely $f$-Ramsey. 

\noindent
4)  $K^\omega_{\tr}$ (and even $K^\kappa_{\tr})$ is Ramsey for 
$\bbL_{\omega,\omega}$.
For definitions of $K^\omega_{\tr}$ see \ref{1.7} above.

\noindent
5) The class $K_{\org}$ of linear ordered graphs is explicitly nicely
Ramsey. The class $K_{\oor,n}$ of linear orders expanded by an $n$-place
relation is explicitly nicely Ramsey.
\end{theorem}

\begin{PROOF}{\ref{1.11}}
1)  This is the content of the Ehrenfeucht-Mostowski proof that
{E.M.} models exist and it use the finitary Ramsey theorem as used in the
proof of \ref{1.8}(1). see \cite[Ch.VII]{Sh:c}.

\noindent
2) By repeating the proof of Morley's omitting type theorem which use 
the Erd\" os-Rado theorem, see \cite[Ch.VII,\S5]{Sh:c};
the to uncountably vocabulary (and many types)
is a generalization  noted by Chang.

\noindent
3) Like \ref{1.11}(2); see \cite[Theorem 2.5]{Sh:16}, and
more in \cite{GrSh:222}, \cite{GrSh:259}.

\noindent
4) By \cite[Ch.VII,\S3]{Sh:c} (we use the compactness of
$\bbL_{\omega,\omega}$ and partition properties of trees).

\noindent
5) By the Nessetril-Rodl theorem (see e.g. \cite{GrRoSp}).
\end{PROOF}

\noindent
By Grossberg-Shelah \cite{GrSh:238} (improving \cite[Ch.VII]{Sh:a}, where
compactness of the logic $\cL$ was used, but no large cardinals):
\begin{theorem}
\label{1.12}
$K^\omega_{\tr}$ is the nicely $*$-Ramsey for $\bbL_{\lambda^+,\omega}$
\Iff \, for example there are arbitrarily large measurable 
cardinals (in fact, large enough cardinals consistent with 
the axiom $\bold V = \bold L$ suffice).
\end{theorem}

\noindent
We shall not repeat the proof.
\begin{lemma}
\label{1.13}
Suppose $K_1,K_2,K_3$ are classes of models, $\Phi$ is proper template for
$(K_1,K_2),\Psi$ proper template for $(K_2,K_3)$ \then \, there is
a unique template $\Theta$ that is proper for $(K_1,K_3)$ and for $I\in K_1$

\[
\EM(I,\Theta) =\EM(\EM(I,\Phi),\Psi)).
\]

\mn
We write ${\Theta}$ as $\Psi \circ \Phi$.
\end{lemma}

\begin{PROOF}{\ref{1.13}}
Straightforward.
\end{PROOF}

\begin{lemma}
\label{1.14}
Suppose $K$ is a class of index models, $\tau=\tau(K)$ and
\mn
\begin{enumerate}
\item[$(*)$]  there is a template $\Psi$ proper for $K$ such that
$|\tau_\Psi| = |\tau_K| + {\aleph_0} $ and for $I \in
K:\EM_\tau(I,\Psi) \in K$ and $J =: \EM_\tau(I,\Psi)$ is strongly
$(\aleph_0,\qf)$-homogeneous over $I$, i.e., if $\bar t=\langle
t_1,\ldots,t_n\rangle,\bar s=\langle s_1,\ldots,s_n\rangle$ realize
the same quantifier free type in $I$, then some automorphism of $J$
takes $\bar a_{\bar t}$ to $\bar a_{\bar s}$.
\end{enumerate}
\mn
We conclude that: if $K$ is $(\mu,\lambda,\chi)$-Ramsey for $\cL$ 
and $|\tau_\Psi|\le \mu$ \then \, $K$ is almost $\cL$-nicely
$(\mu,\lambda,\chi)$-Ramsey for $\cL$.
\end{lemma}

\begin{PROOF}{\ref{1.14}}
Just chase the definitions.
\end{PROOF}

\begin{remark}
\label{1.14A}
1)  E.g. for $\cL \subseteq \bbL_{\omega_1,\omega}$ we get in
\ref{1.14} even $\cL$-nice.

\noindent
2) The assumption $(*)$ of \ref{1.14}(1) holds for $K_{\oor},
K^\omega_{\tr},K^\kappa_{\tr}$ (as well as the other $K$'s from
\cite{Sh:331}).
\end{remark}

\begin{conclusion}
\label{1.15}
Assume that
\mn
\begin{enumerate}
\item[$(a)$]  $K_{\oor}$ is $(\mu,\lambda)$-Ramsey for $\cL$,
\sn
\item[$(b)$]  $T$ is an $\cL$-theory (in the vocabulary $\tau(T)),
|\tau(T)|\leq\mu$,
\sn
\item[$(c)$]  $\varphi_\ell(\bar R_\ell,\bar x,\bar y) \in \cL(\tau(T)\cup
\{\bar R_\ell\})$ for $\ell=1,2$ (and $\bar R_\ell$ is disjoint from
$\tau(T)$ and from $\bar R_{3-\ell}$), and $T\cup\{\varphi_1(\bar{R}_1, \bar
x,\bar y),\varphi_2(\bar{R}_2,\bar x,\bar y)\}$ has no model,
\sn
\item[$(d)$]  for every $I\in K_{\oor}$ there is a model $M_I$ of $T$, and
$\bar{a}_t \in {}^{\omega>}M$ for $t\in I$ such that:
\[
t<s \Rightarrow  M \models(\exists\bar R_1)\varphi_1(\bar R_1,
\bar a_t,\bar a_s)
\]

\mn
and

\[
s<t \Rightarrow M \models(\exists\bar R_2)\varphi_2(\bar R_2,\bar a_s,\bar
a_t).
\]
\end{enumerate}
\mn
\Then \, for $\lambda \ge \mu+ \aleph_1,\isoI(\lambda,T)=2^\lambda$.
\end{conclusion}

\begin{PROOF}{\ref{1.15}}
Obvious by now (mainly \ref{1.11}(3) and \ref{3.10}(3) below).
\end{PROOF}

\begin{conclusion}
\label{1.16}
The parallel of \ref{1.15} for $K^\omega_{\tr}$ instead $K_{\oor}$ holds if
$\lambda>\mu$.
\end{conclusion}

\begin{PROOF}{\ref{1.16}}
By \ref{1.9} (or use \cite{Sh:331}).
\end{PROOF}
\bigskip

\centerline {$* \qquad * \qquad *$}
\bigskip

\begin{discussion}
\label{1.22new}
We return to the general Ramsey properties for other classes (not just
linear orders and trees). For compact logics,  finitary generalization of
Ramsey theorem suffices. More generally, certainly it is nice to have them for
$\cL = \bbL_{\lambda^+,\omega}$, and even $\Delta(\bbL_{\lambda^+,\omega})$, so
we need a partition theorem generalizing Erd\" os-Rado theorem, i.e.,
the case with infinitely many colours. We may for example look at ordered
graph as index models, quite natural one. It consistently holds (\cite{Sh:289})
though unfortunately it does not necessarily hold 
(Hajnal-Komjath \cite{HaKo97}).
However, our main point is that this is enough when the
consistency is by forcing with e.g. complete enough forcing notion. So
the consistency result in \cite{Sh:289} yields a ``real", ZFC theorem
here. The following is an abstract version of the omitting type theorem.
\end{discussion}

\begin{claim}
\label{1.23new}
Assume that
\mn
\begin{enumerate}
\item[$(a)$]  $K$ is a definition of a class of models with vocabulary $\tau$
(the ``index models"); where $\tau$ and the parameters in the definition
belongs to $\cH(\chi^+)$, 
\sn
\item[$(b)$]  $\cL$ is a definition of a logic or logic fragment, the
parameters of the definition belong to $\cH(\chi^+)$ and $\lambda \ge
\chi$,
\sn
\item[$(c)$]  in the definition of ``$\Phi$ is (almost) $\cL$-nice" for
$\Phi$ proper for $K$ with $|\tau_\Phi|<\chi$ (see \ref{1.6}(3), (4); so
without loss of generality $\Phi \in \cH(\chi)$) it suffices to restrict
ourselves to $I$ of cardinality $<\chi$,
\sn
\item[$(d)$]  $\bbP$ is a forcing notion not adding subsets to $\lambda$, and
preserving clauses (a), (b) and (c) (i.e., the definitions of $K$ and  
$\cL$ have these properties) and no new quantifier free complete
$n$-types are realized in $I \in K$,
\sn
\item[$(e)$]  in $\bold V^{\bbP}$, there is a member $I^*$ of
$K$, which is $(\chi,\lambda)$-Ramsey for $\cL$ (or an almost
$\cL$-nicely $(\chi,\lambda)$-large) [or an $\cL$-nicely 
$(\chi,\lambda)$-Ramsey] or such a subset
$K'$ of $K$. For $I\in K$ let $\bold P^n_I=\{p:p$
is complete quantifier-free $\tau_K$-type realized by some
$\bar t\in {}^{n}I\}$. Let $\bold P_n$ be $\bold P^n_{I^*}$ or 
$\cup \{\bold P^n_I:I \in K\}$ according to the case above; 
if $q \in \bold P^n_I$ as exemplified by $\bar t \in {}^{n} I$ 
let $\proj_\ell(q)$ be the quantifier-free type which $t_\ell$
realizes in $I$
\sn
\item[$(f)$]  $\tau_0 \in \cH(\chi^+)$ is a vocabulary, $q_* \in \bold P_1$
and $\langle \Omega_q:q \in \bold P_n$ for some $n<\omega\rangle$
are such that for every $q \in \bold P_n$ we have:
$\Omega_q \subseteq \{p(\bar x_0,\ldots,\bar x_{n-1}):p$ an 
$\cL(\tau_0)$-type in the variables $\bar x_0,\ldots,\bar x_{n-1}$
where $\bar x^\ell=\langle x_{\ell,i}:i<\alpha_{\proj_\ell(q)}\rangle
\in \cH(\chi^+)$ for some $n<\omega\}$,
and in $\bold V^\bbP$, for every $I \in K$ (in the $\bold V^\bbP$'s
sense) or just $I= I^*$ [or just $I \in K'$, according to the case 
in clause (e)], there is a $\tau_0$-model $M_I$ and 
$\bar b^I_t \in {}^{\alpha_t}(M_I)$ for $t \in I$ such that:
\sn
\begin{enumerate}
\item[$(\alpha)$]  $\alpha_t = \alpha_q$ if $q$ is the
quantifier free-$\tau_0$-$1$-type which $t$ realizes in $I$,
\sn
\item[$(\beta)$]  for no $t_0,\ldots,t_{n-1}\in I$, does $\langle t_0,\ldots,
t_{n-1}\rangle$ realize in $I$ the complete quantifier free
$\tau_\kappa-n$-type $q$ and $p=p(\bar x_0,\ldots,\bar x_{n-1}) \in 
\Omega_q$, does $\bar b^I_{t_0} \char 94 \bar b^I_{t_1} \char 94 \ldots 
\char 94 \bar b^I_{t_{n-1}}$ realizes $p$ and $\alpha_{t_\ell} = 
\ell g(\bar x_\ell)$.
\end{enumerate}
\end{enumerate}
\mn
\Then \, we can conclude that there is a $\Phi$ such that:
\mn
\begin{enumerate}
\item[$\aleph$]  $\Phi$ is an (almost) $\cL$-nice
template $\Phi$, proper for $K$, 
\sn
\item[$\beth$]  $\Phi \in \cH(\lambda^+)$ hence also $\tau_\Phi \in 
\cH(\lambda^+)$
\sn
\item[$\gimel$]  if $M = \EM(I,\Phi)$, and $t_0,\ldots,t_{n-1}\in I$,
  and $\bar t=\langle t_0,\ldots,t_{n-1}\rangle$ realizes the 
complete quantifier free $\tau_\kappa-n$-type $q$ \then \,  
$\bar a_{\bar t}$ does not realize in $M$ any $p\in \Omega_q$.
\end{enumerate}
\end{claim}

\begin{PROOF}{\ref{1.25}}
Straightforward.
\end{PROOF}

\begin{claim}
\label{1.24new}
Assume that
\mn
\begin{enumerate}
\item[$(a)$]  $K$ is a class of (index) models,
\sn
\item[$(b)$]  $\kappa$ is a cardinal, for $\alpha<(2^\kappa)^+$ the structure
$I_\alpha\in K$ realizes all quantifier free $\tau_K$-types 
(in $<\omega$ variables) realized in some $I \in K$, and their 
number is $\le \kappa$,
\sn
\item[$(c)$]  if $n < \omega,\alpha<\beta< (2^\kappa)^+,N$ 
is a model, $\tau(N)\le \kappa,\alpha^*_r < \kappa^+$ for a complete
quantifier free $\tau_K-1$-type $r$ realized in $I_\beta,\bar b_r
\in {}^{\alpha^*_r}N$, \then \, we can find $I'_\alpha \subseteq
I_\beta$ isomorphic to $I_\alpha$ such that
\sn
\begin{enumerate}
\item[$(*)$]  if $\bar t,\bar s \in {}^m (I'_\alpha),m \le n$ and they
realize the same quantifier free type in $I'_\alpha$ then $\bar b_{\bar t}=
\langle \bar b_{t_\ell}:\ell<m\rangle$ and $\bar b_{\bar s}=\langle \bar
b_{\bar s_\ell}: \ell<m\rangle$ realizes the same quantifiers free type in
$N$,
\end{enumerate}
\sn
\item[$(d)$]  $\tau$ is a vocabulary, $|\tau| \le \kappa,\psi \in 
\bbL_{\kappa^+,\omega}(\tau)$ and $\alpha^*_p<\kappa^+$ for $p$ a complete
quantifier free $\tau_K-1$-type realized in every $I_\alpha,\cL \subseteq
\bbL_{\kappa^+,\omega}(\tau)$ is a fragment of cardinality $\kappa$ to
which $\psi$ belongs,
\sn
\item[$(e)$]   for every $\alpha<(2^\kappa)^+$, there is a model $N_\alpha$ of
$\psi$ with $\bar b^\alpha_t \in {}^{\alpha^*_t}(N_\alpha)$ for $t\in
I_\alpha$, where $\alpha^*_t=\alpha^*_{\tp_\qf(t,\emptyset,I_\alpha)}$.
\end{enumerate}
\mn
\Then \, there is a $\cL$-nice template $\Phi$, such that:
\mn
\begin{enumerate}
\item[$\otimes$] for $I\in K,m < \omega$ and $\bar t\in {}^m I$ we have: the
$\cL$-type which is $\bar a_{\bar t}$-realized in $\EM(I,\Phi)$ is
realized in some $N_\alpha$ by some $\bar b_{\bar s}$, where $\tp_\qf(\bar s,
\emptyset,I_\alpha)=\tp_\qf(\bar t,\emptyset,I)$.

In other words, $\{I_\alpha:\alpha< (2^\kappa)^+\}$ is
$\kappa$-Ramsey for $\cL$.
\end{enumerate}
\end{claim}

\begin{PROOF}{\ref{1.24new}}
We can expand $N_\alpha$ by giving names to all formulas in $\cL$ 
and adding Skolem functions (to all first order formulas in the new
vocabulary), so we have a $\tau^+$-model $N^+_\alpha,
\tau^+ \supseteq \tau = \tau(\psi),|\tau^+| \le \kappa$,
correspondingly we extend $\cL$ to a fragment $\cL^+$ of
$\bbL_{\kappa^+,\omega}(\tau^+)$ of cardinality $\kappa$.

By induction on $n<\omega$ we choose $A_n$, $f_n$, $\langle I^n_\alpha:
\alpha\in A_n\rangle$ such that:
\mn
\begin{enumerate}
\item[$(i)$]  $A_n$ is an unbounded subset of $(2^\kappa)^+$,
\sn
\item[$(ii)$]  $f_n$ is an increasing function from $(2^\kappa)^+$ onto $A_n$
such that $\alpha<f_n(\alpha)$,
\sn
\item[$(iii)$]  $I^n_\alpha$ is a submodel of $I_\alpha$ isomorphic to
$I_{f_n^{-1}(\alpha)}$,
\sn
\item[$(iv)$]  if $n>m>0$, $\alpha_1,\alpha_2<(2^\kappa)^+$, $\bar t^1\in {}^m
(I^n_{f_n(\alpha_1)})$, $\bar t^2\in {}^m(I^m_{f(\alpha_2)})$, $\tp_\qf(\bar
t^1,\emptyset, I_{f_n(\alpha_1)})=\tp_\qf(\bar t^2,\emptyset,I_{f_n(
\alpha_2)})$, then the quantifier free type of $\bar b_{\bar t^1}$ in
$N_{f_n(\alpha_1)}$ is equal to the quantifier free type of $\bar b_{\bar
t^2}$ in $N_{f_n(\alpha_2)}$,
\sn
\item[$(v)$]  $A_{n+1}\subseteq A_n$ and $\alpha\in A_{n+1} Rightarrow
I^{n+1}_\alpha\subseteq I^{n+1}_\alpha$.
\end{enumerate}
\mn
For $n=0$ let $A_0=(2^\kappa)^+$ and $I^0_\alpha=I_\alpha$.

For $n+1$, for each $\alpha$ we apply assumption (c) to $N_{f_n(\alpha+n+
1)}$, $I^n_{f_n(\alpha+n+1)}$, $\langle\bar b^\alpha_t:t\in I^n_{f_n(\alpha+
n+1)}\rangle$, getting $J^n_{f_n(\alpha+n+1)}$. We define an equivalence
relation $E_n$ on $(2^\kappa)^+: \alpha\; E_n\; \beta$ if and only if
$\tp(\bar b^{f_n(\alpha+n+1)}_{\bar s},\emptyset,N_{f_n(\alpha+n+1)})=
\tp(\bar b^{f_n(\alpha+n+1)}_{\bar t},\emptyset,N_{f_n(\beta+n+1)})$,
whenever $m<\omega,\bar s\in {}^m(J^n_{f_n(\alpha+n+1)}),\bar t 
\in {}^m(J^n_{f_n(\beta+n+1)})$ and $\tp_\qf(\bar s,\emptyset,
I_{f_n(\alpha+n+1)})=\tp_\qf(\bar t,\emptyset,I_{f_n(\beta+n+1)})$.

Clearly $E_n$ has $\leq 2^\kappa$ equivalence classes, so some equivalence
class $B$ is unbounded in $(2^\kappa)^+$. Let

\[
A_{n+1}=\{ f_n(\alpha+n+1):\alpha\in B\},\quad f_{n+1}(\alpha)=f_n(\min(B
\setminus \alpha)+n+1),
\]

\mn
and $I^{n+1}_{f_n(\alpha+n+1)}=J^n_{f_n(\alpha+n+1)}$ for $\alpha\in B$.

Having completed the induction, clearly we have gotten $\Phi$, as the
limit. 
\end{PROOF}

\begin{conclusion}
\label{1.25new}
Assume that
\mn
\begin{enumerate}
\item[$(a)$]  $\cL$ a fragment of $\bbL_{\kappa^+,\omega}$, $T$ is theory in
$\cL(\tau)$, and $\theta \ge \kappa+|T|+|\tau| + |\cL|$,
\sn
\item[$(b)$]  $\varphi_\alpha =
 \varphi_\alpha(x_0,\ldots,x_{k_\alpha-1}) \in \cL(\tau)$ for 
$\alpha<\alpha(*)$ (where $\alpha(*)<\kappa^+$ may be finite),
\sn
\item[$(c)$]  for some $\mu>\theta$, in any forcing extension of
  $\bold V$ by a $\mu$-complete forcing notion the following holds 
for any $\lambda$:
\sn
\begin{enumerate}
\item[${{}}$]  if $R_\alpha$ is a subset of $[\lambda]^{k_\alpha}$ for $\alpha<
\alpha(*)$ then for some model $M$ of $T$ and $a_\alpha\in M$ for $\alpha<
\lambda$ we have: if $\alpha<\alpha(*),\gamma_0 < \ldots <
\gamma_{k_\alpha-1}<\lambda$, then $M \models 
\varphi_\alpha[a_{\gamma_0},\ldots,a_{\gamma_{k_\alpha-1}}] 
\Leftrightarrow \{\gamma_0,\ldots,\gamma_{k_\alpha-1}\} \in R_\alpha$
\end{enumerate}
\sn
\item[$(d)$]  Let $K$ be the 
class of $(I,<,R_0,\ldots,R_\alpha,\ldots)_{\alpha <\alpha(*)},
(I,<)$ linear order, $R_\alpha$ a symmetric irreflexive
$k_\alpha$-place relation on $I$.
\end{enumerate}
\mn
\Then \, we can find a complete $T_1 \supseteq T$ with Skolem
functions, and a template $\Psi$ proper for $K$ and nice, such that:
\mn
\begin{enumerate}
\item[$(\alpha)$]  $\tau\subseteq \tau_\Psi$ (even $\tau_\Psi$ extends
$\tau$), and $|\tau_\Psi| \le \theta$ and $|T_1| \le \theta$,
\sn
\item[$(\beta)$]  $\Psi$ is nice for $\cL$ and $\EM^1(I,\Psi)\vDash T_1$
for $I\in K$,
\sn
\item[$(\gamma)$]  if $\alpha<\alpha(*)$, and $I \models 
t_0 < \ldots <t_{k_\alpha-1}$ then:
\[
\EM(I,\Psi)\models\varphi_\alpha[a_{t_0},\ldots,a_{t_{k_\alpha-1}}]
\text{ iff } I \models R_\alpha(t_0,\ldots,t_{k_\alpha-1}).
\]
\end{enumerate}
\end{conclusion}

\begin{PROOF}{\ref{1.25new}}
We would like to apply \ref{1.24new}, e.g., with $I_\alpha\in K$
being of cardinality $\beth_{\omega \alpha+1}(\theta)$, and being 
$\beth_{\omega \alpha}(\theta)^+$-saturated for quantifier free 
types in the natural sense (such $N_\alpha$ exists by the 
compactness theorem). However why does assumption
(c) of \ref{1.24new} hold? By \cite{Sh:289} there is a $\theta^+$-complete
forcing notion $\bbP$ such that in $\bold V^\bbP$ this will hold;
it would not make a real difference if we replace 
$\beth_{\omega\alpha +1}(\theta)$ by other suitable cardinal. 
But by \ref{1.23new} this suffices (as our assumptions are absolute enough).
\end{PROOF}

\begin{remark}
\label{1.25d}
For first order $T$, this help  in Laskowski-Shelah \cite{LwSh:687}.
\end{remark}

\begin{conclusion}
\label{1.26new}
If $T$ is first order countable with the $OTOP$ (see \cite[Ch.XII]{Sh:c},
the omitting type order property)
\then \, for some sequence $\bar\varphi=\langle\varphi_i(\bar x,
\bar y,\bar z): i<i(*)\rangle$ of first order formulas in
$\bbL_{\omega,\omega} (\tau_T)$ and template $\Phi$ proper for
linear orders we have:
\mn
\begin{enumerate}
\item[$(\alpha)$]  $\tau_T\subseteq\tau_\Phi$, $|\tau_\Phi|=|\tau_T|+
\aleph_0$,
\sn
\item[$(\beta)$]  $\EM_{\tau(T)}(I,\Phi)\models T$ for $I \in
  K_{\org}$,
\sn
\item[$(\gamma)$]  if $I\in K_{\org}$ and $s,t \in I$ then

\[
\EM_{\tau(T)}(I,\Phi) \models (\exists \bar x)
\bigwedge\limits_{i<i(*)}\varphi_i(\bar x,\bar a_s,\bar a_t) \text{ iff }
I \models s R t.
\]
\end{enumerate}
\end{conclusion}

\begin{PROOF}{\ref{1.26new}}
Similarly: $\OTOP$ is defined in \cite[Ch.XII,4.1,p.608]{Sh:c}, in
a way giving clause (e) of \ref{1.24new} above directly, but we
need to know that it is absolute (or just preserved by 
$\lambda$-complete forcing), which holds by \cite[Ch.XII,4.3,p.609]{Sh:c}.
\end{PROOF}

\begin{conclusion}
\label{1.28}
Claim \ref{1.24new} applies to the class of trees with $\omega$ levels.
\end{conclusion}

\begin{PROOF}{\ref{1.28}}
By the proof in \cite[Ch.VII,\S3]{Sh:c}, i.e., looking at what we use
and applying the Erd\" os-Rado theorem.
\end{PROOF}
\newpage

\section {Models Represented in Free Algebras and Applications}

This section presents a framework, which tries to separate the model
theory and combinatorics of \cite[Ch.VIII]{Sh:c} and improve it. We shall prove
the combinatorics in \cite{Sh:309} and \cite{Sh:331};
here we try to show how to apply it.
More applications and combinatorics are in \cite{Sh:511}.

\begin{discussion}
\label{2.1}
We sometimes need $\tau_\Phi$ with function symbols with infinitely many
places and deal with logics $\cL$ with formulas with infinitely many
variables. Why?
\end{discussion}

\begin{example}
\label{2.1A}
We would like to build complete Boolean algebras without non-trivial one-to-one
endomorphisms. How do we get completeness? We build a Boolean algebra, $B_0$
and take its completion. Even when $\bold B_0$ satisfies the c.c.c. we need the
term $\bigcup\limits_{n<\omega} x_n$ to represent elements of the Boolean
algebra from the ``generators" $\{\bar a_t:t\in I\}$.
\end{example}

\begin{discussion}
\label{2.1B}
We also sometimes would like to rely on a 
well ordered construction, i.e., on the
universe of $\cM_{\mu,\kappa}$ there is a well ordering which is involved in the
definition of indiscernibility (see \ref{2.2}). This means that we have in
addition an arbitrary well-order relation. E.g., we would like to build many
non-isomorphic $\aleph_1$--saturated models for a stable not superstable
first order theory, with the DOP (dimensional order property, see
\cite[Ch.X]{Sh:c}) so for some $\varphi(\bar x,\bar y)$ 
(not first order), for any cardinal $\lambda$ for some 
model $M$ of $T$, we have a family $\{\bar a_\alpha:\alpha<\lambda\}$
of sequences of length $\le |T|$ in $M$ with
$M\models\varphi[\bar a_\alpha,\bar{ a}_\beta]$ 
\Iff \,  $\alpha<\beta$. The formula
$\varphi$ says: there are $z_\alpha$ ($\alpha<|T|^+$) such that $\bar x
\char 94 \bar y \char 94 \langle z_\alpha:\alpha<|T|^+\rangle$ realizes a type
$p$. So there is a template $\Phi$ proper for $K_{\oor}$ such that for
$I\in K_{\oor}$ and $s$, $t\in I$ we have

\[
\EM_{\tau(T)}(I,\Phi) \models \varphi[\bar a_{s},\bar a_{t}]
\text{ iff } I \models s<t
\]

\mn
($<$ a relevant order), but we need to make them $\aleph_1$-saturated.
Ultrapowers may well destroy the order. The natural thing is to make
$M_I$ $\aleph_1$-constructible over $EM_{\tau(T)}(I,\Phi)$, that is it's
set of elements is $\{b_\alpha:\alpha<\alpha\}, b_\alpha$ realizing over
$\EM_\tau(I,\Phi)\cup\{b_\beta:\beta<\alpha\}$ in $M_I$ a complete
type which is $\aleph_1$-isolated. So not only are
the $\bar a_t$ infinite and the construction involves infinitary functions,
but {\it a priori} the quite arbitrary order of the constructions may play a
role.

With some work we can eliminate the well order of the construction for this
example (using symmetry, the non-forking calculus) but there is no guarantee
generally and certainly it is not convenient, for example see the
constructions in \cite[\S3]{Sh:136}. Moreover, it is better to delete the
requirement that the universe of the model is so well defined.
\end{discussion}

\noindent
This motivates the following definition.
\begin{definition}
\label{2.2}
\mn
\begin{enumerate}
\item[$(a)$]  $\tau(\mu,\kappa)=\tau_{\mu,\kappa}$ is the vocabulary with
function symbols
\[
\{F_{i,j}:i<\mu,\ j<\kappa\},
\]
where $F_{i,j}$ is a $j$-place function symbol and $\kappa$ is $\aleph_0$
or an uncountable regular cardinal
\sn
\item[$(b)$]  $\cM_{\mu,\kappa}(I)$ is the free $\tau$-algebra generated
by $I$ for $\tau = \tau_{\mu,\kappa}$.
\end{enumerate}
\mn
We use the following notation in the remainder of this definition.

Let $f:M \longrightarrow \cM_{\mu,\kappa}(I)$. For $\bar a=\langle a_i:
i<\alpha\rangle\in {}^\alpha M$ let for $i<\alpha$, $f(a_i)=\sigma_i(\bar
t_i)$, where $\bar t_i$ is a sequence of length $<\kappa$ from $I$ and
$\sigma_i$ is a term from $\tau_{\mu,\kappa}$. 

Now if $\alpha<\kappa$ then there is one sequence $\bar t$ of 
members of $I$ of length $<\kappa$ such that

\[
\bigwedge\limits_i\Rang (\bar t_i)\subseteq \Rang(\bar t);
\]

\mn
so we can find terms $\sigma'_i$ satisfying $f(a_i) = \sigma'_i(\bar
t)$, so without loss of generality $\bar t_i=\bar t$, we let 
$\bar\sigma=\langle\sigma_i:i< \alpha\rangle$ and 
$\bar\sigma(\bar t)$ be $\langle\sigma_i(\bar t):i<\alpha \rangle$, 
so $f(\bar a)=\bar\sigma (\bar t)$.
\mn
\begin{enumerate}
\item[$(c)$]   We say that $M$ is $\Delta$-represented in
$\cM_{\mu,\kappa}(I)$ \Iff \, there is a function $f:M
\longrightarrow \cM_{\mu,\kappa}(I)$ such that the $\Delta$-type of 
$\bar a \in {}^{\kappa>}M$ (i.e., $\tp_\Delta(\bar a,\emptyset,M))$
can be calculated from the sequence of terms $\langle
\sigma_i:i < \alpha\rangle$ and $\tp_{\qf}(\langle\bar t_i:i<\alpha\rangle,
\emptyset,I)$ where $f(\bar a)=\langle\sigma_i(\bar t_i):i<\alpha\rangle$ (from
(b), so if $f(\bar a)=\bar\sigma(\bar t)$ from then can be 
calculated $\bar\sigma$ and $\tp_{\qf}(\bar t,\emptyset,I)$). 
We may say ``$M$ is $\Delta$-represented in $\cM_{\mu,\kappa}(I)$ 
by $f$"; similarly below.
\sn
\item[$(d)$]  We say that $M$ is weakly $\Delta$-represented in 
$\cM_{\mu,\kappa}(I)$ \Iff \, for some function $f:M \longrightarrow 
\cM_{\mu,\kappa}(I)$, there is a well-ordering $<$ of the universe 
of $\cM_{\mu,\kappa}(I)$ such that for $\bar a \in {}^\alpha M$ 
the $\Delta$-type of $\bar a$ can be computed
from the information described in (c) and the order $<$ restricted to the
family of subterms of the terms $\langle\sigma_i(\bar t_i):i<\alpha\rangle$.
\end{enumerate}
\mn
[We introduce weak representability to deal with the dependence on the order
of a construction, (cf. \ref{2.1B})].
\mn
\begin{enumerate}
\item[$(e)$]   For $i=1,2$ if $\bar a_i=\langle\sigma^i_j(\bar t^i_j):j<\alpha
\rangle$, $\sigma^1_j=\sigma^2_j$ and
\[
\tp_{\qf}(\langle\bar t^1_j:j<\alpha\rangle,\emptyset,I) = \tp_{\qf}(\langle
\bar t^2_j:j<\alpha\rangle,\emptyset,I)
\]
we write $\bar a^1 \sim \bar a^2 \mod \cM_{\mu,\kappa}(I)$
and may say $\bar a^1,\bar a^2 $ are similar in $\cM_{\mu,\kappa}(J)$.
For the case of weak representability we write $\bar a^1 \sim \bar a^2
\mod (\cM_{\mu,\kappa}(I),<)$ and may say $\bar a^1,\bar a^2$ are
similar in $(\cM_{\mu,\kappa}(J),<)$ \when \, in addition the mapping
\[
\{\langle\sigma (\bar t^1_i),\sigma (\bar t^2_i)\rangle:i<\alpha,\
\sigma \text{ is a subterm of } \sigma^1_i=\sigma^2_i\}
\]
is a $<$-isomorphism (and both sides are linear orders). We write $\bar a^1
\sim_A\bar a^2 \mod \ldots$ if $\bar a^1 \char 94 \bar b \sim \bar a^2
\char 94 \bar b \mod \ldots$ whenever $\bar b \in {}^{\kappa>}A$ where
$A \subseteq \cM_{\mu,\kappa}(I)$. (This latter is especially 
important when we work over a set of parameters). We might, for
instance, insist that $\bar t^1_i$ and $\bar t^1_j$ realize the same 
Dedekind cut over $I_0 \subseteq I$. (So ``$M$ is
$\Delta$-represented in $\cM_{\mu,\kappa}(I)$" means: $f(\bar a^1)$
similar to $f(\bar a^2) \mod \cM_{\mu,\kappa}$ implies $\bar a^1$ and
$\bar a^2$ realize the same $\Delta$-type in $M$.)
\sn
\item[$(f)$]  We say the representation is full \when \,
\[
c_1 \sim c_2 \mod \cM_{\mu,\kappa}(I)) \text{ implies }
[c_1 \in \Rang(f) \Leftrightarrow c_2 \in \Rang(f)].
\]
We say the weak representation is full if we replace $\cM_{\mu,\kappa}
(I)$ by $(\cM_{\mu,\kappa}(I),<)$, where $<$ is a given well ordering from
clause (d).  
\sn
\item[$(g)$]  If $\Delta$ is the family of quantifier free formulas it
  may be omitted.
\sn
\item[$(h)$]  For $f:M \longrightarrow \cM_{\mu,\kappa}(I)$, let $\bar a\sim
\bar b \mod (f,\cM_{\mu,\kappa}(I))$ means
\[
f(\bar a) \sim f(\bar b) \mod \cM_{\mu,\kappa}(I).
\]
Similarly, $\bar a\sim\bar b \mod (f,\cM_{\mu,\kappa}(I),<)$ means
\[
f(\bar a)\sim f(\bar b) \mod (\cM_{\mu,\kappa}(I),<).
\]
\item[$(i)$]  There is no harm in allowing $f$ (in clauses (c),(d)) 
to be multi-valued, but we shall mention explicitly when we allow 
multi-valued functions.
\sn
\item[$(j)$]  We may restrict ourselves to well orderings $<$
of $\cM_{\mu,\kappa}(I)$ which respect subterms; this means that if 
$\sigma_1(\bar t_1)$ is a subterm of $\sigma_2(\bar t_2)$
then $\sigma_1(\bar t_1) \le \sigma_2(\bar t_2)$.
\end{enumerate}
\end{definition}

\noindent
Now we define a very strong negation (when $\varphi$ is ``right") to even
weak representability.
\begin{definition}
\label{2.3}
1) $I$ is strongly $\varphi(\bar x,\bar y)$-unembeddable for $\tau(\mu,
\kappa)$ into $J$ \Iff \, for every $f:I \longrightarrow 
\cM_{\mu,\kappa}(J)$ and well ordering $<$ (of $\cM_{\mu,\kappa}(J)$)
there are sequences $\bar x,\bar y$ of members of $I$ such that $I\models
\varphi [\bar x,\bar y]$ and $\bar x,\bar y$ have ``similar"
(\ref{2.2}(e)) images in $\cM_{\mu,\kappa}(J),<)$.
If we delete the well ordering, we get only ``$I$ is $\varphi
(\bar x,\bar y)$-unembeddable". If $\varphi$ clear from the context we may
omit it.  Note that the formula $\varphi(\bar x,\bar y)$
should be in the vocabulary $\tau_I$; here almost always we have
$\tau_J = \tau_I$ but this is not really necessary. 

\noindent
2) $K$ has the [strong] $(\chi,\lambda,\mu,\kappa)$-bigness property for
$\varphi(\bar x,\bar y)$ \Iff \, there are $I_\alpha\in K_\lambda$ for $\alpha<
\chi$ such that for $\alpha \ne \beta$ we have $I_\alpha$ is [strongly]
$\varphi(\bar x,\bar y)$-unembeddable for $\tau(\mu,\kappa)$ into
$I_\beta$.

\noindent
3) $K$ has the full [strong] $(\chi,\lambda,\mu,\kappa)$-bigness
property for $\varphi(\bar x,\bar y)$ \Iff \, there are $I_\alpha\in
K_\lambda$ for $\alpha<\chi$ such that, for $\alpha<\chi$, $I_\alpha$ is
[strongly] $\varphi(\bar x,\bar y)$-unembeddable for $\tau(\mu,\kappa)$
into $\sum\limits_{\beta<\chi,\beta \ne \alpha} I_\beta$ (where
$\sum\limits_{\beta\in u} I_\beta$, when all the $I_\beta$ are $\tau$-models
for some fixed vocabulary $\tau$, is a $\tau$-model $I$ with universe
$\bigcup\limits_{\beta\in u}|I_\beta|$; if those universes are not pairwise
disjoint we use $\bigcup\limits_{\beta\in u}(\{\beta\}\times(I_\beta))$; for
a predicate $P\in\tau $, $P^I=\bigcup\limits_{\beta\in u}P^{I_B}$, for every
function symbol $F\in \tau$, $F^I$ is the (partial) function
$\bigcup\limits_{\beta\in u}F^{I_\beta}$).

\noindent
4) Saying ``$I$ is [strongly] $\varphi(\bar x,\bar y)$-unembeddable into
$J$ for function $f$ satisfying $\Pr$" means we restrict ourselves (in
\ref{2.3}(1)) to function $f$ from $I$ to $\cM_{\mu,\kappa}(J)$
satisfying $\Pr$.

\noindent
5)  The most popular restriction is ``$f$ finitary on some $P$" which
means that for every $\eta\in P^I$ for some $n<\omega$,
$\tau_{\mu,\kappa}$-term $\sigma$ and $\eta_0,\ldots,\eta_{n-1}\in J$ we
have $f(\eta)=\sigma(\eta_0,\ldots,\eta_{n-1})$. We say $f$ is strongly
finitary if in addition $\sigma$ has only finitely many subterms.

\noindent
6) Clearly (4) induces parallel variants of \ref{2.3}(2), \ref{2.3}(3).
\end{definition}

\begin{remark}
\label{2.3A}
1) This definition is used in proving that the model constructed from $I$ is
not isomorphic to (or not embeddable into) the model constructed from $J$.
For existence see \cite[2.15(2)]{Sh:331} (which we deduce from
\cite[1.7(2)]{Sh:331}).

\noindent
2)  We may in \ref{2.3}(1) and the other variants, add: 
moreover, given $A \subseteq J$ of cardinality $< \kappa$ we demand 
that $\bar x,\bar y$ are similar over $A$.
This does not make a real difference so far. 

\noindent
3) About the connection to $\dot I \dot E(\lambda,T_1,T)$ see \cite{Sh:331}.
\end{remark}

\begin{claim}
\label{2.3B}
If $\Phi$ is proper for $I$ and $\mu=|\tau_\Phi|$ \then \, $\EM(I,\Phi)$ can be
represented in $\cM_{\mu,\aleph_0}$.
\end{claim}

\begin{PROOF}{\ref{2.3B}}
Easy.
\end{PROOF}
\bigskip

\centerline {$* \qquad * \qquad *$}
\bigskip

\begin{discussion}
\label{2.4}
The following example illustrates the application of this method. We first
fix $K^\omega_{\tr}$ (see \ref{1.7}) as the class of index models and fix a
formula $\varphi_{\tr}$ (see \ref{2.4A} below); note that we shall prove
later that for many pairs
$I,J\in K^\omega_{\tr},I$ is $\varphi_{\tr}(\bar x,\bar y)$-unembeddable
in $J$. In \ref{2.5A} below we choose for each $I \in K^\omega_{\tr}$ a
reduced separable Abelian $\dot p$-group $\bbG_I$ which is
representable in $\cM_{\omega,\omega}(I)$. 
In \ref{2.5B} below we show that: [$I$ is
$\varphi_{\tr}$-unembeddable in $J$ implies $\bbG_I \ncong \bbG_J$];
thus the number of reduced separable Abelian $\dot p$-groups
of cardinality $\lambda$ is at least as great as the number 
of trees in $K^\omega_{\tr}$ with cardinality $\lambda$ which
are pairwise $\varphi_{\tr}$-unembeddable. We showed in \cite{Sh:136} that
this number is $2^\lambda$ for regular $\lambda$ and many singulars. But as
said in \ref{1.9} for every uncountable $\lambda$ we get $2^\lambda$ pairwise
non-isomorphic such groups in $\lambda$, using $\bbG_I$ as below.

We may like to strengthen ``$\bbG_I \not\cong \bbG_J$" to
``$\bbG_I$ not embeddable in $\bbG_J$". This depends on the exact notion of
embeddability we use (we shall return to this in \cite[3.22]{Sh:331}).
\end{discussion}

\begin{example}
\label{2.4A}
For the class of $I \in K^\omega_{\tr}$

\begin{equation*}
\begin{array}{clcr}
\varphi_{\tr}(x_0,x_1:y_0,y_1) := &[x_0=y_0] \text{ and }
P_\omega(x_0) \text{ and} \\
  &\bigvee\limits_{n<\omega} [P_n(x_1) \text{ and } P_n(y_1) \text{
    and } P_{n-1}(x_1 \cap y_1)] \text{ and} \\
  &[x_1 \triangleleft x_0 \wedge y_1 \ntriangleleft y_0] \text{ and }
y_1<_{\lx} x_1]
\end{array}
\end{equation*}

\mn
in other words, when for transparency we restrict ourselves to
standard $I \subseteq {}^{\omega \ge}\lambda:x_0 = y_0 \in
{}^\omega\lambda$, and for some $n<\omega$ and $\alpha<\beta<\lambda$ we have

\[
x_1 = (x_0\rest n) \char 94 \langle\alpha\rangle \triangleleft x_0
\]

\mn
and

\[
y_1 = (x_0 \rest n) \char 94 \langle\beta\rangle
\]
\end{example}

\noindent
The connection of the bigness properties from \ref{2.3} to the results 
on $\dot I \dot E(\lambda,T_1,T)$ is done by: 
\begin{claim}
\label{2.4B}
Assume that
\mn
\begin{enumerate}
\item[$(a)$]  $\Phi,\varphi_n$ are as in the conclusion of \ref{1.8}(1), $\mu=
|\tau_\Phi|$,
\sn
\item[$(b)$]  $I,J\in K^\omega_{\tr}$, $I$ is strongly
$\varphi_{\tr}$-unembeddable into $J$ for a $\tau_{\mu,\aleph_0}$,
\sn
\item[$(c)$]  $\tau_0\subseteq \tau_\Phi$ is a vocabulary including that of the
$\varphi_n$'s.
\end{enumerate}
\mn
\Then \, $\EM_{\tau_0}(I,\Phi)$ cannot be elementarily embedded into
$\EM_{\tau_0}(J,\Phi)$.  Moreover, no function from $\EM(I,\Phi)$ into
$\EM(J,\Phi)$ preserves the formulas $\pm \varphi_n$ (for $n<\omega$).
\end{claim}

\begin{PROOF}{\ref{2.4B}}
Straightforward, reread the definitions.
\end{PROOF}

\begin{se}
\label{2.5}
Separable reduced Abelian $\dot p$-groups. 

\noindent
(See more in \cite[\S3]{Sh:331}; as $p$ denote types we use $\dot{p}$
for prime numbers.)
\end{se}

\begin{definition}
\label{2.5A}
1) A separable reduced Abelian $\dot p$-group $\bbG$ is a group
$\bbG$ which satisfies (we use additive notation):
\mn
\begin{enumerate}
\item[$(a)$]  $\bbG$ is commutative (that is ``Abelian''),
\sn
\item[$(b)$]  for every $x \in \bbG$ for some $n$, $x$ has order
$\dot p^n$ (i.e., $\dot p^n x$ is the zero and $n$ is minimal),
\sn
\item[$(c)$]  $\bbG$ has no divisible non-trivial subgroup (=
  reduced),
\sn
\item[$(d)$]  every $x \in \bbG$ belongs to some 1-generated subgroup which is
a direct summand of $\bbG$ (= separable).
\end{enumerate}
\mn
2)  Any such group is a normed space:

\[
\|x\| = \inf\{2^{-n}:(\exists y \in \bbG) \dot p^n y = x\}.
\]

\mn
3) For a tree $I\in K^\omega_{\tr}$ we define the
$\dot p$-group $\bbG_I$ as follows, $\bbG_I$ is generated 
(as an Abelian group) by

\[
\{x_\eta:\eta\in\bigcup\limits_{n<\omega} P^I_n\} \cup \{y^n_\eta:\eta\in
P^I_\omega \text{ and } n<\omega\},
\]

\mn
freely except for the relations:

\[
\dot p^{n+1} x_\eta = 0 \text{ for } \eta \in P^I_n;
\]

\mn
and

\[
\dot p y^{n+1}_\eta-y^n_\eta = x_{\eta \rest n} \text{ and }
\dot p^{n+1} y^n_\eta = 0 \text{ for } \eta \in P^I_\omega.
\]

\mn
4)  It is well known that $\bbG_I$ is a reduced separable Abelian 
$\dot p$-group.  Also note that we have essentially say

\[
y^n_\eta = \sum\{\dot p^{\ell-n} x_{\nu_\ell}:\ell \text{ satisfies }
n \le \ell<\omega,\nu_\ell \in P^I_\ell \text{ and } \nu_\ell
\vartriangleleft \eta\}
\]

\mn
(the infinitary sum may be well defined as $\bbG_I$ is a normed
space).
\end{definition}

\noindent
It is easy to see that
\begin{fact}
\label{2.5B}
$\bbG_I$ is a reduced separable Abelian $\dot p$-group which is
represented in $\cM_{\omega,\omega}(I)$.
\end{fact}

\noindent
We shall prove now
\begin{fact}
\label{2.5C}
If $I$ is $\varphi_{\tr}$-unembeddable into $J$ \then \, $\bbG_I
\not\cong \bbG_J$.
\end{fact}

\begin{proof}
Let $g:\bbG_I \cong \bbG_J \longrightarrow_{h} g$ be an isomorphism 
from $\bbG_I$ onto $\bbG_J$ and 
$h:\bbG_J \longrightarrow \cM_{\omega,\omega}(J)$, where $h$ 
witnesses that $\bbG_J$ is representable in $\cM_{\omega,\omega}(J)$. 

Let $f:I \longrightarrow \bbG_I$ be:

\[
f(\eta) = \begin{cases} 
\sum\limits_{1\le\ell\le\ell g(\eta)} \dot p^{\ell-1}x_{\eta\rest\ell}
\quad &\text{ if } \quad \eta\in\bigcup\limits_{n<\omega}P^I_n,\\
y^1_\eta \quad &\text{ if } \quad \eta\in P^I_\omega.
\end{cases}
\]

\mn
So $(h \circ g \circ f):I \longrightarrow \cM_{\omega,\omega}(J)$. Now we
use the fact that $I$ is $\varphi_{\tr}$-unembeddable into $J$. 

So suppose
\[
I \models \varphi_{\tr}[\eta_0,\nu_0;\eta_1,\nu_1] \text{ and }
h \circ g \circ f(\eta_0,\nu_0) \sim h \circ g \circ f(\eta_1,\nu_1).
\]

\mn
Invoking the definition of $\varphi_{\tr}$: for some $\eta := \eta_0 =
\eta_1 \in P^I_\omega$ and for some $n$,

\[
\nu_1 \vartriangleleft \eta_1,\nu_1 \in P^I_n,\nu_0 \in P^I_n,
\]

\[
\nu_1 \rest (n-1) = \nu_0 \rest (n-1),\nu_0(n-1) < \nu_1(n-1).
\]

\mn
For $i=0,1$ let

\[
z_{\nu_i} = \sum\{\dot p^{\ell-1} x_\nu:\, \nu \vartriangleleft \nu_i,
\nu \in P^I_\ell \text{ and } 1 \le \ell \le n\}.
\]

\mn
Now $\bbG_I \models ``\dot p^n$ divides $(y^1_\eta-z_{\nu_0})$",
hence, as $g$ is an isomorphism, $\bbG_J \models ``\dot p^n$ divides 
$(g(y^1_\eta)-g(z_{\nu_0}))"$, which means $\bbG_J \models 
``\dot p^n$ divides $(g \circ f(\eta)-g \circ f(\nu_0))"$.

Similarly, $\bbG_J \models ``\dot p^n$ does not divide
$(g \circ f(\eta)-g \circ f(\nu_1))"$, but

\[
h \circ g \circ f(\langle\eta_0,\nu_0\rangle) \sim h \circ g \circ f(\langle
\eta_1,\nu_1\rangle) \mod \cM_{\omega,\omega}(J),
\]

\mn
a contradiction, proving \ref{2.5C}. 
\end{proof}
\bigskip

\centerline {$* \qquad * \qquad *$}
\bigskip

\begin{discussion}
\label{2.6}
We still can get considerable amounts of information by the general
theory. When we try to construct many models of $K$ (no one 
embeddable into the others) we need
\mn
\begin{enumerate}
\item[$(*)$]  there are $2^\lambda$ index models $I$ of cardinality $\lambda$
each $\varphi_K(\bar x,\bar y)$-unembeddable into any other.
\end{enumerate}
\mn
But when you intend to construct rigid, indecomposable, etc., you need:
\mn
\begin{enumerate}
\item[$(**)$]  there are $\{I_\alpha\in K:\alpha<\lambda\},I_\alpha,
\varphi_K$-unembeddable into $\sum\limits_{\beta\ne\alpha} I_\beta$ (and
$I_\alpha$ has cardinality $\lambda$).
\end{enumerate}
\end{discussion}

\noindent
Why?
\begin{example}
\label{2.7}
{\rm Constructing Rigid Boolean Algebras.}
(See more, and for more details, in \cite[\S2]{Sh:511}.) For $I \in
K^\omega_\tr$ let $\BA_{\tr}(I)$ be the Boolean Algebra freely generated by
$\{a_\eta:\eta\in I\}$ except the relations

\[
a_\eta \le a_\nu \text{ when } \nu \in P^I_\omega,n<\omega,\eta=\nu
\rest n.
\]

\mn
We shall choose a sequence $\langle \bold B_i,a_j:i \le \lambda,j < 
\lambda \rangle$ such that $\bold B_i$ is a Boolean algebra, 
$\subseteq$-increasing with $i,a_i \in \bold B_i$ and if $i < \lambda$ 
and $a \in \bold B_i$ then $a=a_j$ for some $j \in [i,\lambda)$. 
Start with $\bold B_0 = \BA_\tr(I_0)$, successively for some $a_i \in 
\bold B_i,0 < a_i <1$, take

\[
\bold B_{i+1} = (\bold B_i \rest (1-a_i))+((\bold B_i \rest a_i) *
\BA_{\tr}(I_i)),
\]

\[
\bold B_\lambda = \bigcup\limits_{i<\lambda} \bold B_i = 
\{a_i:i<\lambda\},|I_\alpha|=\lambda.
\]

\mn
(In such situations we say that $\bold B_{i+1}$ is a result of the
$\BA_{\tr}(I_i)$-surgery of $\bold B_i$ at $a_i$ that is, below
$1-a_i$ we add nothing and below $a_i$ we use the free product of
$\bold B_i \rest a_i$ and $\BA_{\rm tr}(I_i)$.)

\noindent
Of course, we choose $\{I_\alpha:\alpha<\lambda\}$ such that $I_\alpha$ is
$\varphi_{\tr}$-unembeddable into $\sum\limits_{\beta\ne\alpha}I_\beta$.
The point is that each $a\in \bold B_\lambda\setminus\{0,1\}$ was ``marked" by
some $I_\alpha$, (the $\alpha$ such that $a_\alpha=a$). Now
$\BA_{\tr}(I_\alpha)$ is embeddable into $\bold B_\lambda \rest a_\alpha$; but
$\bold B_{\lambda} \rest (1-a_\alpha)$ is weakly
$\bbL_{\omega,\omega}$-represented in $\cM_{\omega,\omega}(\sum_{\beta
  \ne \alpha}I_\beta)$.
So for no automorphism $f$ of $\bold B_\lambda$ do we have,
$f(a_\alpha) \le 1-a_\alpha$, which suffices to get
``$\bold B_\lambda$ is rigid"; in fact, it has no one-to-one endomorphism. If
we are trying to get stronger rigidity and/or $\bold B_\lambda
\models$ {c.c.c.}, and/or $\bold B_\lambda$ is complete, we may have to change
$K^\omega_{\tr}$ and/or $\varphi_{\tr}$.
\end{example} 

This illustrates the need for some of the complications in definition
\ref{2.1}. E.g., the weak representation and the uncountable $\kappa$ (for
complete Boolean Algebras).

The definition below (variants of closure under sums) are
satisfied by the cases we shall deal with and enable us to translate
results e.g. from the full
(strong) $(\lambda,\lambda,\mu,\kappa)$-bigness to the (strong)
$(2^\lambda,\lambda,\mu,\kappa)$-bigness.

\noindent
Of course: 
\begin{definition}
\label{2.15}
We say that the class $K$ of $\tau$-structures; 
with $\tau$ a relational vocabulary for transparency, is closed under sums
\when \, for every sequence $\langle I_s: s \in S \rangle$
of members of $K$, pairwise disjoint for simplicity, also $I$ 
belongs to $K$ where $I$ is the $\tau$- structure which is the 
union of $\langle I_s: s \in S \rangle$; that is the set of elements
of $I$ is the union of the sets of elements of $I_s$
for $s \in S$ and $P^I= \cup\{P^{I_s}: s \in S\}$ for every predicate
$P$ from $\tau$.
\end{definition}

\noindent
But in many cases which interest us, this is only almost true, hence
we define:
\begin{definition}
\label{2.16}
1)  We say that $K$ is almost $(\mu,\kappa)$-closed under sums for
$\lambda$ and $\psi$ where $\psi = \psi(\bar x,\bar y),\ell g(\bar x)
= \ell g(\bar y)$, \Iff \, for every $I_\alpha\in K$ (for
$\alpha<\alpha_0\le\lambda),I_\alpha$ of cardinality $\le \lambda$, there are
$J,g,h_\alpha (\alpha<\alpha_0)$ such that:
\mn
\begin{enumerate}
\item[$(a)$]  $J \in K,|J| \le \lambda$,
\sn
\item[$(b)$]  $h_\alpha:I_\alpha \longrightarrow J$, and for any 
$x_0,\ldots,y_0,\ldots\in I_\alpha,
I_\alpha \models \psi[\langle x_0,\ldots\rangle,\langle y_0,\ldots\rangle]$
implies $J \models \psi[\langle h_\alpha(x_0),\ldots\rangle,
\langle h_\alpha(y_0),\ldots\rangle]$,
\sn
\item[$(c)$]  $g:J \longrightarrow \sum\limits_{\alpha<\alpha_0}\cM_{\mu,
\kappa}(I_\alpha)$ satisfies, for any $\gamma<\kappa,\bar x,\bar y \in
{}^\gamma J$ and $A \subseteq J$ of cardinality $<\kappa$,
\sn
\begin{enumerate}
\item[$\boxdot_0$]  if $g(\bar x) \approx g(\bar y) \mod 
\cM_{\mu,\kappa}(\sum\limits_{\alpha< \alpha_0} I_\alpha)$
then $\bar x \approx \bar y \mod \cM_{\mu,\kappa}(J)$.
\end{enumerate}
\end{enumerate}
\mn
2)  We replace ``almost" by ``semi", if in clause (c) above we 
weaken $\boxdot_0$ to:
\mn
\begin{enumerate}
\item[$\boxdot_1$]  if $g(\bar x) \approx g(\bar y) \mod 
(\cM_{\mu,\kappa}(\sum\limits_{\alpha<\alpha_0} I_\alpha),R)$
then $\bar x \approx \bar y  \mod \cM_{\mu,\kappa}(J)$, where we define

$R = \{\langle\langle\eta,i\rangle,\langle\nu,j\rangle\rangle:
\eta \in I_i,\nu \in I_j \text{ and } i< j\} \subseteq 
(\sum\limits_{\alpha<\alpha_0} I_\alpha) \times 
(\sum\limits_{\alpha<\alpha_0}I_\alpha)$.
\end{enumerate}
\mn
3)  We add ``strongly" to close in part (1) if we strengthen clause (c)
to: 
\mn
\begin{enumerate}
\item[$(c)^+$]  $g:J \longrightarrow \cM_{\mu,\kappa}(\sum\limits_{\alpha<
\alpha_0} I_\alpha)$ such that for any well ordering $<_0$ of 
$\cM_{\mu,\kappa}(J)$ (as in \ref{2.2}(d)), there is a well ordering $<_1$ of
$\cM_{\mu,\kappa}(\sum\limits_{\alpha<\alpha_0} I_\alpha)$ such that:
for any $\gamma<\kappa$ and $\bar x,\bar y \in {}^\gamma J$
and $A \subseteq J$ of cardinality $<\kappa$,
\sn
\begin{enumerate}
\item[$\boxdot_2$]  if $g(\bar x) \approx g(\bar y) \mod 
(\cM_{\mu,\kappa}(\sum\limits_{\alpha<\alpha_0} I_\alpha),<_1)$,
\then \, $\bar x \approx \bar y \mod (\cM_{\mu,\kappa}(J),<_0)$.
\end{enumerate}
\end{enumerate}
\mn
4) We add strongly in part (2) \Iff \, we strengthen (c) to (c)$^+$, only
using $(\cM_{\mu,\kappa}(\sum\limits_{\alpha<\alpha_0}
I_\alpha),<_1,R)$.

\noindent
5)  We may omit ``$(\mu,\kappa)$" above if $\Rang(g)\subseteq J$.

\noindent
6) We say that $K$ is essentially closed under sums for $\lambda$
\Iff \, in part (1) in addition, $\Rang(h_\alpha)$, $\Rang(g)$ are unions of
equivalence classes of $(R$ is from part (2))

\[
\approx\ \mod \ J,\qquad\approx\ \mod\ (\sum\limits_{\alpha<\alpha_0}
I_\alpha,R),\qquad \mbox{ respectively.}
\]
\end{definition}

\begin{remark} 
We could have made, for example $h_\alpha:I_\alpha
\longrightarrow \cM_{\mu,\kappa}(J)$,  or in the definition of sum expand
by $R$, without serious changes in the paper.
\end{remark}

\begin{claim}
\label{2.8}
0) ``$K$ is closed under sums'' implies ``$K$ is essentially closed
under sums", which implies ``$K$ is almost closed under sums", which
implies ``$K$ is almost $(\mu,\kappa)$-closed under sums". If $\mu_1 \le
\mu_2,\kappa_1 \le \kappa_2$ then ``$K$ is almost $(\mu_1,
\kappa_1)$-closed under sums" implies ``$K$ is $(\mu_2,\kappa_2)$-closed
under sums". 

In all above implications we can add ``strongly" to both sides (when
relevant, related).  

\noindent
1) If $K$ is closed under sums, \then \, the full (strong)
$(\chi,\lambda,\mu,\kappa)-\psi$-bigness property implies the (strong)
$(\chi_1,\lambda,\mu,\kappa)-\psi$-bigness property, where $\chi_1=\min
\{2^\chi,2^\lambda\}$.

\noindent
2) In (1), instead of ``$K$ closed under sums" it is enough to
assume that $K$ is (strongly) almost closed under sums for $\lambda$,
$\psi$.

\noindent
3) The classes defined in \ref{1.7} above $K^\kappa_{\tr},K_{\oor}$ are 
almost closed under sums and almost strongly closed under sums.

\noindent
4) The relations defined in \ref{2.3}(2), (3), (6) have obvious
monotonicity properties in $\chi,\mu,\kappa$; and for all our $K$, for
$\lambda$ too. For example

\[
\chi \le \chi \Rightarrow [(\chi',\lambda,\mu,\kappa)\text{-bigness}
\Rightarrow (\chi,\lambda,\mu,\kappa)\text{-bigness}]
\]

\[
\mu \le \mu' \& \kappa \le \kappa' \Rightarrow [(\chi,\lambda,\mu',\kappa')
\text{-bigness} \Rightarrow (\chi,\lambda,\mu,\kappa)\text{-bigness}].
\]
\end{claim}

\begin{PROOF}{\ref{2.8}}
0) Obvious.

\noindent 
1) So we assume $K$ has the full
$(\chi,\lambda,\mu,\kappa)-\psi$-bigness property.
\Wilog \, $ \langle I_\alpha:\alpha < \chi \rangle$ are pairwise disjoint.

As $K$ has the [strong] full $(\chi,\lambda,\mu,\kappa)-\psi$-bigness
property, there are $I_\alpha\in K$ (for $\alpha<\chi$), each of
cardinality $\lambda$, such that $I_\alpha$ is $\psi$-unembeddable intao
$\sum\limits_{\beta \ne \alpha} I_\beta$.  
\medskip

\noindent
\underline{Case 1}:  $\chi \le \lambda$.

For $U\subseteq\chi$ let $J_U=\sum\limits_{\alpha\in U}I_\alpha$.
Let $\cP$ be a collection of subsets of $\chi$ such that $|\cP|=
2^\chi$ and $U \ne  V \in \cP \Rightarrow U \nsubseteq V$. Suppose $U,V
\in \cP,f:J_U \longrightarrow M(J_V)$.  Choose $\alpha\in U\setminus
V$. Thus $f \rest I_\alpha:I_\alpha \longrightarrow \cM_{\mu,\kappa}
(\sum\limits_{\beta\ne\alpha} I_\beta)$ and the desired conclusion follows.
\medskip

\noindent
\underline{Case 2}: $\lambda<\chi$.

Take a family $\cW$ of subsets of $\lambda$, each of cardinality 
$\lambda$, such that 

\[
U \ne V \in H \Rightarrow U \nsubseteq V
\]
\mn
and proceed as in Case 1.
\medskip

\noindent
2) As $K$ has the [strong] full
$(\chi,\lambda,\mu,\kappa)-\psi$-bigness property, there are $I_\alpha\in K$
(for $\alpha<\chi$), each of cardinality $\lambda$, such that $I_\alpha$ is
$\psi$-unembeddable into $\sum\limits_{\beta\ne\alpha}I_\beta$. By the
assumption of (2) (that $K$ is almost (strongly) closed under sums) for
every $U \subseteq \chi,|U| \le \lambda$ let $J_U,g_U,h^U_\alpha$
($\alpha\in U$) satisfy  clauses (a), (b), (c) of Definition \ref{2.16}(1) for
$\sum\limits_{\alpha\in U}I_\alpha$. As in the proof of (1), it suffices to
show:
\mn
\begin{enumerate}
\item[$(*)$]  if $U,V \subseteq\chi,|U| \le \lambda,|V| \le \lambda,U
\setminus V \ne \emptyset$ and $f:J_U \longrightarrow \cM_{\mu,\kappa}(J_V)$,
\then \, for some $\bar a,\bar b \in {}^{\ell g(\bar x)}(J_U),
J_U \models \psi[\bar a,\bar b] \text{ and } f(\bar a) \approx_A
f(\bar b) \mod \cM_{\mu,\kappa}(J_V)$; or $\mod
(\cM_{\mu,\kappa}(J_V),<)$ for the strong version.
\end{enumerate}
\mn
Choose $\alpha\in U\setminus V$.

In the strong case let $<_0$ be a well ordering of
$\cM_{\mu,\kappa}(J_V)$ (as in \ref{2.2}(d),\ref{2.16}(3)); choose a 
well ordering $<_1$ of $\cM_{\mu,\kappa}(\sum\limits_{\alpha<\alpha_0}
I_\alpha)$ as guaranteed by Definition \ref{2.16}(3); in the
non-strong case let $<_0$, $<_1$ be the empty relations.

Now define

\[
g^*_V:\cM_{\mu,\kappa}(J_V) \longrightarrow \cM_{\mu,\kappa}(\sum_{i\in V}I_i)
\]

\mn
by

\[
g^*_V(\tau (x_0,\ldots))=\tau(g_V(x_0),\ldots).
\]

\mn
Consider the sequence of mappings:

\[
I_\alpha\mathop{\longrightarrow}\limits_{h^U_\alpha}J_U
\mathop{\longrightarrow}\limits_f\cM_{\mu,\kappa}(J_V)
\mathop{\longrightarrow}\limits_{g^\ast_V}\cM_{\mu,\kappa}\big(\sum_{i\in V}
I_i\big).
\]

\mn
So $g^*_V \circ f \circ h^U_\alpha:I_\alpha \longrightarrow \cM_{\mu,\kappa}
(\sum\limits_{i\in V}I_i)$. As $\sum\limits_{i\in V}I_i$ is a submodel of
$\sum\limits_{i\ne\alpha}I_i$, also without loss of generality $\cM_{\mu,
\kappa} (\sum\limits_{i\in V}I_i)$ is a submodel of $\cM_{\mu,\kappa}
(\sum\limits_{i\ne\alpha}I_i)$. But we know that $I_\alpha$ is
$\psi$-unembeddable into $\sum\limits_{i\ne\alpha}I_i$. Hence there are
$\bar x,\bar y\in I_\alpha$ such that:
\mn
\begin{enumerate}
\item[$(i)$]  $I_\alpha \models \psi [\bar x,\bar y]$,
\sn
\item[$(ii)$]  $g^*_V \circ f \circ h^U_\alpha(\bar x) \approx
g^*_V \circ f \circ h^U_\alpha(\bar y) \mod (\cM_{\mu,\kappa}
(\sum\limits_{i\in V}I_i),<_1)$.
\end{enumerate}
\mn
By (i) and clause $(b)$ from \ref{2.16}(1),
\mn
\begin{enumerate}
\item[$(iii)$]  $J_U \models \psi[\bar x',\bar y']$, where 
$\bar x'= h^U_\alpha(\bar x),\bar y'=h^U_\alpha(\bar y)$.
\end{enumerate}
\mn
By (ii) and the definition of $\bar x',\bar y'$,
\mn
\begin{enumerate}
\item[$(iv)$]  $g^*_V (f(\bar x')) \approx g^*_V(f(\bar y')) \mod 
(\cM_{\mu,\kappa}(\sum\limits_{i\in V}I_i),<_1)$.
\end{enumerate}
\mn
By (iv), clause (c) of \ref{2.16}(1) or clause $(c)^+$ 
of \ref{2.16}(3), the definition of $\cM_{\mu,\kappa}
(\sum\limits_{i\in V}I_i)$, and of $g^*_V$,
\mn
\begin{enumerate}
\item[$(v)$]  $f(\bar x') \approx f(\bar y') \mod (\cM_{\mu,\kappa}(J_V),<_0)$.
\end{enumerate}
\mn
So we have proved $(*)$ (by (iii) and (v)), which suffices.

\noindent
3)-6)  Left to the reader.
\end{PROOF}

\begin{claim}
\label{2.17}
The following classes are almost (and also semi) $(\mu,\kappa)$-closed
under sums for $\lambda$
\mn
\begin{enumerate}
\item[$(a)$]  $K_{\oor}$ (the class linear orders)  
\sn
\item[$(b)$]  $K^{\omega}_{\tr}$ (trees with $\omega+1$ levels)
\sn
\item[$(c)$]  $K^\kappa_{\tr}$ (trees with $\kappa+1$ levels)
\sn
\item[$(d)$]  $K_{\org}$ (ordered graphs). 
\end{enumerate}
\end{claim}

\begin{PROOF}{\ref{2.17}}
Case (a)

If $\langle I_\alpha:\alpha<\alpha_0\rangle$ is a sequence of
linear orders then we let:
\mn
\begin{enumerate}
\item[$(i)$]  $J = \cup\{\{\alpha\} \times I_\alpha:\alpha<\alpha_0\}$
\sn
\item[$(ii)$]  $(\alpha_1,t_1) <_J (\alpha_2,t_2)$ \underline{if and
    only if} $\alpha_1<\alpha_2 \vee (\alpha_1=\alpha_2\, \& \, 
t_1 <_{I_{\alpha_1}} t_2)$
\sn
\item[$(iii)$]  $h_\alpha:I_\alpha \rightarrow J$ is
  $h_\alpha(t)=(\alpha,t)$
\sn
\item[$(iv)$]  $g:J \rightarrow \sum\limits_{\alpha<\alpha_0} I_\alpha$
 is the identity.
\end{enumerate}
\mn
Now check
\medskip

\noindent
Case (b):

Given  $\langle I_\alpha:\alpha<\alpha_0\rangle$ the unique we identify the
member of $P^{J_\alpha}_0$ for $\alpha<\alpha_0$ but make then
otherwise disjoint and take the union.
\medskip

\noindent
Case (c):

Similar to case (b).
\medskip

\noindent
Case (d):

Similar to case (a).
\end{PROOF}

\noindent
Another way to present those matters is
to do it around the following definition and claim.
\begin{definition} 
\label{2.18n}
We say that $J_2$ does $(\mu,\kappa)$-dominate $J_1$ \when \,
there is a function $g$ from $\cM_{\mu,\kappa}(J_1)$ into
$\cM_{\mu,\kappa}(J_2)$ such that:
if $\rho \varphi \xi < \kappa$ and $\bar a,\bar b \in  
{}^{\xi}(\cM_{\mu,\kappa}(J_1))$ and $g(\bar a) \cong g(\bar b) \mod
\cM_{\mu,\kappa}(J_2)$ then $\bar a \cong \bar b \mod \cM_{\mu,\kappa}(J_1)$.

We say that $J_2$ strongly $(\mu,\kappa)$-dominate $J_1$ \underline{when}\
there is  a function $g$ from $\cM_{\mu,\kappa}(J_1)$ into
$\cM_{\mu,\kappa}(J_2)$ such that:
if $\xi < \kappa$ and $\bar a,\bar b \in {}^{\xi}(\cM_{\mu,\kappa }(J_1))$
and $g(\bar a) \cong g(\bar b) \mod \cM_{\mu,\kappa}(J_2)$
and $<_2$ is a well ordering of $(\cM_{\mu,\kappa}(J_2),<_2)$
then there is a well ordering $<_1$ of $\cM_{\mu,\kappa}(J_1)$
such that $\bar a \cong \bar b \mod (\cM_{\mu,\kappa}(J_1,<_1))$.

We say $J_1,J_2$ are [strongly] $(\mu,\kappa)$-equivalent when $J_2$ 
[strongly] dominate $J_1$ and vice versa.
\end{definition}

\begin{claim}   
\label{2.19n}
If $I$ is [strongly] $\varphi(\bar x,\bar y)$-unembeddable into $J_2$
and $J_2$ [strongly] $(\mu,\kappa)$-dominate $J_1$ \then \, $I$ is [strongly]
$\varphi(\bar x,\bar y)$-unembeddable into $J_2$.
\end{claim}
\bigskip

\centerline {$* \qquad * \qquad *$}
\bigskip

As we have remarked in the introduction to this paper,
results on trees can be translated to results on linear
orders;  this is done seriously in \cite{Sh:363}.
Originally this was neglected as the results on unsuperstable $T$ (and
trees with $\omega+1$ levels) give the results on unstable
theories (and linear orders). Anyhow, now we deal with the simplest case
parallel to \cite[Ch.2.1]{Sh:c}. 
\begin{definition}
\label{2.20}
1) For any $I \in K^\kappa_{\tr}$ we define
{\bf or}$(I)$ as the following linear order (See Def \ref{1.8}(4)).

\underline{set of elements} is chosen as
$\{(t,\ell):\ell\in \{1,-1\}, t\in I\}$

\underline{the order} is defined by $(t_1,\ell_1)<
(t_2,\ell_2)$ if and only if $t_1 \triangleleft t_2 \wedge \ell_1=1$
or $t_2 \triangleleft t_1 \wedge \ell_2 = -1$ or $t_1=t_2 \wedge
{\ell}_1 = -1 \wedge {\ell}_2 = 1$ or $t_1<_{\lx} t_2 \wedge (t_1,t_2$ 
are $\triangleleft$-incomparable. 

\noindent
2) Let $\varphi_{\oor} = \varphi_{\oor}(x_0,x_1; y_0,y_1)$
be the formula $x_0<x_1 \wedge y_1<y_0$.

\noindent
3) Let $\varphi^\kappa_{\tr} = \varphi^\kappa_{\tr}(x_0,x_1;y_0,y_1)$
be (this is for $K^\kappa_{\tr}$, for $\kappa=\aleph_0$ 
see example \ref{2.4A})

\begin{equation*}
\begin{array}{clcr}
\varphi_{\tr}(x_0,x_1:y_0,y_1) := &[x_0=y_0] \text{ and } 
P_\kappa(x_0)\wedge \bigvee\limits_{\epsilon <\kappa } [P_{\epsilon +1}(x_1)
\wedge P_{\epsilon +1}(y_1) \wedge \\
  & P_{\epsilon}(x_1 \cap y_1)] \wedge [x_1 \triangleleft x_0 \wedge
 \neg(y_1 \triangleleft y_0)] \text{ and } y_1<_{\lx} x_1].
\end{array}
\end{equation*}
\end{definition}

\begin{claim}
\label{2.21}
1)  Assume that $I,J \in K^\kappa _{\tr}$
\mn
\begin{enumerate}
\item[$(a)$]  If $I$ is strongly $\varphi^\kappa_{\tr}$-unembeddable
for $\tau_{\mu,\kappa}$ into $J$ then {\bf or}$(I)$ is strongly 
$\varphi^\kappa_{\tr}$-unembeddable for $\tau_{\mu,\kappa}$ into {\bf
  or}$(J)$
\sn
\item[$(b)$]  similarly without "strongly".
\end{enumerate}
\mn
2) If $K^\kappa _{\tr}$ has the strong
$(\chi,\lambda,\mu,\kappa)$-bigness property \then \, $K_{\oor}$ has the strong
$(\chi,\lambda,\mu,\kappa) $-bigness property.

\noindent
3) In part (2) we may  add "full" and/or omit "strong"
in the assumption and the conclusion.
\end{claim}

\begin{PROOF}{\ref{2.21}}
The main point is that:
\mn
\begin{enumerate}
\item[$(*)$]   if $I \models \varphi^\kappa_{\tr}(x_0,x_1;y_0,y_1)$
then {\bf or} $\models \varphi((x_0,1),(x_1,1); (y_0,1),(y_1,1))$.
\end{enumerate}
\end{PROOF}

\begin{remark}
\label{2.22}
1) We deal mainly with $K^\omega_{\tr}$, see
 \cite[3.1]{Sh:331}, so by it we know that $K^\omega_{\oor}$ has the full strong
$(\lambda,\lambda,\mu,{\aleph_0})$-bigness property when
$\mu<\lambda$.

\noindent
2)  For $\kappa$ regular uncountable, there are parallel results,
noting that obviously $K^\kappa_{\oor}$ have the full strong 
$(\chi,\lambda,\mu,\kappa)$ when $\lambda$ is regular $> |\alpha|^{<
  \kappa} + \mu$ for every $\alpha < \lambda$ and $\lambda \le \chi$.

It seems reasonable to conjecture that the parallel of 
\cite[3.1(2)]{Sh:331} holds, but we have not tried to work on it, see
part (3) of the remark.

\noindent
3)  The results below (on $\varphi_{\oor,\alpha,\beta,\pi}$)
seem to me a natural step but have actually set down to phrase and
prove them for Usvyatsov-Shelah \cite{ShUs:928}.

\noindent
4)  Even for $\kappa=\aleph_0$ we do not deal with $\lambda$
singular below, it seems reasonable that this, i.e., the parallel 
of \cite[\S1]{Sh:331} holds, but the results below are more 
than sufficient for its purpose, as for $\chi>\mu$ singular 
we can use the result here for
$(\chi,\lambda,\mu,\kappa)$ for any regular $\lambda\in (\mu,\chi)$.

\noindent
5) In \ref{2.15} we use $\alpha,\beta$ well orders.

It seems reasonable that we can say more for a more general case but
again this was not required.

\noindent
6) We use freely the obvious observation \ref{2.23}.
\end{remark}

\begin{observation}
\label{2.23}
1)  $K_{\oor}$ is essentially closed under sums for $\lambda$ and 
$\varphi_{\oor}$, recalling Definitions \ref{2.16}, \ref{2.17}.

\noindent
2)  Similar for $\varphi_{\oor,\alpha,\beta,\pi}$ defined below.
\end{observation}

\begin{definition}
\label{2.24}
We define the following (quantifier free infinitary)
formulas for the vocabulary $\{<\}$. For any ordinal
$\alpha,\beta$ and a one-to-one function $\pi$ from $\alpha$
onto $\beta$, and we let $\varphi_{\oor,\alpha,\beta,\pi}
(\bar x,\bar y)$ where $\bar x=\bar x^\alpha=\langle x_i:i<\alpha\rangle$
and $\bar y=\bar y^\alpha= \langle y_i:i<\alpha\rangle$, be

\[
\bigwedge \{x_i<x_j:i<j<\alpha\} \text{ and } 
\bigwedge \{y_i<y_j:i,j<\alpha \text{ and } \pi(i)<\pi (j)\}.
\]
\end{definition}

\begin{claim}
\label{2.25}
Assume $\chi \ge \lambda=\cf(\lambda)>\mu^{<\kappa},
\kappa=\cf(\kappa)$ and $\gamma<\lambda \Rightarrow |\gamma|^{<\kappa}
<\lambda$.

\noindent
1) For $(\alpha,\beta,\pi)$ as in Definition \ref{2.24}, such that
$\alpha,\beta \le \lambda$, the class $K_{\oor}$ has the full strong
$(\lambda,\chi,\mu,\kappa)$-bigness property for $\varphi_{\oor,
\alpha,\beta,\pi}(\bar x,\bar y)$.

\noindent
2) For $(\alpha,\beta,\pi)$ as in Definition \ref{2.24} such that
$\alpha,\beta \le \lambda$, the class $K_{\oor}$ has the strong
$(2^\lambda,\chi,\mu,\kappa)$ bigness property for 
$\varphi_{\oor,\alpha,\beta,\pi}$.

\noindent
3)  In fact in both part (1) and (2) we can find examples which 
satisfies the conclusion for all triples
$(\alpha,\beta,\pi)$ as there simultaneously.
\end{claim}

\begin{PROOF}{\ref{2.25}}
1) By \ref{2.26} below because there are $\lambda$ pairwise disjoint
stationary sets $S \subseteq S^\lambda_{\aleph_0}$.

\noindent
2) By part (1) and \ref{2.23}(1) and \ref{2.8}(1).

\noindent
3) Check the proof.
\end{PROOF}

\begin{claim}
\label{2.26}
Assume $\kappa = \cf(\kappa) \le \mu,\mu^{<\kappa}<\lambda=
\cf(\lambda) \le \lambda_1,\kappa \le \partial =\cf(\partial) 
< \lambda$ and $\gamma<\lambda \Rightarrow |\gamma|^{<\kappa}<\lambda$.

If $I,J \in K^\kappa_{\oor}$ satisfies $\circledast$ below and
$\alpha_*,\beta_*\le \lambda$ and $\pi$ is a one-to-one function 
from $\alpha_*$ onto $\beta_*$ \then \, (recalling Definition
\ref{2.20}) {\bf or}$(I)$ is strongly
$\varphi_{\oor,\alpha_*,\beta_*,\pi} (\bar x^{\alpha_*},
\bar y^{\alpha_*})$-unembeddable for $(\mu,\kappa)$ into {\bf or}$(J)$
where
\mn
\begin{enumerate}
\item[$\circledast$]  $(a) \quad 
S_1,S_2 \subseteq S^\lambda_\partial$ such that 
$S_1 \setminus S_2$ is a stationary subset of $\lambda$
\sn
\item[${{}}$]  $(b) \quad \bar \eta=\langle \eta_\delta:\delta\in S_1\cup S_2
\rangle$ where $\eta_\delta$ is an increasing sequence of ordinals

\hskip25pt $<\delta$ with limit $\delta$ of length $\partial$
\sn
\item[${{}}$]  $(c) \quad$ for every $\alpha<\lambda$ the set $\{\eta_\delta
\rest i:\delta\in S, i< \partial$ and $\sup \Rang(\eta_\delta \rest i)
\le \alpha\}$ 

\hskip25pt  has cardinality $<\lambda$
\sn
\item[$(d)$]  $I \in K^\kappa_{\tr}$ is 
$\{\eta_\delta \rest i:i \le \partial,\delta \in S_1\}
\cup\{\langle\alpha\rangle:\alpha<\lambda_1\}$
\sn
\item[$(e)$]  $J \in K^\kappa_{\tr}$ is $\{\eta_\delta \rest i:i
  \le \partial,\delta \in S_1\} \cup
  \{\langle\alpha\rangle:\alpha<\lambda_1\}$.
\end{enumerate}
\end{claim}

\begin{PROOF}{\ref{2.26}}
So let $f$ be a function from {\bf or}$(I)$ into
$\cM_{\mu,\kappa}$({\bf or}$(J)$) so actually a function from
$I\times \{1,-1\}$ into $\cM_{\mu,\kappa}(J \times \{1,-1\})$, and 
$<_*$ a well ordering of $\cM_{\mu,\kappa}(J)$  but we ``forget" 
to deal with it, as there are no problems, and let $\chi$
be large enough. Let $\bar N=\langle N_\alpha:\alpha<\lambda\rangle$ 
be an increasing continuous sequence of elementary submodels of 
$(\cH(\chi),\in)$ such that $I,J,\lambda,\bar \eta,\cM_{\mu,\kappa}(J),f,<_*$
belong to $N_0$ and $N_\alpha\cap \lambda\in \lambda,
\bar N \rest (\alpha+1)\in N_{\alpha +1}$ for every 
$\alpha<\lambda$; as it happens ``$\alpha_*,\beta_*,\pi\in N_0$" is not needed.
So $E:=\{\delta<\lambda:N_\delta \cap \lambda=\delta\}$ is
club of $\lambda$ hence we can choose $ \delta \in E \cap S_1 
\setminus S_2$.

For any $\eta\in I$, clearly $f((\eta,1))$ is well defined and
$\in \cM_{\mu,\kappa}(J)$ so let $f((\eta,1))=\sigma_\eta
(\bar \nu_\eta),\bar \nu_\eta=\langle (\nu_{\eta,\epsilon}, 
\iota_{\eta,\epsilon}):\epsilon< \epsilon_\eta\rangle,\nu_{\eta,i}\in
J$ and $\iota_{\eta,\epsilon} \in \{1,-1\},\epsilon_\eta < \kappa$.

Let $\epsilon_*=\epsilon_{\eta_\delta},\iota_\epsilon=
\iota_{\eta_\delta,\epsilon},i^*_\epsilon = {\rm lg}
(\nu_{\eta_\delta,\epsilon})$, so $i^*_\varp \le \partial$ 
for $\epsilon<\epsilon_*$ and let
$j^*_\epsilon = \sup\{j \le i^*_\epsilon:\sup \Rang
(\nu_{\eta_\delta,\epsilon} \rest j)<\delta\}$. By our assumption
$j^*_\epsilon = \partial$ implies that $i_\epsilon = \partial$
hence as $\delta \notin S_2$ it follows that $\sup
\Rang(\nu_{\eta_\delta,\epsilon}) <\delta$ hence
by clause (c) of the assumption $\nu_{\eta_\delta,\epsilon} \in N_\delta$.
Also $\alpha<\delta\Rightarrow J \cap {}^{\kappa >} \alpha 
\subseteq  N_{\alpha+1}$ because it has cardinality $<\lambda$
and it belongs to $N_{\alpha+1}$; also let $\nu^*_\epsilon =
\nu_{\eta_\delta,\epsilon} \rest j^*_\epsilon$, it too belongs to $N_\delta$.

So $\{\nu^*_\epsilon:\epsilon < \epsilon_*\} \subseteq N_\delta$, and 
it has cardinality $< \kappa$ as $\alpha < \lambda \rightarrow 
|\alpha|^{< \kappa } <  \lambda$ and $\cf(\delta) = \partial \ge 
\kappa$ it follows that $\bar\nu^*=  \langle \nu^*_\epsilon:\epsilon 
< \epsilon_* \rangle \in N_\delta$.

Let $u_*=\{\epsilon<\epsilon_*:j^*_\epsilon<i^*_\epsilon\}$.
For $\epsilon\in u_*$ let $\alpha^*_\epsilon = \min(N_\delta\cap 
(\lambda+1) \setminus \nu_{\eta_\delta,\epsilon}(j^*_\epsilon))$, 
so also $\bar{\alpha}^* := \langle \alpha_\epsilon:\epsilon \in 
u_* \rangle$ belongs to $N_\delta$.

Now for $\eta \in {}^{\partial >}\lambda$ we define
$\cU_\eta$ as the set of $\beta\in S_1$ such that:
\mn
\begin{enumerate}
\item[$(*)_{\eta,\beta}$]  $(a) \quad \eta \triangleleft \eta_\beta$
\sn
\item[${{}}$]  $(b) \quad \sigma_{\eta_\beta}=\sigma_*$ so
$\epsilon_{\eta_\beta}=\epsilon_*$
\sn
\item[${{}}$]  $(c) \quad$ {\rm lg}$(\nu_{\eta_\beta,\epsilon})= i^*_\epsilon$
for $\epsilon<\epsilon_*$
\sn
\item[${{}}$] $(d) \quad \nu_{\eta_\beta,\epsilon} \rest j^*_\epsilon
  = \nu^*_\epsilon$ for $\epsilon<\epsilon_*$
\sn
\item[${{}}$] $(e) \quad \iota_{\eta_\beta,\epsilon} = \iota_\epsilon$
for $\epsilon<\epsilon_*$
\end{enumerate}
\mn
Note
\begin{enumerate}
\item[$\circledast$]  if $\eta\triangleleft \eta_\delta$ then
\sn
\begin{enumerate}
\item[$(a)$]  $\delta \in \cU_\eta$ and $\cU_\eta\in N_\delta$
\sn
\item[$(b)$]  $\cf(\alpha^*_\epsilon)=\lambda$ for $\epsilon\in u_*$
\sn
\item[$(c)$]  if $\bar \alpha\in \prod\limits_{\epsilon\in u_*}
\alpha^*_\epsilon$ \then \, for arbitrarily large 
$\beta \in \cU_\eta$ we have
$\epsilon\in u_*\Rightarrow \nu_{\eta_\beta,\epsilon}
(j^*_\epsilon)\in (\alpha_\epsilon,\alpha^*_\epsilon)$
\sn
\item[$(d)$]  $\cU_\eta$ is an unbounded subset of $S_1$.
\end{enumerate}
\end{enumerate}
\mn
[Why?  Clause (a) directly.  Why clause (d)?  Otherwise $\sup(\cU_\eta)$ is
$< \lambda$ and it belongs to $N_\delta$ because $\cU_\eta \in
N_\delta$, hence $\sup(\cU_\eta) \in N_\delta \cap \delta$ so
$\sup(\cU_\eta) < \delta$ contradicting clause (a).  The other clauses
follows by them.]

Next let $\Lambda$ be the set of
$\eta \in {}^{\partial >} \lambda$ such that
\mn
\begin{enumerate}
\item[$\odot_\eta$]  for every $\bar \alpha\in \prod\limits_{\epsilon\in
u_*} \alpha^*_\epsilon$ there is $\beta \in \cU_\eta$
such that $\epsilon\in u_*\Rightarrow \nu_{\eta_\beta,\epsilon}
(j^*_\epsilon)\in (\alpha_\epsilon,\alpha^*_\epsilon)$.
\end{enumerate}
So
\mn
\begin{enumerate}
\item[$(*)_1$]  $\eta_1\triangleleft \eta_2 \in \Lambda 
\Rightarrow \eta_1\in \Lambda$
\sn
\item[$(*)_2$]  $\epsilon<\kappa \Rightarrow \eta_\delta 
\rest \epsilon\in \Lambda$.
\end{enumerate}
\mn
Hence
\mn
\begin{enumerate}
\item[$(*)_3$]  for some $\eta_*\in \Lambda$ the set
$\cW = \{\gamma<\lambda:\eta_* \char 94 
\langle \gamma\rangle\in \Lambda\}$ is an
unbounded subset of $\lambda$.
\end{enumerate}
\mn
Let $\langle \gamma_\zeta:\zeta<\lambda\rangle$ list
$\cW$ in increasing order, and let
$\alpha,\beta \le \lambda$ and $\pi$ be a one-to-one function from
$\alpha$ onto $\beta$.

Now first we choose $\delta (1,\zeta) \in S_1$ by induction
on $\zeta<\alpha$ such that
\mn
\begin{enumerate}
\item[$(*)_4$]  $(a) \quad \delta (1,\zeta) \in \cU_{\eta_* \char 94 
\langle \gamma_\zeta\rangle}$ i.e. $\gamma_\zeta\in \cW$
\sn
\item[${{}}$]  $(b) \quad$ if $\epsilon\in u_*$ then
$\nu_{\eta_{\delta(1,\zeta)},\epsilon} (j^*_\epsilon)$
is $<\alpha^*_\epsilon$ but is 
$> \sub\{\nu_{\eta_{\delta(1,\xi)},\epsilon}(j^*_\epsilon):\xi<\zeta\}$.
\end{enumerate}
\mn
This is easy.  

Second we choose $\delta(2,\zeta)\in S_1$ by induction on
$\zeta<\beta$ such that:
\mn
\begin{enumerate}
\item[$(*)_5$]  $(a) \quad \delta(2,\zeta)\in \cU_{\eta_* \char 94 \langle
\gamma_\xi\rangle}$ when $\pi(\xi)=\zeta$
\sn
\item[${{}}$]  $(b) \quad$ if $\epsilon \in u_*$ then 
$\nu_{\eta_{\delta(2,\zeta)},\epsilon}(j^*_\epsilon)$ is 
$<\alpha^*_\epsilon$ but is $> \sup\{\nu_{\eta_{\delta(2,\xi)},\epsilon}
(j^*_\epsilon):\xi<\zeta\}$.
\end{enumerate}
\mn
Let $\bar a=\langle a_\zeta:\zeta<\alpha\rangle$,
$\bar b=\langle b_\zeta:\zeta <\alpha\rangle$ from ${}^{\alpha}I$
be chosen as follows:
$a_\zeta=(\eta_{\delta(1,\zeta)},1),
b_\zeta=(\eta_{\delta(1,\pi(\zeta))},1)$ for $\zeta<\alpha$.

Now check, e.g.:
\mn
\begin{enumerate}
\item[$(*)_6$]  $a_{\zeta(1)} <_{\bold{\oor}(I)} a_{\zeta(2)}$ \Iff \,
$\gamma_{\zeta(1)}<\gamma_{\zeta(2)}$ \Iff \, $\zeta(1)<\zeta(2)$
\sn
\item[$(*)_7$]  $b_{\zeta(1)} <_{\bold{\oor}(I)} b_{\zeta(2)}$
\Iff \, $\gamma_{\pi(\zeta)(1)}<\gamma_{\pi(\zeta)(2)}$
\Iff \, $\pi(\zeta)(1)<\pi(\zeta)(2)$.
\end{enumerate}
\end{PROOF}

\begin{conclusion}
\label{2.27}
For $(\kappa,\mu,\lambda,\lambda_1,\alpha_*,\beta_*,\pi)$
as in \ref{2.26}, the class $K_{\oor}$ has the
full strong $(\lambda,\lambda_1,\mu,\kappa)-\varphi_{\oor,\alpha_*,
\beta_*,\pi}$-bigness property and the strong
$(2^\lambda,\lambda_1,\mu,\kappa)-
\varphi_{\oor,\alpha_*,\beta_*,\pi}$-bigness property.
\end{conclusion}

\begin{PROOF}{\ref{2.27}}
By \ref{2.26}.
\end{PROOF}
\newpage

\section {Order Implies Many Non-Isomorphic Models}

In this section (in a self contained way) we prove that not only the old
result that any unstable (first order) $T$ has in any $\lambda \ge |T|+
\aleph_1$, the maximal number $(2^\lambda)$ of pairwise non-isomorphic
models holds, but for example that for any template $\Phi$ proper for linear
orders, if the formula $\varphi(\bar x,\bar y)$ with vocabulary $\tau$,
linearly orders $\{\bar a_s:s\in I\}$ in $\EM_{\tau}(I,\Phi)$
(Ehrenfeucht-Mostowski model, see \S1) for every $I$, then the number of
non-isomorphic models of the form 
$\EM_{\tau}(I,\Phi)$ of cardinality $\lambda$ up to
isomorphism is $2^\lambda$ when $\lambda \ge |\tau_\Phi|+\aleph_1$.

Dealing with this problem previously, the author (in the first attempt
\cite{Sh:12}) excluded some of the cardinals $\lambda$ which satisfy
$\lambda = |\tau_\Phi|+\aleph_1$ and in the second \cite[Ch.VIII\S3]{Sh:a},
replaced the $\EM_\tau(I,\Phi)$ with some kind of restricted ultrapower (of
itself). Subsequently (\cite{Sh:100}) we proved that for some unsuperstable
first order complete theory $T$, and a first order theory $T_1$ extending
$T$, $|T_1|=\aleph_1$, $|T|=\aleph_0$ the class

\[
\PC(T_1,T)=\{M \rest \tau(T):M\models T_1\}
\]

\mn
may be categorical in $\aleph_1$, ``may be categorical" mean that
some forcing extension this holds for some $T,T_1$; in fact if
the original universe ${\bold V}$ satisfies $\CH$, we may choose
$T,T_1$ in ${\bold V}$.

We also prove there for $T=$ the theory of dense linear order,
that we may, i.e. in some forcing extension, have a
universal model in $\aleph_1$ even though $\CH$ fails. We then thought that the
use of ultrapower in \cite[Ch.VIII,\S3]{Sh:a} was necessary. This is not
true. (We thank  Rami Grossberg for 
a stimulating discussion which directed me to this problem again).

By the present theorem we can get the theorem also for the number of models
of $\psi \in \bbL_{\lambda^+,\aleph_0}$ in $\lambda$ ($>\aleph_0$) 
when $\psi$ is unstable. Incidentally the proof is considerably easier.

Note that we do not need to demand $\varphi(\bar x,\bar y)$ to be
first-order; a formula in any logic is O.K.; it is enough to demand
$\varphi(\bar x,\bar y)$ to have a suitable vocabulary. This is because an
isomorphism from $N$ onto $M$ preserves satisfaction of such $\varphi$ and
its negation. However, the length of $\bar x$ (and $\bar y$) is
crucial. Naturally we first concentrate on the finite case (in
\ref{3.1}--\ref{3.10}). But when we are not assuming this, we can, ``almost
always" save the result. In first reading, it may be advisable to
concentrate on the case ``$\lambda$ is regular".

For this section, the notion ``$\langle\bar a_t: t\in I\rangle$ is weakly
$(\kappa,\varphi(\bar x,\bar y)$-skeleton like inside $M$" is central
and in Definition \ref{3.1} the reader can concentrate on it.

\begin{definition}
\label{3.1}
Let $M$ be a model, $I$ an index model; for $s\in I$, $\bar a_s$ is a
sequence from $M$, the length of $\bar a_s$ depends
on the quantifier-free type of $s$ over $\emptyset$ in $I$ only; 
$\Lambda$ is a set of
formulas of the form $\varphi(\bar x,\bar a)$, $\bar a$ from $M$, $\varphi$
has a vocabulary contained in $\tau(M)$.

\noindent
1) We say that $\langle\bar a_s:s\in I\rangle$ is weakly
$\kappa$-skeleton like inside $M$ 
for\footnote{The simplest example is: $\Lambda$
the set of first order formulas with parameters from $M$.}
$\Lambda$ when: for every $\varphi(\bar x,\bar a) \in \Lambda$,   
there is $J \subseteq I$, $|J|<\kappa$ such that:
\mn
\begin{enumerate}
\item[$(*)$]  if $s,t \in I$ and $\tp_{\qf}(t,J,I)=\tp_{\qf}(s,J,I)$ then
\[
M \models ``\varphi[\bar a_s,\bar a]\equiv\varphi[\bar a_t,\bar a]".
\]
\end{enumerate}
\mn
2) If $\Lambda=\{\varphi(\bar x,\bar a):\varphi(\bar x,\bar y_\varphi)\in
\Delta,\\bar a\in \dot{\bold J}\}$ we may write 
$(\Delta,\dot{\bold J})$ instead of $\Lambda$; 
if $\Delta=\{\varphi(\bar x,\bar y)\}$ 
we write $\varphi(\bar x,\bar y)$ instead of $\Delta$. If

\[
\dot{\bold J} = \{\bar a:\bar a \text{ from } A,\text{ and for some }
\varphi(\bar x,\bar y) \in\Delta,\ \ell g(\bar a) = \ell g(\bar y)\}
\]

\mn
we write $A$ instead of $\Lambda$. If $|M|=A$ we write $M$ instead $A$, and we
omit it if clear from the context.

\noindent
3) Supposing $\psi(\bar x,\bar y)=:\varphi(\bar y,\bar x)$, $I$ a linear
order, we say $\langle \bar a_s:s\in I\rangle$ is weakly $(\kappa,\varphi
(\bar x,\bar y))$-skeleton like inside $M$ for $\dot{\bold J}$
\Iff \,: $\varphi(\bar x,\bar y)$ is asymmetric (at least in $M$)
with vocabulary contained in $\tau(M)$, $\ell g(\bar a_s) = 
\ell g(\bar x) = \ell g(\bar y),\langle \bar a_s:s\in I\rangle$ is weakly
$\kappa$-skeleton like inside $M$ for $(\{\varphi(\bar x,\bar y),\psi (\bar
x,\bar y)\},\dot{\bold J})$ and for $s,t\in I$ we have:

\[
M \models \varphi[\bar a_s,\bar a_t] \text{ iff }  I\models s<t.
\]

\mn
4) In (1), (3), if $M$ is clear from the context \then\ we may omit
``inside $M$".  In part (3), if $\dot{\bold J} = {}^\alpha |M|,
\alpha = \ell g(\bar x) = \ell g(\bar y)$ then we may omit it.
\end{definition}

\begin{discussion}
\label{3.1A}
Note that Definition \ref{3.1} requires considerably more than ``the $\bar
a_s$ are ordered by $\varphi$" and even than ``the $\bar a_s$ are order
indiscernibles ordered by $\varphi$", but much less than
``$M=\EM_\tau(I,\Phi)$".
\end{discussion}

We now would like to assign invariants to linear orders. We prove that there are
enough linear orders with well defined pairwise distinct invariants. This is
related to proofs from the Appendix to \cite{Sh:a}=\cite{Sh:c},
 where different terminology was employed. Speaking very roughly, 
we discussed there only $\inv^\alpha_\kappa$ where
$\kappa=\aleph_0$. The assertion in the appendix of \cite{Sh:c}
that two linear orders are contradictory corresponds to the assertion here
that the invariants are defined and different.

\begin{notation}  
\label{3.1B}  
In the following, for any regular cardinal $\mu>\aleph_0,D_\mu$ 
denotes the filter on $\mu$ generated by the closed unbounded sets.

\noindent
2)  If $D$ is a filter on $\mu$ and $X \subseteq \mu$
intersects each member of $D$, then $D+X$ denotes the filter generated by $D
\cup\{X\}$.

\noindent
3) For a linear order $I=(I,<_I)$ the cofinality $\cf(I)$ of $I$ is

\[
\Min\{|J|:J\subseteq I \text{ and }(\forall s \in I)(\exists t\in J) I\models
s<t\}.
\]

\mn
4) $I^*$ is the inverse linear order and $\cf^*(I)$ is the cofinality
of $I^*$.

\noindent
5)  For a linear order $I$ and a cardinal $\kappa$, let

\[
D = \cD(\kappa,I) := \cD_{\cf(I)} + \{\delta < \cf(I):\kappa \le 
\cf(\delta)\}.
\]

\mn
6) Two functions $f$ and $g$ from $\cf(I)$ to some set $X$, are equivalent
$\mod \, D$ if $\{\delta:f(\delta)=g(\delta)\}\in D$.

\noindent
7)  We write $f/D$ for the equivalence class of $f$ for this 
equivalence relations.
\end{notation}

\begin{definition}
\label{3.2}
1) For a regular cardinal $\kappa$ (for example $\aleph_0$) and an
ordinal $\alpha$ we define $\inv^\alpha_\kappa(I)$ for linear orders $I$
(sometimes undefined), by induction on $\alpha$, by cases:

$\alpha=0,\inv^\alpha_\kappa(I)$ is the cofinality of $I$ if $\cf(I)$ is
$\ge \kappa$, and is undefined otherwise.

$\alpha=\beta+1$

\noindent
Let $I=\bigcup\limits_{i<\cf(I)}I_i$, where $I_i$ is increasing and
continuous in $i$ and $I_i$ is a proper initial segment of $I$. For
$\delta<\cf(I)$ let $J_\delta=(I\setminus I_\delta)^*$
(where $X^*$ denotes the inverse order of $X$).
recalling \ref{3.1B}(4).

If $\cf(I)>\kappa$ and for some club $\cC$ of $\cf(I)$:
\mn
\begin{enumerate}
\item[$(*)_\cC$]  $[\delta\in \cC$ and $\cf(\delta)\ge \kappa]
\Rightarrow \inv^\beta_\kappa(J_\delta)$ is defined,
\end{enumerate}
\mn
\then \, we let

\[
\inv^\alpha_\kappa(I)=\langle\inv^\beta_\kappa(J_\delta):\cf(\delta)\ge
\kappa,\delta<\cf(I)\rangle/\cD(\kappa,I).
\]

\mn
Otherwise (i.e., there is no such $\cC$ or $\cf(I) \le \kappa$)
$\inv^\alpha_\kappa(I)$ is not defined.

$\alpha$ is limit

$\inv^\alpha_\kappa(I)=\langle\inv^\beta_\kappa(I): \beta<\alpha\rangle$.

\noindent
2) If $\bold d = \inv^\alpha_\kappa(I)$ \then \, ``the cofinality of
$\bold d$" means $\cf(I)$, clearly well defined.
\end{definition}

\begin{remark}
\label{3.2A}
1)  Really just $\alpha=0,1,2$ are used. For regular $\lambda$,
$\alpha=1$ suffices, but for singular $\lambda$, $\alpha=2$ is used (see
\ref{3.4}).

\noindent
2)  To understand the aim of \ref{3.3} below, think of $J$ as a linear
order such that for some linear order $U$, and $\langle \bar c_t:t\in U\rangle$
we have $\bar c_t\in {}^{\ell g(\bar x)}M$
and $\langle \bar a_{s}:s\in I \rangle \char 94 \langle \bar c_t:
t \in U\rangle$ and $\langle
\bar b_{t}:t\in U\rangle \char 94 
\langle\bar c_t:t\in U\rangle$ are both weakly
$(\kappa,\varphi(x,y))$-skeleton like in $M$ and $\cf(U^*) \ge
\kappa$.

\noindent
3)  We can omit assumption (c) in \ref{3.3}, so the conclusion will
tell us that if one of $\inv^\alpha_\kappa(I)$, $\inv^\alpha_\kappa(J)$ is
well defined then both are, but presently there is no real gain.

\noindent
4) The following lemma will be helpful as we will try to deal with cases of
inv inside models and try to prove that it is quite independent of a
(relevant) choice of representatives.
\end{remark}

\begin{observation}
\label{3.2B}
1)  If $\beta \le \alpha$ and $\inv^\alpha_\kappa (I)=
\inv^\beta_\kappa (J)$, and both are well defined \then \,
$\inv^\beta_\kappa(I),\inv^\beta_\kappa (J)$ are well defined and
equal.

\noindent
2) If $I,J$ are linear orders, $\inv = \inv^\alpha_\kappa(I)$ is well
defined, {\bf E} is a convex equivalence relation on $J$, $f:J\stackrel{\rm
onto}{\longrightarrow} I$ preserves $\le$, and $(f(x)=f(y))\equiv (x
{\bf E} y)$, \then \, $\bold d = \inv^\alpha_\kappa(J)$.

\noindent
3) Assume that $\psi(\bar x,\bar y)=\varphi(\bar y,\bar x)$ and 
$\varphi_{\ell}(\bar x,\bar y) \in \{\varphi(\bar x,\bar y),
\neg \varphi(\bar x,\bar y),\psi(\bar x,\bar y),\neg \psi(\bar x,\bar y)\}$
for ${\ell} = 1,2$. \Then \, $\langle \bar a_s:s \in I \rangle$
is weakly $(\kappa,\varphi_1(\bar x,\bar y))$-skeleton like in $M$ if
and only if $\langle \bar a_s:s \in I^* \rangle$ is weakly 
$(\kappa,\varphi_2(\bar x,\bar y))$-skeleton like in $M$;
also in $M$ we have $\varphi(\bar x,\bar y) \vdash \neg \psi(\bar
x,\bar y)$ and $\psi(\bar x,\bar y) \vdash \neg \varphi (\bar x,\bar y)$.
\end{observation}

\begin{lemma}
\label{3.3}
Suppose that $\kappa$ is a regular cardinal, $I,J$ are linear orders, and
$\bar a_s$ (for $s\in I$), $\bar b_t$ (for $t\in J$) are from $M$, and
$\varphi (\bar x,\bar y)$ is a $\tau (M)$-formula ($\kappa > \ell g(\bar x)=
\ell g(\bar y) =\ell g(\bar a_s) = \ell g(\bar b_t)$), and 
$\psi(\bar x,\bar y) := \varphi(\bar y,\bar x)$. 

Assume:
\mn
\begin{enumerate}
\item[$(a)$]   $(\alpha) \quad$ for every $s\in I$ for every 
large enough $t\in J,M \models \varphi[\bar a_{s},\bar b_{t}]$,
\sn
\item[${{}}$]  $(\beta) \quad$ for every $t\in J$ for every 
large enough $s\in I,M \models \neg \varphi[\bar a_s,\bar b_t]$,
\sn
\item[$(b)$]   $(\alpha) \quad \langle\bar a_s:s \in I \rangle$ is 
weakly $(\kappa,\varphi(\bar x,\bar y))$-skeleton like inside $M$,
\sn
\item[${{}}$]  $(\beta) \quad \langle\bar b_t:t\in J\rangle$ is weakly $(\kappa,
\varphi(\bar x,\bar y))$-skeleton like inside $M$,
\sn
\item[$(c)$]  $\inv^\alpha_\kappa(I)$, $\inv^\alpha_\kappa(J)$ are defined.
\end{enumerate}
\mn
\Then \, $\inv^\alpha_\kappa(I)=\inv^\alpha_\kappa(J)$.
\end{lemma}

\begin{PROOF}{\ref{3.3}}
By induction on $\alpha$.
\medskip

\noindent
\underline{First Case}: $\alpha=0$

Assume not, so $\inv^0_\kappa(I) \ne \inv^0_\kappa(J)$. Then $\cf(I),
\cf(J)$ are distinct (and $\ge \kappa$).
By symmetry, \wilog \, $\cf(I)>\cf(J)$, so $\cf (I)>\kappa$.

Let $\langle t_\zeta:\zeta<\cf(J)\rangle$ be increasing unbounded in
$J$. For each $\zeta<\cf(J)$ (by clause (a)$(\beta)$ of \ref{3.3}
and \ref{3.2B}) there is $s_\zeta\in I$ such that:

\[
s_\zeta \le s \in I \Rightarrow  M \models \neg\varphi[\bar a_s,b_{t_\zeta}].
\]

\mn
As $\cf(I)>\cf(J)$ there is $s\in I$ such that $\bigwedge\limits_\zeta
s_\zeta <s$.  Now, the set

\[
\{t \in J:M \models \neg \varphi[\bar a_s,\bar b_t]\}
\]

\mn
includes each $t_\zeta$ (as $s_\zeta< s\in I$), and hence it is unbounded in
$J$, contradicting clause (a)$(\alpha)$ of \ref{3.3}.
\medskip

\noindent
\underline{Second Case}:  $\alpha=\beta+1$

By the first case and Observation \ref{3.2B}, $\cf(I) = \cf(J) \ge 
\kappa$.  Let $\lambda=\cf(I)=\cf(J)$; let

\[
I = \bigcup\limits_{i<\lambda} I_i,
\]

\mn
where $I_i$ is increasing continuous in $i$, $I_i$ a proper initial segment
of $I$ and $[i \ne j \Rightarrow I_i \ne I_j]$. 

Similarly let

\[
J=\bigcup\limits_{i<\lambda} J_i.
\]

\mn
Choose $s_i\in I_{i+1}\setminus I_i$ and $t_i\in J_{j+1}\setminus J_j$. By
assumption (a), for every $i<\lambda$ there is $j_i<\lambda$ such
that:
\mn
\begin{enumerate}
\item[$(\alpha)'$]  if $t\in J\setminus J_{j_i}$ then $M\models\varphi [\bar
a_{s_i},\bar b_t]$,
\sn
\item[$(\beta)'$]  if $s\in I\setminus I_{j_i}$ then $M\models\neg\varphi
[\bar a_s,\bar b_{t_i}]$.
\end{enumerate}
\mn
Let

\[
\cC = \{\delta<\lambda:\delta \text{ is a limit ordinal and } i<\delta
\Rightarrow j_i<\delta\};
\]

\mn
it is a club of $\lambda$. For $\delta\in \cC$ let $I^\delta=(I\setminus
I_\delta)^*$ and let $J^\delta=(J\setminus J_\delta)^*$. By Definition
\ref{3.2} above it suffices to prove, for $\delta\in \cC$ satisfying
$\cf(\delta) \ge \kappa$ such that $\inv^\beta_\kappa(I^\delta),
\inv^\beta_\kappa(J^\delta)$ are defined, that:
\mn
\begin{enumerate}
\item[$(*)_\delta$]  $\inv^\beta_\kappa(I^\delta)=\inv^\beta_\kappa(J^\delta)$.
\end{enumerate}
\mn
For this we use the induction hypothesis, but we have to check that the
assumptions (a), (b), (c) hold for this case.

Now clause (c) is part of the assumption of $(\ast)_\delta$, and clause (b)
is inherited from the same property of $\langle\bar a_s:s\in I\rangle$,
$\langle \bar b_t:t\in J\rangle$; lastly clause (a) follows from
$(\alpha)' + (\beta)'$ above as $\delta \in \cC$. In detail, if
$t\in J^\delta$ then $J \models ``t_j<t"$ for $j<\delta$. Hence, for
$i<\delta,M\models \varphi[\bar{a}_{s_i},\bar{b}_t]$ (by clause
$(\alpha)'$ above). So by clause (b)$(\beta)$ from the assumptions, 
for every large enough $s\in I^\delta$ we have $M \models 
\varphi[\bar{a}_s,\bar{b}_t]$, which means that $\langle\bar{a}_s:
s\in I^\delta\rangle$, $\langle\bar{a}_t:t\in J^\delta\rangle$ satisfy
clause (a)$(\alpha)$.  Similarly clause (a)$(\beta)$ holds.
\medskip

\noindent
\underline{Third Case}:  $\alpha$ is limit

Immediate by Definition \ref{3.2}. 
\end{PROOF}

\begin{lemma}
\label{3.4}
1)  If $\lambda,\kappa$ are regular, $\lambda>\kappa$,
\then \, there are $2^\lambda$ linear orders $I_\alpha$ (for
$\alpha<2^\lambda$), each of cardinality $\lambda$, with pairwise distinct
$\inv^1_\kappa(I_\alpha)$ (for $\alpha<2^\lambda$), each well defined.

\noindent
2)  If $\lambda>\kappa$, $\kappa$ is regular, \then \, there
are linear orders $I_\alpha$ (for $\alpha<2^\lambda$), each of cardinality
$\lambda$ with pairwise distinct $\inv^2_\kappa(I_\alpha)$ (for $\alpha<
2^\lambda$), each well defined.

\noindent
3)  If in (2) we have $\lambda\geq\theta=\cf(\theta)>\kappa$,
\then \, we can have $\cf(I_\alpha)=\theta$ if we use 
$\inv^3_\alpha$. Similarly, if in part (1) we have $\lambda \ge 
\theta=\cf(\theta)>\kappa$, \then\ we can have
$\cf(I_\alpha)=\theta$ if we use $\inv^2_\kappa$; of course can use
$\inv^\alpha_\kappa$ for $\alpha \ge 2$ (similarly elsewhere).

\noindent
4) Assume $\Phi$ is an almost $\cL$-nice template proper for
linear orders (see Definition \ref{1.6}). \Then \, for any linear
order $I$, the sequence $\langle \bar a_t: t\in I\rangle$ is
$\aleph_0$-skeleton like for $\cL$ inside $\EM(I,\Phi);\cL$ can
be any set of formulas in the vocabulary $\tau_\Phi$.

\noindent
5) In part (4), if $I$ is $\aleph_0$-homogeneous (i.e., for any 
$n<\omega$ and $t_0<_I\ldots<_I t_{n-1},s_0<_I\ldots <_I s_{n-1}$, 
there is an automorphism of $I$ mapping $t_\ell$ to $s_\ell$ 
for $\ell<n$), then we can omit ``almost $\cL$-nice".
\end{lemma}

\begin{remark}
\label{3.4d}
1)  The construction of the linear orders is ``hinted" by the proof 
\ref{3.3}, and by the properties of stationary sets. 
Alternatively see the inductive
construction in Claims 3.7, 3.8 of the
Appendix of \cite{Sh:a} or see \cite{Sh:12} where $\inv^\alpha_\kappa
(1), \alpha<\lambda^+, \lambda=|I|$ are used.

\noindent
2)  Note that part (4) says that being skeleton-like really is 
a property of the skeleton of $\EM$-models.

\noindent
3)  Note that \ref{3.4}(4) apply to $\EM_\tau(I,\Phi)$ 
whenever $\tau \subseteq \tau_\Phi$.
\end{remark}

\begin{PROOF}{\ref{3.4}}
1) So $\lambda>\kappa$ are regular. The set $S=\{\delta<\lambda:
\cf(\delta)=\kappa\}$ is stationary and hence we can find a partition
$\langle S_\epsilon:\epsilon<\lambda\rangle$ of $S$ into pairwise
disjoint stationary subsets (well known, see
Solovay theorem). For $u \subseteq \lambda$ we define
$I_u$ as the set

\[
\{(\alpha,\beta):\alpha < \lambda \text{ and }\alpha \in
\bigcup\limits_{\epsilon\in u} S_\epsilon \Rightarrow \beta < \kappa^+
\text{ and } \alpha\in\lambda\setminus
\bigcup\limits_{\epsilon\in u} S_\epsilon 
\Rightarrow \beta<\kappa\}
\]
linearly ordered by

\[
(\alpha_1,\beta_1) <_I (\alpha_2,\beta_2) \text{ iff } \alpha_1<
\alpha_2 \vee (\alpha_1 = \alpha_2 \text{ and } \beta_1>\beta_2).
\]

\mn
By the proof of \ref{3.3} above clearly $\langle I_u:u\subseteq
\lambda \rangle$ is as required.

\noindent
2) So we have $\lambda>\kappa$, $\kappa=\cf(\kappa)$.

\noindent 
Let $\lambda=\sum\limits_{i<\cf(\lambda)}\lambda_i$, $\lambda_i$
increasing continuous $>\kappa$, let $\theta=\cf(\lambda)+\kappa^+$, or just
$\kappa^++\cf(\lambda) \le \theta=\cf(\theta) \le \lambda$. Let $h:\theta
\longrightarrow\cf(\lambda)$ be such that for any $i<\cf(\lambda)$ the set
$\{\delta<\theta:\cf(\delta)=\kappa$ and $h(\delta)=i\}$ is stationary.

For each $i$, let $\langle I_{i,\epsilon}:\epsilon<2^{\lambda^+_i}
\rangle$ be as in the proof of (1) (for $\lambda^+_i$). For any $\nu\in
\prod\limits_{i<\cf(\lambda)} 2^{\lambda^+_i}$ let $J_\nu=\sum\limits_{
\alpha<\theta} J^*_{\nu,\alpha}$ with $J_{\nu,\alpha}\cong I_{h(\alpha),
\nu(\alpha)}$.

\noindent
3) Let $\langle I_\epsilon:\epsilon<2^{\lambda}\rangle$
be as guaranteed in part (2) (or part (1) if $\lambda$ is regular). 
For each $\epsilon<2^\lambda$, let $J_\epsilon = 
\sum\limits_{i<\theta} J^*_{\epsilon,i}$ where
$J_{\epsilon,i}\cong I_\epsilon$; now the sequence $\langle
I_\epsilon: \epsilon < 2^\lambda \rangle$ is as required.

\noindent
4) Let $\varphi=\varphi(\bar x,\bar b) \in \cL(\tau_\Phi)$, so
for some finite sequence $\bar t$ from $I$ and a sequence $\bar\sigma$ of
$\tau_\Phi$-terms we have $\bar b=\bar{\sigma}(\bar t)$. So if $s_1$, $s_2$
realize the same quantifier free type over $\bar t$ in $I$, 
by indiscernibility (i.e., almost $\cL$-niceness) $\EM(I,\Psi)\vDash
``\varphi[\bar a_{s_1},\bar b] = \varphi[\bar a_{s_2},\bar b]"$. 
So $\rang(\bar t)$ is as required.

\noindent
5) Should be clear.  
\end{PROOF}
\bigskip

\centerline {$* \qquad * \qquad *$}
\bigskip

Now we would like to attach the invariants of a linear order $I$ to a model $M$
which has a skeleton-like sequence indexed by $I$. In $(\alpha)$ (in
Definition \ref{3.5} below) we define what it means for a sequence indexed
by $I$ to $(\kappa,\theta)$-represent the $(\varphi,\psi)$-type of
$\bar c$ over $A$.
\begin{definition}
\label{3.5}
Let $A \subseteq M,\bar c\in M$ and $\varphi(\bar x,\bar y)$ be an
asymmetric formula with vocabulary contained in $\tau(M)$ and $\psi(\bar x,
\bar y)=:\varphi(\bar y,\bar x)$
\mn
\begin{enumerate}
\item[$(\alpha)$]  We say that $\langle\bar a_s:s\in I\rangle$ does $(\kappa,
\theta)$-represents $(\bar c,A,M,\varphi(\bar x,\bar y))$ \Iff \,:
$I$ is a linear order, $\cf(I)\ge\kappa$ and for some linear order $J$ of
cofinality $\theta\ge\kappa$ disjoint to $I$, there are $\bar a_t\in
{}^{\ell g(\bar x)}A$ for $t \in J$, such that:
\sn
\begin{enumerate}
\item[$(i)$]  for every large enough $t\in I$, $\bar a_t$ realizes $\tp_{\{
\varphi({\bar x},{\bar y}),\psi({\bar x},{\bar y})\}}(\bar{c},A,M)$,
and
\sn
\item[$(ii)$]  $\langle\bar a_s:s\in J+(I)^\ast\rangle$ is weakly $(\kappa,
\varphi (\bar x,\bar y))$-skeleton like inside $M$ ($I^*$ denotes the
inverse of $I$).
\end{enumerate}
\sn
\item[$(\beta)$] We say that $(\bar c,A,M,\varphi (\bar x,\bar y))$ has a
$(\kappa,\theta,\alpha )$-invariant \when \,:
\sn
\begin{enumerate}
\item[$(i)$]  if for $\ell=1,2,\langle\bar a^\ell_s:s\in I_\ell\rangle$
does $(\kappa,\theta)$-represents $(\bar c,A,M,\varphi(\bar x,\bar y))$ and
$\inv^\alpha_\kappa(I_\ell)$ are defined\footnote{but see \ref{3.8A}(2)}
for $\ell = 1,2$ \then \,
$\inv^\alpha_\kappa(I_1)=\inv^\alpha_\kappa(I_2)$,
\sn
\item[$(ii)$]  some $\langle\bar a_s:s\in I\rangle$ does $(\kappa,
\theta)$-represent $(\bar c,A,M,\varphi(\bar x,\bar y))$ and
$\inv^\alpha_\kappa(I)$ is well defined.
\end{enumerate}
\sn
\item[$(\gamma)$]  Let  $\INV^\alpha_{\kappa,\theta}(\bar c,A,M,\varphi(\bar
x,\bar y))$ be $\inv^\alpha_\kappa(I)$ when $(\bar c,A,M,\varphi(\bar x,\bar
y))$ has $(\kappa,\theta,\alpha)$-invariant and $\langle\bar a_s:s\in I
\rangle$ does $(\kappa,\theta)$-represent it
\sn
\item[$(\delta)$]  Let ``$(\kappa,\alpha)$-invariant" means ``$(\kappa,
\theta,\alpha)$-invariant for some regular $\theta\ge\kappa$".  Similarly
for ``$\kappa$-represents" and $\INV^\alpha_\kappa(\bar c,A,M,\varphi(\bar
x,\bar y))$ (justified by Fact \ref{3.5A} below).
\end{enumerate}
\end{definition}

\begin{fact}
\label{3.5A}
Suppose that for $\ell=1,2$, the sequence $\langle\bar a^\ell_s:s\in I_\ell
\rangle$ does $(\kappa,\theta_\ell)$-represent $(\bar c,A,M,\varphi (\bar x,
\bar y))$. \Then \, $\theta_1=\theta_2$.
\end{fact}

\begin{PROOF}{\ref{3.5}}
So let for $\ell=1,2$ the sequence $\langle\bar a^\ell_s:s\in J_\ell
\rangle$ witness that $\langle\bar a^\ell_s:s\in I_\ell\rangle$ does $(\kappa,
\theta_\ell)$-represent $(\bar c,A,M,\varphi(\bar x,\bar y))$, i.e., they
are as in $(\alpha)$ of \ref{3.5}. Assume toward contradiction that
$\theta_1 \ne \theta_2$ and by symmetry \wilog \, $\theta_1<\theta_2$. Let
$\langle s_\ell(\alpha):\alpha<\theta_\ell\rangle$ be an increasing
unbounded sequence of members of $J_\ell$ for $\ell=1,2$. 
So for each $\alpha<\theta_1$ we have

\[
t\in I_1\quad \Rightarrow\quad M \vDash \varphi[\bar a^1_{s_1(\alpha)},\bar
a_t^1]
\]

\mn
and hence by clause (i) of $(\alpha)$ of Definition \ref{3.5} we have
$M \vDash \varphi[\bar a^1_{s_1(\alpha)},\bar{c}]$
recalling $\bar{a}^1_{s_1(\alpha)} \subseteq A$, so for every 
large enough $t\in I_2,M \vDash \varphi[\bar a^1_{s_1(\alpha)},
\bar a^2_t]$. But $\langle \bar a^2_t:t\in J_2+(I_2)^*\rangle$ 
is weakly $(\kappa,\varphi (\bar x,\bar y))$-skeleton like inside $M$,
hence for some $\beta_\alpha<\theta_2$ we have

\[
s_2(\beta_\alpha) \le t\in J_2 \Rightarrow M \vDash 
\varphi[\bar a^1_{s_1(\alpha)},\bar a^2_t]
\]

\mn
and so $\beta(*)=\sup\{\beta_\alpha+1:\alpha<\theta_1\}<\theta_2$ (as
$\theta_1<\theta_2 = \cf(\theta_2)$). So $M \vDash \varphi[\bar a^1_{s_1(
\alpha)},\bar a^2_{s_2(\beta(*))}]$ for $\alpha<\theta_1$.

But $t\in I_2\ \Rightarrow\ M\vDash\neg\varphi [\bar a^2_t,\bar a^2_{s_2(
\beta)}]$ and hence $M\vDash\neg \varphi[\bar c,\bar a^2_{s_2(\beta)}]$.
Therefore, for every large enough $t\in I_1,M \vDash \neg\varphi[\bar
a^1_t,\bar a^2_{s_2(\beta)}]$ and hence for every large enough $t\in J_1$,
$M \vDash \neg \varphi[\bar a^1_t,\bar a^2_{s_2(\beta)}]$. Hence this holds
for $t=s_1(\alpha)$, $\alpha$ large enough, a contradiction to the
previous paragraph.
\end{PROOF}

\begin{discussion}
\label{3.6y}
Each of Definition \ref{3.6}, Lemmas \ref{3.7} and \ref{3.8}, and the proof
of Theorem \ref{3.9} have 3 cases. In the easiest case $\lambda=\|M\|$ is
regular. When $\lambda$ is singular the computation of $\inv^\alpha_\kappa
(\bar{c},A,M,\varphi(\bar{x},\bar{y}))$ is easier when $\cf(\lambda
)>\kappa$ (second case). 
The third case arises when $\lambda>\kappa > \cf(\lambda)$.

The relative easiness of the regular case is caused by the fact that any two
increasing representations of a model with cardinality $\lambda$ must ``agree"
on a club. In the second case we are able to restrict the first argument to
a cofinal sequence of $M$. For the third case we must construct a ``dual
argument'', noticing that much of a long sequence must concentrate on one
member of the representation.
\end{discussion}

\begin{definition}
\label{3.6}
Let $\varphi(\bar x,\bar y)$ be an asymmetric formula with vocabulary
$\subseteq\tau(M)$ (where $\ell g(\bar x) = \ell g(\bar y)$ is
finite),  and let $M$ be a model of cardinality
$\lambda,\lambda>\kappa,\kappa$ regular, $\alpha$ be an ordinal.

\noindent
0)  A representation of the model $M$ is an increasing continuous
sequence $\bar M=\langle M_i: i<\cf(\lambda)\rangle$ such that $\|M_i\|<
\lambda$, and $M=\bigcup\limits_{i<\cf(\lambda)} M_i$.

Similarly for sets.

\noindent
1) For a regular cardinal $\lambda$: 

\begin{equation*}
\begin{array}{clcr}
{\bf INV}^\alpha_\kappa(M,\varphi(\bar x,\bar y)) = \{\bold d: &\text{
  for every representation } \langle A_i:i<\lambda\rangle \text{ of }|M|,\\
  &\text{ there are }\delta < \lambda \text{ and } \bar c \in M 
\text{ (of course, } \ell g(\bar c) = \ell g(\bar x) \\
  &\text{ such that } \cf(\delta) \ge \kappa \text{ and } \bold d=
\INV^\alpha_\kappa(\bar c,A_\delta,M,\varphi(\bar x,\bar y)) \\
  &\text{ (in particular so the latter is well defined) }\}.
 \end{array}
\end{equation*}

\mn
2)  For regular cardinals $\theta>\kappa$ such that $\lambda>\cf(\lambda)=
\theta$ we let

\[
\cD_{\theta,\kappa} = \cD_\theta + \{\delta<\theta:\cf(\delta) \ge \kappa\}
\]

\mn
and

\begin{equation*}
\begin{array}{clcr}
\bf{INV}^\alpha_{\kappa,\theta}(M,\varphi(\bar x,\bar y)) = 
\{\langle \bold d_i:i<\theta\rangle/\cD_{\theta,\kappa}: &\text{ for every
representation } \langle A_i:i<\theta\rangle \text{ of }|M|,\\
  &\text{ there is } S \in \cD_{\theta,\kappa} \text{ satisfying:}\\
  &\text{ for every } \delta \in S \text{ there is } \bar c_\delta \in
  M \text{ such that}\\
  & \bold d_\delta = \INV^\alpha_\kappa(\bar c_\delta,A_\delta,M,\varphi
(\bar x,\bar y))\\
  &\text{ so is well defined and the cofinality of } \bold d_\delta 
\text{ is } > |A_\delta|\}.
  \end{array}
 \end{equation*}

\mn
3)  For regular cardinals $\kappa>\theta$, $\lambda>\theta>\kappa+\cf(
\lambda)$ and a function $h$ with domain a stationary subset of $\{\delta<
\theta:\cf(\delta)\ge\kappa\}$ and range a set of regular cardinals
$<\lambda$, we let 

\[
\cD_{\theta,h} = \cD_\theta + \{\{\delta<\theta:h(\delta) \ge \mu
\quad (\text{hence } \delta \in \Dom(h))\}:\mu<\lambda\},
\]

\mn
and assuming that $\cD_{h,\lambda}$ is a proper filter we let:

\begin{equation*}
\begin{array}{clcr}
\bf{INV}^{\alpha,h}_{\kappa,\theta}(M,\varphi(x,y)) = \{\langle 
\bold d_i:i < \theta\rangle/\cD_{\theta,h}:&\text{ for every
  representation } \langle A_i:i<\cf(\lambda)\rangle \text{ of }|M|,\\
 &\text{ there are }\gamma < \cf(\lambda) \text{ and } S \in
\cD_{h,\lambda},S\subseteq\Dom(h),\text{ satisfying} \\
  &\text{ the following for each } \delta \in S,\text{ if } 
h(\delta)> |A_\gamma| %%\text{ \then \, for some } \bar c_\delta\in M,
\\
  & \text{ \then \, for some } \bar c_\delta\in M,\\
&\text{ we have } {\bf d}_\delta = \INV^\alpha_\kappa(\bar c_\delta,
A_\gamma,M,\varphi(\bar x,\bar y))\\
  &\text{ so is well defined and the cofinality of } e_\delta
\text{ is } >|A_\gamma|\}.
 \end{array}
\end{equation*}
\end{definition}

\begin{remark}
\label{3.6A}
1)  Of course, also in \ref{3.6}(1) we could have used $\langle \bold d_i:i<
\lambda\rangle/\cD_\lambda$ as the invariant.

\noindent
2) In \ref{3.6}(3), we  may demand ``$\cf(\bold d_\delta) > |A_\delta|$".
\end{remark}

\begin{lemma}
\label{3.7}
Suppose $\varphi(\bar x,\bar y)$ is a formula in the vocabulary of $M$,
$\ell g(\bar x) = \ell g(\bar y)<\omega$.

\noindent
1)  If $\lambda>\aleph_0$ is regular, $M$ a model of cardinality
$\lambda$, $\kappa$ regular $<\lambda$, \then \, 
$\bf{INV}^\alpha_\kappa(M,\varphi(\bar x,\bar y))$ has cardinality ! $\le\lambda$.

\noindent
2) If $\lambda$ is singular, $\theta=\cf(\lambda)>\kappa$,
\then \, $\bf{INV}^\alpha_{\kappa,\theta}(M,\varphi(\bar x,\bar y))$
almost has cardinality $\le \lambda$, which means: there are no 
$\bold d^\zeta_i$ (for $i<\theta$, $\zeta<\lambda^+$) such that:
\mn
\begin{enumerate}
\item[$(i)$]   for $\zeta < \lambda^+,\langle \bold d^\zeta_i:
i<\theta\rangle/\cD_{\theta,\kappa} \in \bf{INV}^\alpha_{\kappa,\theta}
(M,\varphi(\bar x,\bar y))$,
\sn
\item[$(ii)$]  for $i<\theta,\zeta<\xi<\lambda^+$, we have 
$\bold d^\zeta_i\ne \bold d^\xi_i$.
\end{enumerate}
\mn
3) If $\lambda $ is singular, $\theta,\kappa$ are regular, $\kappa+\cf(
\lambda)<\theta<\lambda$, $h$ is a function from some stationary subset of
$\{i<\theta:\cf(i) \ge \kappa\}$ into

\[
\{\mu<\lambda:\mu \text{ is a regular cardinal }\}
\]

\mn
such that $\cD_{\theta,h}$ is a proper filter, \then \, 
$\bf{INV}^{\alpha,h}_{\kappa,\theta}(M,\varphi(\bar x,\bar y))$ 
almost has cardinality $\le \lambda$, which means: there are no
$\bold d^\zeta_i(i<\theta,\zeta<\lambda^+)$ such that:
\mn
\begin{enumerate}
\item[$(i)$]  for $\zeta<\lambda^+,\langle \bold d^\zeta_i:i<\theta\rangle/
\cD_{\theta,h} \in \bf{INV}^{\alpha,h}_{\kappa,\theta}(M,\varphi
(\bar x,\bar y))$,
\sn
\item[$(ii)$]  for $i<\theta,\zeta<\xi<\lambda^+$, we have 
$\bold d^\zeta_i \ne \bold d^\xi_i$.
\end{enumerate}
\end{lemma}

\begin{PROOF}{\ref{3.7}}
Straightforward.
\end{PROOF}
\bigskip

\centerline{$* \qquad * \qquad *$}
\bigskip

We now show that (for example for the case $\lambda$ regular) if $|I| \le
\lambda$ and $\inv^\alpha_\kappa(I)$ is well defined then there is a 
linear order $J$ such that:
if a model $M$ has a weakly $(\kappa,\varphi)$-skeleton like sequence
inside $M$ of order-type $J$ then $\inv^\alpha_\kappa(I)\in
\bf{INV}^\alpha_\kappa(M,\varphi)$.

Again, the proof splits into three cases depending on the cofinality of
$\lambda$. The following result provides a detail needed for the proof.
\begin{claim}
\label{3.7A}
Suppose that $\kappa$ is a regular cardinal and $\langle\bar a_t:t\in
J\rangle$ is a weakly $(\kappa,\varphi)$-skeleton like inside $M$ and
$I\subseteq J$. If for each $s\in J \setminus I$ either $\{t\in I:t<s\}$
or the inverse order on $\{t\in I:t>s\}$ has cofinality less than $\kappa$
(for example 1) \then \, $\langle\bar a_t:t\in I\rangle$ is weakly
$(\kappa,\varphi)$-skeleton like for $M$.
\end{claim}

\begin{PROOF}{\ref{3.7A}}
As usual let $\psi(\bar x,\bar y) =
\varphi(\bar y,\bar x)$.
We must show that for every $\bar a\in {}^{\ell g(\bar x)}M$ there is an
$I_{\bar a} \subseteq I$ with $|I_{\bar a}|<\kappa$ such that: if $s,t\in I$
and $\tp_{\qf}(s, I_{\bar a},I)=\tp_{\qf}(t,I_{\bar a},I)$ then

\[
M \models ``\varphi(\bar a_s,\bar a)\equiv \varphi (\bar a_t,\bar a)"
\text{ and } M \models `` \psi(\bar a_s,\bar a) \equiv \psi(\bar a_t,\bar a)".
\]

\mn
We know that there is such a set $J_{\bar a}$ for $J$ and $\bar a$ and for
each $s \in J_{\bar a}$ choose a set $X_s$ of $<\kappa$ elements of $I$ such
that $X_s$ tends to $s$, i.e., to the cut that $s$ induces in $I$ (either
from above or below).  (So if $s\in I$, $X_s=\{s\}$; otherwise use the
assumption).  Let $I_{\bar a}=\bigcup\limits_{s\in J_{\bar a}}X_s$;
as $ \kappa $ is regular, $|X_s| < \kappa$ for $s \in J_{\bar a}$ and
$|J_{\bar a}| < \kappa$ clearlly $I_{\bar a}$ has cardinality $<
\kappa$;  also trivially $J_{\bar a} \subseteq I$.

Now it is easy to see that if $t_1$ and $t_2 \in I$ have the same
quantifier free type over $I_{\bar a}$, then they have the same
quantifier free type over $J_{\bar a}$, and the claim follows.
\end{PROOF}

\begin{lemma}
\label{3.8}
Assume $\ell g(\bar x) = \ell g(\bar y) <\aleph_0$ and  
$\varphi=\varphi(\bar x,\bar y)$.

\noindent
1) Let $\lambda>\aleph_0$ be regular. If $I$ is a linear order of
cardinality  $\le \lambda$, and $\inv^\alpha_\kappa(I)$ is well defined,
\then \, for some linear order $J$ of cardinality $\lambda$ the following
holds:
\mn
\begin{enumerate}
\item[$(*)$]  if $M$ is a model of cardinality $\lambda$,
$\bar a_s \in {}^{\ell g(x)} M,\langle \bar a_s:s \in J \rangle$ 
is weakly $(\kappa,\varphi(\bar x,\bar y))$-skeleton like inside 
$M$ (hence $\varphi(\bar x,\bar y)$ is asymmetric), \then \, 
$\inv^\alpha_\kappa(I) \in \bf{INV}^\alpha_\kappa(M,\varphi(\bar x,\bar y))$.
\end{enumerate}
\mn
2) Let $\lambda$ be singular, $\theta=\cf(\lambda)>\kappa$, $\lambda=
\sum\limits_{i<\theta}\lambda_i$, where the sequence $\langle\lambda_i:i<
\theta\rangle$ is increasing continuous. Suppose that for $i<\theta$, $I_i$
is a linear order of cofinality $>\lambda_i$ and cardinality $\le \lambda$
such that $\inv^\alpha_\kappa(I_i)$ is well defined. \Then\ for
some linear order $J$ of cardinality $\lambda$ the following holds:
\mn
\begin{enumerate}
\item[$(**)$]  if $M$ is a model of cardinality $\lambda$,
$\bar a_s \in {}^{\ell g(x)} M$ for $s \in J$,
$\langle \bar a_s:s\in J\rangle$ is weakly $(\kappa,\varphi(\bar x,\bar
y))$-skeleton inside $M$, (so $\varphi(\bar x,\bar y)$ asymmetric), 
\then \, $\langle\inv^\alpha_\kappa(I_i):i<\theta\rangle/\cD_{\theta,
\kappa}$ belongs to $\bf{INV} ^\alpha_\kappa(M,\varphi(\bar x,\bar y))$.
\end{enumerate}
\mn
3) Let $\lambda$ be singular, $\theta,\kappa$ be regular, $\lambda>
\theta>(\cf(\lambda)+\kappa)$, $\lambda=\sum\limits_{i<\cf(\lambda)}
\lambda_i$, $\lambda_i$ increasing continuous. If, for $i<\theta$, $I_i$ is
a linear order such that $\inv^\alpha_\kappa(I_i)$ is well defined,
\then \, for some linear order $J$ of cardinality $\lambda$ the following
holds:
\mn
\begin{enumerate}
\item[$(***)$]  if $M$ is a model of cardinality $\lambda$,
$\bar a_s \in {}^{\ell g(x)} M$ for $s \in J$,
$\langle \bar a_s:s\in J\rangle$ is weakly $(\kappa,\varphi(\bar x,\bar
y))$-skeleton like inside $M$, (so $\varphi(\bar x,\bar y)$) asymmetric), $h$
is a function from a stationary subset of $\{\delta<\theta:\cf(\delta)\ge
\kappa\}$ with range a set of regular cardinals $<\lambda$ but $>\theta$
such that $\cf(I_i) \ge h(i)$ and $\cD_{\theta,h}$ is a proper filter
\then \, $\langle\inv^\alpha_\kappa(I_i):i<\theta\rangle/\cD_{\theta,
h}$ belongs to $\bf{INV}^{\alpha,h}_{\kappa,\theta}(M,\varphi(\bar x,\bar y))$.
\end{enumerate}
\end{lemma}

\begin{PROOF}{\ref{3.8}}
1 We must choose a linear order $J$ of cardinality $\lambda$ such that:
if $J$ indexes a weakly $(\kappa,\varphi(\bar x,\bar y))$-skeleton like
sequence inside $M$, a model of cardinality $\lambda$, then

\[
\inv^\alpha_\kappa(I) \in \bf{INV}^\alpha_\kappa(M,\varphi(\bar x,\bar y)).
\]

\mn
For this, for any continuous increasing decomposition $\bar A$ of $|M|$, we
must find a sequence $\bar c\in M$ and an ordinal $\delta$ with

\[
\INV^\alpha_\kappa(\bar c,A_\delta,M,\varphi(\bar x,\bar y))=
\inv^\alpha_\kappa(I).
\]

\mn
To obtain $\bar c$, we shall use a function from $\lambda$ to
$J$. Let $I_\alpha$ for $\alpha<\lambda$ be 
pairwise disjoint linear orders isomorphic to $I$.

Let $J=\sum\limits_{\alpha<\lambda}I^*_\alpha$ (where $I^*$ means we
use the inverse of $I$ as an ordered set). Suppose $\langle\bar a_s:s\in J
\rangle$ is weakly $(\kappa,\varphi(\bar x,\bar y))$-skeleton like inside
$M$, (hence $\varphi(\bar x,\bar y))$ is asymmetric), $M$ has cardinality
$\lambda$. For $\alpha<\lambda$ let $s(\alpha)\in I_\alpha$ and let
$\langle A_\alpha:\alpha<\lambda\rangle$ be an increasing continuous
sequence such that $M = \bigcup\limits_{\alpha<\lambda} A_\alpha$,
$|A_\alpha|< \lambda$. By the definition of weak $(\kappa,\varphi(\bar x,
\bar y))$-skeleton like (Definition \ref{3.1}(1)), for every
$\bar{a}\in {}^{\ell g(\bar x)} M$, here is a subset 
$J_{\bar a}$ of $J$ of cardinaltiy $<\kappa$ such that:
if $s,t\in J\setminus J_{\bar a} $ induces
the same Dedekind cut on $J_{\bar a}$,
then $M\models ``\varphi [\bar a_s,\bar a]\equiv  \varphi[\bar
a_t,\bar a]"$ and $M \models  ``\varphi [\bar a,\bar a_s]\equiv
\varphi[\bar a,\bar a_t]"$.

\noindent 
Since $\lambda$ is regular, for some closed unbounded subset
$\cC^*$ of $\lambda$, for every $\delta\in \cC$ we have:
\mn
\begin{enumerate}
\item[$(*)$]
\begin{enumerate}
\item[$(i)$]  $\bar a_{s(\alpha)} \in {}^{\ell g(\bar x}(A_\delta)$
 for $\alpha<\delta$,
\sn
\item[$(ii)$]  $J_{\bar a}\subseteq \sum\limits_{\beta<\delta}I^*_\beta$ for
$\bar a \in {}^{\ell g(\bar x}(A_\delta)$.
\end{enumerate}
\end{enumerate}
\mn
So it is enough to prove that for any $\delta\in \cC^*$ of cofinality
$ \le  \kappa $ we have 

\[
\inv^\alpha_\kappa(I)=\INV^\alpha_\kappa
(\bar a_{s(\delta)},A_\delta,M,\varphi(\bar x,\bar y)).
\]

\mn
Let $\cC \subseteq \delta$ be closed unbound of order types
$\cf(\delta)$.  It is easy to see that $\langle\bar a_s:s \in 
I_\delta\rangle$ does $\kappa$-represents $(\bar a_{s(\delta)},
A_\delta,M,\varphi(\bar x,\bar y))$ as: the required $\theta$ and $J$
in Definition \ref{3.5}$(\alpha)$ are $\cf(\delta)$ and 
$\langle\bar a_{s(\beta)}:\beta\in \cC \rangle$, and now use 
claim \ref{3.7A} with $J,\{s(\beta):\beta\in \cC\} \cup I^*_\delta$
here standing for $J,I$ there. 

So (see Definition \ref{3.5}$(\gamma)$) it is enough to 
show that $(\bar a_{s(\delta)},A_\delta,M,\varphi(\bar x,\bar y))$ 
has a $(\kappa,\alpha)$-invariant. Now in Definition
\ref{3.5}$(\beta)$,  part (ii) is obvious by the above; so it
remains to prove (i).

Let $\theta=:\cf(\delta)$.
So assume that for $\ell=1,2$,

\[
\langle\bar a^\ell_s:s\in I^\ell\rangle \text{ weakly } (\kappa,\theta)
\text{-represents }(\bar a_{s(\delta)},A_\delta,M,\varphi(\bar x,\bar y)).
\]

\mn
Let $J^\ell$, $\langle a^\ell_t:t\in J^\ell\rangle $ exemplify this (so each
$\bar a^\ell_t$ belongs to $A_\delta$) and let $J^\ast_\ell=J^\ell+
(I^\ell)^*$ and assume $\inv^\alpha_\kappa(I^\ell)$ are well defined.
We have to prove that $\inv^\alpha_\kappa(I^1)=\inv^\alpha_\kappa(I^2)$.
This follows by \ref{3.8A}(2) below.
\end{PROOF}

\begin{fact}
\label{3.8A}
1) Suppose $\langle\bar a_s:s\in J+I^\ast\rangle$ is weakly $(\kappa,
\varphi(\bar x,\bar y) )$-skeleton like inside $M$ and both $J$ 
and $I$ have cofinality $\ge \kappa$. \Then \, for every $\bar b \in
M$ there exist $s_0 \in J$ and $s_1\in I^*$ such that if 
$s_0<t_\ell<s_1$ (in $J+I^\ast)$ for $\ell=0,1$, \then
$M \models ``\varphi(\bar a_{t_0},\bar b)\equiv 
\varphi(\bar a_{t_1},\bar b)",M \models ``\psi(\bar a_{t_0},\bar b)\equiv
\psi(\bar a_{t_1},\bar b)"$.

\noindent
2) Suppose that, for $\ell=1,2$, $\langle \bar{a}^\ell_s:s\in I^\ell
\rangle$ does $(\kappa,\theta)$-represent $(\bar{c},A,M,\varphi(\bar{x},
\bar{y}))$ and $\langle\bar{a}^\ell_s:s\in J^\ell\rangle$ witnesses this.
\Then\ $\inv^\alpha_\kappa(I^1)=\inv^\alpha_\kappa(I^2)$.
\end{fact}

\begin{PROOF}{\ref{3.8A}}
1) Easy.

\noindent
2)  As we can replace $I^\ell$ by any end segment, \wilog \,
\mn
\begin{enumerate}
\item[$(*)$]   for $\ell=1,2$ for every $t\in I^\ell, \bar a_t$ realizes
$\tp_{\{\varphi(\bar x,\bar y),\psi(\bar x,\bar y)\}} (\bar c,A,M)$.
\end{enumerate}
\mn
We shall use Lemma \ref{3.3} (with $I^1,I^2$ here standing for
$I,J$ there and $\psi$ for $\varphi$). Conditions (b),(c) from \ref{3.3}
are met trivially, for (b) using \ref{3.2B} and by similar arguments 
in condition (a) it is enough to prove clause $(\alpha)$.

Let us prove (a)$(\alpha)$ from \ref{3.3}. So suppose it fails, i.e., $s\in
I^1$ but for arbitrarily large $t\in I^2$, $M\models\neg\varphi[\bar{a}^1_s,
\bar{a}^2_t]$.

Since $\langle\bar a^2_t:t\in J^2+(I^2)^*\rangle$ is weakly
$(\kappa,\varphi)$-skeleton like inside $M$, the preceding Fact
\ref{3.8A}(1) yields that for arbitrarily large $t\in J^2$, $M\models
\neg\varphi[\bar a^1_s,\bar a^2_t]$. Since $\bar a^1_s$ and $\bar c$
realize the same $\{\varphi,\psi\}$-type over $A_\delta$
(see definition \ref{3.5}$(\alpha)$ and (*) above),
and as $\bar a^2_t\subseteq A_\delta$
for $t\in J^2$, this implies $M\models\neg\varphi [\bar c, \bar
a^2_t]$, so this holds for arbitrarily large $t\in J^2$. Choose such $t_0\in
J^2$, this quickly contradicts the choice of $J^2$ and $I^2$. For, it
implies that for every $t\in I^2$ (as $\bar c,\bar a^2_t$
realize the same $\{\varphi,\psi\}$-type over $A_\delta$) we have

\[
M \models \neg\varphi [\bar a^2_t,\bar a^2_{t_0}],
\]

\mn
which is impossible as $\langle \bar{a}_s:s\in J^2+(I^2)^\ast\rangle$ is
weakly $(\kappa,\varphi)$-skeleton like (see Definition \ref{3.1}(3) the
last phrase).
\bigskip

\noindent
\underline{Continuing the proof of \ref{3.8}(2),(3)}: \, Left to the
reader (or see the proof of case (d) and formulation of case (e) in Theorem
\ref{3.11}). Take $J=\sum\limits_{i<\theta} (I_i)^*$ where $I_i\cong I$
are pairwise disjoint. 
\end{PROOF}

\begin{theorem}
\label{3.9}
Suppose that $\lambda>\kappa$, $K_\lambda$ is a family of $\tau$-models,
each of cardinality $\lambda,\varphi(\bar x,\bar y)$ is an asymmetric formula
with vocabulary $\subseteq\tau$, and $\ell g(\bar x) = 
\ell g(\bar y)<\aleph_0$.
Further, suppose that for every linear order $J$ of cardinality $\lambda$
there are $M\in K_\lambda$ and $\bar a_s\in M$ for $s\in J$ such that
$\langle \bar a_s:s\in J\rangle$ is weakly $(\kappa,\varphi(\bar x,
\bar y))$-skeleton like in $M$.

\Then \,, in $K_\lambda$, there are $2^\lambda$ pairwise non-isomorphic models.
\end{theorem}

\begin{PROOF}{\ref{3.9}}
First let $\lambda>\aleph_0$ be regular.

By \ref{3.4}(1) there are linear order $I_\zeta$ (for $\zeta<2^\lambda)$
each of cardinality $\lambda$, such that $\inv^1_\kappa(I_\zeta)$ are 
well defined and distinct. Let $J_\zeta$ relate to $I_\zeta$ as guarantee by
\ref{3.8}(1). Let $M_\zeta\in K_\lambda$ be such that there are $\bar
a^\zeta_s\in M_\zeta$ for $s\in J_\zeta$ such that $\langle\bar a^\zeta_s:
s\in J_\zeta\rangle$ is weakly $(\kappa,\varphi(\bar x,\bar y))$-skeleton
like inside $M_\zeta$ (exists by assumption). By \ref{3.8}(1), that is our
choice of $J_\zeta$, we have

\[
\inv^1_\kappa(I_\zeta)\in\bf{INV}^1_\kappa(M_\zeta,\varphi(\bar x,\bar y)).
\]

\mn
Clearly,

\[
M_\zeta \cong M_\xi \Rightarrow \bf{INV}^1_\kappa(M_\zeta,
\varphi(\bar x,\bar y))=\bf{INV}^1_\kappa(M_\xi,\varphi(\bar x,\bar y)),
\]

\mn
and hence

\[
M_\zeta\cong M_\xi\quad\Rightarrow\quad\inv^1_\kappa(I_\zeta)\in
\bf{INV}^1_\kappa(M_\xi,\varphi(\bar x,\bar y)).
\]

\mn
So if for some $\xi<2^\lambda$, the number of $\zeta<2^\lambda$ for which
$M_\zeta\cong M_\xi$ is $>\lambda$, then $\bf{INV}^1_\kappa(M_\xi,\varphi(\bar
x,\bar y))$ has cardinality $>\lambda$ (remember $\inv^1_\kappa(I_\zeta)$ were
pairwise distinct for $\zeta <2^\lambda$). But this contradicts
\ref{3.7}(1). 

So

\[
\{(\zeta,\xi):\zeta,\xi<2^\lambda\mbox{ and }M_\zeta\cong M_\xi\},
\]

\mn
which is an equivalence relation on $2^\lambda$, satisfies: each equivalence
class has cardinality  $\le \lambda$. Hence there are $2^\lambda$ equivalence
classes and we finish.

For $\lambda$ singular the proof is similar. If $\cf(\lambda)>\kappa$, we
can choose $\theta = \cf(\lambda)$ and use $\INV^2_{\kappa,\theta}$,
\ref{3.4}(2), \ref{3.8}(2), \ref{3.7}(2) instead of $\bf{INV}^1_{\kappa,
\theta}$, \ref{3.4}(1), \ref{3.8}(1), \ref{3.7}(1) respectively.

If $\cf(\lambda)\le\kappa$, let $\theta=\kappa^+$ so $\lambda>\theta>
\kappa+\cf(\lambda)$. Hence we can find a mapping

\[
h:\{\delta<\theta:\cf(\delta)\geq\kappa\}\longrightarrow\{\mu:\mu=\cf(\mu)<
\lambda\}
\]

\mn
such that for each $\mu=\cf(\mu)<\lambda$ the set

\[
\{\delta<\theta:\cf(\delta)\geq\kappa\mbox{ and }h(\delta)% 10.05.31 =
\ge \mu\}
\]

\mn
is stationary. Now we can use $\bf{INV}^{2,h}_{\kappa,\theta}$, \ref{3.4}(2),
\ref{3.8}(3), \ref{3.7}(3) instead $\bf{INV}^1_\kappa$, \ref{3.4}(1),
\ref{3.8}(1), \ref{3.7}(1) respectively.

Alternatively, for singular $\lambda$ see the proof of \ref{3c.16} and
\ref{3.11} case (d) below. 
\end{PROOF}

\begin{conclusion}
\label{3.10}
1) If $T_1$ is a first order $T\subseteq T_1$, $T$ is unstable and
complete, $\lambda \ge |T_1|+\aleph_1$, \then \, there are
$2^\lambda$ pairwise non-isomorphic models of $T$ of cardinality 
$\lambda$ which are reducts of models of $T_1$.

\noindent
2) If $T\subseteq T_1$ are as above, $\lambda \ge |T_1|+\kappa^+$,
$\lambda=\lambda^{<\kappa}$, $\kappa$ is regular, \then \, there are
$2^\lambda$ pairwise non-isomorphic models of $T$ of cardinalty $\lambda$ which
are reducts of models $M^1_i$ of $T_1$ such that $M_i,M^1_i$ are
$\kappa$-compact and $\kappa$--homogeneous. [Really we can get strongly
homogeneous; see \cite[\S1]{Sh:363}].

\noindent
3) Assume that $\psi\in \bbL_{\kappa^+,\omega}(\tau_1)$, $\tau\subseteq
\tau^1$, $\psi$ has the order property for $\bbL_{\kappa^+,\omega}(\tau)$
[i.e., for some formula $\varphi(\bar x,\bar y)\in \bbL_{\kappa^+,\omega}
(\tau)$ for arbitrarily large $\mu$ there is a model $M$ of $\psi$ and $\bar
a_i\in M$ for $i<\mu$ such that

\[
M\models\varphi[\bar a_i,\bar a_j] \text{ iff } [i<j \text{ and }
\ell g(\bar x) = \ell g(\bar y)<\aleph_0].
\]

\mn
\Then \, for $\lambda\ge\kappa+\aleph_1$, $\psi$ has $2^\lambda$
models of cardinality $\lambda$, with pairwise non-isomorphic $\tau$-reducts.
\end{conclusion}

\begin{PROOF}{\ref{3.10}}
1) Let $\varphi=\varphi(\bar x,\bar y)$ be a first order formula
exemplifying ``$T$ is unstable'' (see Definition \ref{1.2}). By \ref{1.8}(1)
there is a template $\Phi$ proper for linear orders such that $|\tau_\Phi|=
|\tau_1|$ and for any linear order $I$, $EM(I,\Phi)$ is a model of $T_1$
satisfying $\varphi[\bar a_s,\bar a_t]$ \underline{if and only if} $I 
\vDash s<t$. Clearly
$\EM_{\tau(T_1)}(I,\Phi)$ has cardinality $\ge |I|$ but $\le |\tau_\Phi|+
|I|+\aleph_0$. So for every $\lambda \ge |T_1|+\aleph_0 =|\tau_\Phi|+
\aleph_0$ and linear order $I$ of cardinality $\lambda$ the model 
$M= EM_\tau(I,\Phi)$ is a $\tau$--model, a reduct of a model 
of $T_1$, hence $M$ is a model of $T$ of cardinality
exactly $\lambda$, and by \ref{3.4}(4) the sequence $\langle\bar a_t: t\in
I\rangle$ is weakly $\kappa$-skeleton like. So we have the assumption of
\ref{3.9}, hence its conclusion as required.

\noindent
2) By \cite[Ch.VII 3.1]{Sh:c}, or case II of the proof of Theorem
3.2 (there) we have the assumption of \ref{3.9};
but \cite[\S1]{Sh:363} supersedes upon this.

\noindent 
3) See \ref{1.11}(3) and Definition \ref{1.10} why the assumption
of \ref{3.9} holds. 
\end{PROOF}

\begin{remark}
\label{3.10A}
Also \ref{1.15} is a similar result.
\end{remark}
\bigskip

\centerline {$* \qquad * \qquad *$}
\bigskip

Now we turn our attention to the case in which the sequences on which
$\varphi(\bar{x},\bar{y})$ speaks are infinite.

\begin{theorem}
\label{3.11}
Suppose $\partial<\kappa<\lambda$ are cardinals, $\kappa$ regular. Assume $K$
is a class of $\tau$-models, $\varphi=\varphi(\bar x,\bar y)$ is a formula
with vocabulary $\subseteq\tau$, and $\partial=\lg(\bar x)=\lg(\bar
y)$, and
\mn
\begin{enumerate}
\item[$(*)$]  $K=K_\lambda$ and for every linear order $I$ of cardinality
$\lambda$ there are $M_I\in K_\lambda$ and a sequence $\langle\bar a_t: t\in
I\rangle$ which is weakly $(\kappa,\varphi(\bar x,\bar y))$-skeleton like
inside $M_I$.
\end{enumerate}
\mn
We can conclude that $\dot{\bbI}(K)=2^\lambda$ \Iff \, at least one of the
following conditions holds:
\mn
\begin{enumerate}
\item[$(a)$]  $\lambda=\lambda^\partial$
\sn
\item[$(b)$]  $\lambda^\kappa<2^\lambda$
\sn
\item[$(c)$]  We replace the assumption $(*)$ by:
\sn
\begin{enumerate}
\item[$(*)_0$] $K=K_\lambda$,
\sn
\item[$(*)_1$]  $\lambda^\partial<2^\lambda$, $\cf(\lambda)>\partial$,
\sn
\item[$(*)_2$] for every linear order $J$ of cardinality $\lambda$ there
are $M_J\in K_\lambda$ and a weakly $(\kappa,<\lambda,
\varphi(\bar x,\bar y))$-skeleton like inside $M_J$ sequence $\langle\bar
a_s:s\in J\rangle$ (where $\bar a_s\in{}^\partial |M_J|$), see Definition
\ref{3.12} below.
\end{enumerate}
\sn
\item[$(d)$] We replace the assumption $(*)$ by:
for some $\lambda(0) \le \lambda(1) \le \lambda\le\lambda(3)<2^\lambda$,
$\mu(0) \le\ mu(1) \le 2^\lambda$ with $\lambda(1)$ and $\mu(1)$ are
regular, we have:
\sn
\begin{enumerate}
\item[$(*)_0$]  $K=K_{\lambda(3)}$,
\sn
\item[$(*)_1$] $\lambda^\partial<2^\lambda$,
\sn
\item[$(*)_2$] for every linear order $J$ of cardinality $\lambda$ there
is $M_J\in K_{\lambda(3)}$ (of cardinality $\lambda(3)$) and 
$\langle\bar a_s:s\in J\rangle$ (where $\bar a_s \in {}^\partial 
|M_J|$) which is weakly $(\kappa,\lambda(0),<\lambda(1),
\varphi(\bar x,\bar y))$-skeleton like inside $M_J$
(see Definition \ref{3.12} below),
\sn
\item[$(*)_{3,\mu(0),\lambda(0)}$] for $J\in K^{\oor}_\lambda
(=(K_{\oor})_{\lambda})$ and a set $A\subseteq M_J$ ($M_J$ is from $(*)_2$)
if $|A|<\lambda(0)$ then:
\sn
\begin{enumerate}
\item[$(i)$]  $\mu(0) > |\bbS^\partial_{\{\varphi,\psi\}}(A,M_J)|$, or at least
\sn
\item[$(ii)$]  $\mu(0)>|\big\{\Av_{\{\varphi,\psi\}}(\langle\bar b_i:i<\kappa
\rangle,A,M_J:\bar b_i\in A$ for $i<\kappa$, the average is well defined and
is realized in $M\big\}|$, where
\begin{equation*}
\begin{array}{clcr}
\Av_\Delta(\langle b_i:i<\kappa\rangle,A,M_J) :=
\{\varphi(\bar x,\bar a)^{\bold t}:&\varphi(\bar x,\bar y)\in\Delta,
{\bold t}\text{ a truth value,} \\
  &\bar a\in A\text{ and } \text{ for all but a bounded set of }
i<\kappa,\ M_J\models \varphi[\bar b_i,\bar a]^{\bold t}\},
\end{array}
\end{equation*}
\end{enumerate}
\sn
\item[$(*)_{4,\lambda,\mu(1),\mu(0),\lambda(0)}$] if $\dot{\bold I}_i\subseteq
{}^\partial\lambda(3)$ and $|\dot{\bold I}_i|=\lambda$ for $i<\mu(1)$, \then\
for some $B\subseteq\lambda(3)$ we have:
\[
|B|<\lambda(0) \text{  and } |\{i:|\\dot{\bold I}_i \cap {}^\partial
B| \ge \kappa\}|\ge \mu(0).
\]
\end{enumerate}
\sn
\item[$(e)$]   We replace assumption $(*)$ by:
for some $\lambda_{0,\epsilon}\le\lambda_{1,\epsilon} \le \lambda \le
\lambda_3,\mu_{0,\epsilon} \le \mu_1 \le 2^\lambda$, for $\epsilon<
\epsilon(*),\mu_1$ is regular and:
\sn
\begin{enumerate}
\item[$(*)_0$] $K= K_{\lambda_3}$,
\sn
\item[$(*)_1$] $\lambda^\partial<2^\lambda$,
\sn
\item[$(*)_2$] for every linear order $J$ of cardinality $\lambda$ there
is $M_J\in K_{\lambda(3)}$ and $\langle\bar a_s:s\in J\rangle$ (where $\bar
a_s\in {}^\partial |M_J|)$ which for each $\epsilon<\epsilon(*)$ is
weakly $(\kappa_1,<\lambda_{0,i},<\lambda_{1,i},
\varphi(\bar x,\bar y))$-skeleton like inside $M_J$,
\sn
\item[$(*)_{3,\mu_{0,\epsilon},\lambda_{0,\epsilon}}$] if
$\epsilon<\epsilon(*)$ and $J\in K^{\oor}_\lambda(=(K_{\oor})_\lambda)$
and a set $A\subseteq M_J$ ($M_J$ is from $(*)_2)$ 
if $|A|<\lambda_{0,\epsilon}$ then:
\sn
\begin{enumerate}
\item[$(i)$]  $\mu _{0,\epsilon}>|{\bf S}^\partial_{\{\varphi,\psi\}}
(A,M_J)|$ or at least
\sn
\item[$(ii)$]  $\mu_{0,\epsilon}>|\big\{\Av_{\{\varphi,\psi\}}(\langle\bar
b_i: i<\kappa\rangle,A,M_J):\bar b_i\in A$ for $i<\kappa$, the average is
well defined and is realized in $M\big\}|$, where

\begin{equation*}
\begin{array}{clcr}
\Av_\Delta(\langle b_i:i<\kappa\rangle,A,M_J):= 
\quad\{\varphi(\bar x,\bar a)^{\bf t}:&\varphi(\bar x,\bar y)\in\Delta,\
{\bold  t} \text{ a truth value,} \\
  &\bar a \in A \text{ and } \text{for all but a bounded set of }
i<\kappa,\ M_J\models \varphi[\bar b_i,\bar a]^{\bold t}\},
\end{array}
\end{equation*}
\end{enumerate}
\sn
\item[$(*)_4$]  there are $h_\alpha:\lambda\longrightarrow\{\theta:\theta$
regular, $\kappa \le\ theta\le\lambda\}$ for $\alpha<2^\lambda$ such that:
if $S\subseteq 2^\lambda$, $|S|\ge\mu(1)$ and $f_\alpha:\lambda\longrightarrow
{}^\partial(\lambda_3)$ for $\alpha\in S$, \then \, we can find
$\epsilon<\epsilon(*),B \subseteq\lambda_3$ satisfying:
$|B|<\lambda_{0,\epsilon}$ and the set
$	
\{\alpha:\text{ the closure of }\{\zeta<\lambda:f_\alpha(\zeta)\subseteq
B\}$ has a member $\delta$ of cofinality $\kappa$ such that 
$h_\alpha(\delta)\ge\lambda_{1,\epsilon}\}$
has $\ge \mu_{0,\epsilon}$ members. [Note: $\cf(\delta)=\kappa'\ge
\kappa$ can be allowed if
$(*)_{3,\mu_{0,\epsilon},\lambda_{0,\epsilon}}$ is changed accordingly].
\end{enumerate}
\sn
\item[$(f)$]  For some $\mu<\lambda$, there is a linear order of
cardinality $\mu$ with $\ge \lambda$ Dedekind cuts each with 
upper and lower cofinality $\ge \kappa$ and
$2^{\mu+\partial}<2^\lambda$.
\sn
\item[$(g)$]  there is $\cP \subseteq [\lambda^\partial]^\kappa$ of
cardinality $<2^\lambda$ such that every $X\subseteq \lambda^\partial$
of cardinality $\lambda$ contains at least one of them (and $(*)$);
(can use similar considerations in other places).
\end{enumerate}
\end{theorem}

\begin{definition}
\label{3.12}
We say $\langle\bar a_s:s\in I\rangle$ is weakly $(\kappa,\mu,<\lambda,
\varphi(\bar x,\bar y))$-skeleton like inside $M$;
if $\mu=\lambda$ we may omit $\mu$; \Iff \,:
\mn
\begin{enumerate}
\item[$(i)$]  for $s,t\in I$ we have
\[
M \models \varphi[\bar a_s,\bar a_t] \text{\underline{ if and only if }} 
I\models s<t,
\]
\sn
\item[$(ii)$]  for every $\bar{c}\in {}^{\ell g(\bar a_s)}M$ for some
$J\subseteq I$, $|J|<\kappa$ and $(*)$
of \ref{3.1}(1) holds, and
\sn
\item[$(iii)$]  moreover, for each $A\subseteq M$, $|A|<\mu$,
there is $J\subseteq I$, $|J|<\lambda$ such that for every $\bar c \in
{}^{\ell g(\bar x)} A$,the statement $(*)$ of \ref{3.1} holds for $J$.
\end{enumerate}
\end{definition}

\begin{PROOF}{\ref{3.12}}
\medskip

\noindent
\underline{Case (a)}:

In Definition \ref{3.5} we can replace
$A$ by ${\bf \dot{J}}$, a set of sequences of length $\partial$ from $M$,
which means that clause (i) in $(\alpha)$ of \ref{3.5} now becomes (i)' for
every large enough $t\in I$, for every $I \in \dot{\bold J}$ we have $M\models
\varphi [\bar a,\bar b]=\varphi [\bar{a}_t,\bar b]$ and $M\models \psi
[\bar c,\bar b]\equiv \varphi[\bar{a}_t,\bar b]$.

Thus in Definition \ref{3.6}, replace 
$\langle A_i:i<\lambda\rangle$ by $\langle
\dot{\bold J}_i: i< \cf(\lambda))\rangle,{}^\partial |M|=
\bigcup\limits_i \dot{\bold J}_i,|\dot{\bold J}_i| < \lambda$, 
$\dot{\bold J}_i$ increasing continuous. No further changes in
\ref{3.1}-\ref{3.9} is needed.

Alternatively, we can define $N=F_\partial(M)$ as the model with universe $|M|
\cup {}^\partial |M|$, assuming of course $|M|$ is disjoint to
${}^\partial|M|$,

\[
\tau(N)=\tau(M)\cup\{F_i:i<\partial\},
\]

\[
R^N=R^M \text{ for } R \in \tau(M),
\]

\[
G^N(x_1,\ldots,x_n) = \begin{cases}
G^M(x_1,\ldots,x_n) \quad &\text{ if } x_1,\ldots,x_n\in |M|,\\
x_1 \quad &\text{ otherwise }
\end{cases}.
\]

\mn
for function symbol $G\in\tau(M)$ which has $n$-places and

\[
F^N_i(x) = \begin{cases}
x(i) \quad&\text{ if } x \in {}^\partial M,\\
x \quad &\text{ if } x \in M
\end{cases}
\]

\mn
for $i<\partial$, so $F_i$ is a new, unary function symbol for
$i<\partial$.

Note that [$M_1 \cong M_2$ if and only if $F_\partial(M_1)\cong F_\partial
(M_2)$], and $\|F_\partial(M)\|=\|M\|^\partial$, etc. So we can apply \ref{3.9}
to the class $\{F_\partial(M):M\in K_\lambda\}$ and we can get the desired
conclusion.
\bigskip

\noindent
\underline{Case (b)}:  We use weakly 
$(\kappa,\varphi(\bar x,\bar y))$-skeleton like sequences $\langle\bar a_s:s\in
\kappa +(I_\zeta)^*\rangle$ in $M_\zeta\in K_\lambda$ for
$\zeta<2^\lambda$, with $\langle\inv^2_\kappa(I_\zeta):\zeta<2^\lambda
\rangle$ pairwise distinct, and count the number of models $(M_\zeta,\langle
\bar a_s:s\in \kappa\rangle)$ up to isomorphism. Then ``forget the $\bar
a_s$, $s\in \kappa$", i.e., use \ref{3.13} below.
\bigskip

\noindent
\underline{Case (c)}:  We revise \ref{3.5}--\ref{3.10}; we use
this opportunity to present another reasonable choice in clause $(\alpha)$
of \ref{3.5}.

\noindent
\underline{Change 1}:  In \ref{3.5}($\alpha$) we replace
(i), (ii) by 
\mn
\begin{enumerate}
\item[$(i)'$]  for every formula $\vartheta(\bar{x},\bar{d})\in\tp_{\{\varphi(
\bar{x},\bar{y}),\psi(\bar{x},\bar{y})}(\bar{c},A,M)$, for every large
enough $t\in I$ we have $M\models\vartheta[\bar{c},\bar{d}]\equiv
\vartheta[\bar{a}_t,\bar{d}]$,
\sn
\item[$(ii)'$]  $\langle\bar{a}_s:s\in J+(J)^*\rangle$ is weakly $(\kappa,
\varphi(\bar{x},\bar{y}))$-skeleton like inside $M$,
\sn
\item[$(iii)'$] $\theta>\cf(J)$ (actually $\theta\neq\cf(J)$ would suffice,
but no real need)
(not actually needed, but natural).
\end{enumerate}
\mn
Of course, the meaning of Definition \ref{3.5}($\beta$)-($\delta$) changes,
and the reader can check that, e.g., the proof of the Fact is still valid.
\bigskip

\noindent
\underline{Change 2}:  In Definition \ref{3.6}(1), inside
the definition of $\bf{INV}^\alpha_\kappa$, we demand $\cf(\bold d) = \lambda$
recalling $\lambda$ is regular.
\bigskip

\noindent
\underline{Change 3}:  In Definition \ref{3.6}(2), inside
the definition of $\INV^\alpha_{\kappa, \theta}$ add $\cf(\bold d_\delta)> \cf
(\delta)$ (necessitate by change 1, actually $\cf(\bold d_\delta) \ne 
\cf(\delta)$ suffices).
\bigskip

\noindent
\underline{Change 4}: In Definition \ref{3.6}(3) demand
$\cf(\lambda)>\partial$.
\bigskip

\noindent
\underline{Change 5}:  In \ref{3.7}, in all cases the
``cardinality $\le \lambda$'' is replaced by ``cardinality $\le
\lambda^\partial$" and part (2) becomes like part (3).
\bigskip

\noindent
\underline{Change 6}: We replace ``$(\kappa,
\varphi(\bar x,\bar y))$-skeleton likeq" by $(\kappa,<\lambda,\varphi(\bar
x, \bar y))$-skeleton like.
In \ref{3.8}(3) add the demand $\cf(\lambda)>\partial$, $h(i)>\cf(i)$.
\bigskip

\noindent
\underline{Change 7}: Inside the proof of \ref{3.8}(1), now
not for every $\bar a\in {}^{\ell g(\bar x)}M$ we define $J_{\bar a}$, but for
every $A\subseteq M$ of cardinality $<\lambda$ we choose $J_A\subseteq J$,
$|J_A|<\lambda$ by Definition \ref{3.12}, and in $(*)(ii)$
in the proof there we demand

\[
(\forall\alpha<\delta)(\exists\beta<\delta)[\bigcup\limits_{s\in J_{A_\alpha}}
\bar a_s \subseteq A\beta].
\]

\mn
\underline{Change 8}:  In the proof of \ref{3.8}(2) let
$\langle I_i:i<\theta\rangle$ be as in the statement of \ref{3.8}(2), and
let $J=\sum\limits_{i<\theta} I^*_i$, and assume $\langle\bar a_s: s\in J
\rangle$ is $(\kappa,<\lambda,\varphi(\bar x,\bar y))$-skeleton like inside
$M\in K_\lambda$. So let $\langle A_i: i<\theta\rangle$ be a representation
of $M$, and for each $i<\theta$ let $J_{A_i}\subseteq J$,
$|J_{A_i}|<\lambda$ be as in Definition \ref{3.12}. 

Define

\begin{equation*}
\begin{array}{clcr}
\cC = \{\delta<\theta:&\delta \text{ is a limit ordinal such that for every }
\alpha<\delta \\
  & \text{ the cardinality of } J_{A_i} \text{ is } <\lambda_\delta\}.
\end{array}
\end{equation*}

\mn
So let $\delta\in C,\cf(\delta)\ge\kappa$. Recall that $\cf(I_\delta)
>\lambda_\delta$ so clearly we can find $s(\delta)\in I_\delta$ such that

\[
I_\delta\models s(\delta)\le s \Rightarrow s \notin
\bigcup\limits_{i<\delta} J_{A_i}.
\]

\mn
Now $(\bar c_{s(\delta)},A_\delta,M,\varphi(\bar x,\bar y))$ is as required.
\bigskip

\noindent
\underline{Change 9}:  In the proof of \ref{3.12}(3) let
$J=\sum\limits_{\alpha<\theta} I^*_\alpha$ and $M$, $\langle\bar a_s: s\in J
\rangle$, $\langle A_i:i<\cf(\lambda)\rangle$, $J_{A_i}\subseteq J$ be as
above, and let $s(\alpha)\in I_\alpha$. As $\cf(\lambda)>\partial$ by
$(*)_1$ of the assumption, for each
$s\in J$ for some $i(s)<\cf(\lambda)$ we have $\bar c\subseteq A_{i(s)}$,
but $\theta=\cf(\theta)>\cf(\lambda)$ hence for some $i(*)<\cf(\lambda)$ the
set $W=\{\alpha<\theta: i(\alpha)\le i(*)\}$ is unbounded in $\theta$. Let
$\cC = \{\delta<\theta:\delta=\sup(\delta \cap W)\}$. We can choose $\delta\in
\cC$ of cofinality $\geq\kappa$ such that $h(\delta)>|J_{A_{i(*)}}|$, and
continue as in the previous case.
\bigskip

\noindent
\underline{Change 10}: Proof of \ref{3.8A}(2) (necessitated by change 1)

We shall use Lemma \ref{3.3} (with  $I^1,I^2$ here standing for $I,J$ there
and $\psi$ for $\varphi$). Conditions (b), (c) from \ref{3.3} are met
trivially and  by similar arguments in condition (a) it is enough to prove
clause $(\alpha)$.

Let us prove (a)$(\alpha)$ from \ref{3.3}.
Let $I^\ell_* \subseteq I^\ell$ be unbounded of order type
$\cf(I^\ell)=\theta$ and let $J^\ell_\ast\subseteq J^\ell$ be
unbounded of order type $\cf(J^\ell)$, which is $\ne \theta$. Possibly
shrinking those sets the truth values of
$\varphi[\bar{a}^1_s,\bar{a}^2_t]$ when $s\in I^1_*,y \in J^2\wedge
(\exists t') (t'\in J^2_*$ and $t'<_{J^2} t)$ is constant.
We can continue as before.

Note that if $\cf(\lambda)>\kappa$ this follows from case (d). If $\lambda$
is regular, choose $\lambda(0)=\lambda(1)=\lambda(3)=\lambda$ and $\mu(0)=
\mu(1)=(\lambda^\partial)^+$ and now the assumptions hold. If $\lambda$ is
singular, let $\epsilon(*)=\cf(\lambda),\chi=(\cf(\lambda)+\kappa)^+
\le \lambda$, $\mu_0=\mu_{1,\epsilon}=(\lambda^\partial)^+$ and let
$\{(\lambda_{0,\epsilon},\lambda_{1,\epsilon}):\epsilon<
\epsilon(*)\}$ list $\{(\lambda^+_i,\lambda^+_j):i<j<
\cf(\lambda)\}$ and
choose $h_\lambda=h:\lambda\longrightarrow \{\theta:\theta$ regular, $\kappa
\leq\theta \le \lambda\}$ such that $\epsilon<\epsilon(*)=\cf(\lambda)$
implies $\{\delta<\chi:\cf(\delta)=\kappa$ and $h(\delta)=\epsilon\}$ is
stationary. Now we can apply case (e).
\bigskip

\noindent
\underline{Case (d)}: 
Let $\langle I_\alpha:\alpha< 2^\lambda\rangle$ be a sequence of linear
orders of cofinality $\cf(\lambda(1))=\lambda(1)$, each of cardinality
 $\lambda$, with pairwise distinct $\inv^2_\kappa(I_\alpha)$ if 
$\lambda$ is regular, $\inv^3_\kappa(I_\alpha)$ if $\lambda$ is 
singular exists by \ref{3.4}. Let $J_\alpha =
\sum\limits_{\zeta\le\lambda} I^*_{\alpha,\zeta}$, where 
$I_{\alpha,\zeta}$ are pairwise disjoint, $I_{\alpha,\zeta}
\cong I_\alpha$. Let $M_{J_\alpha}$ be a model as guaranteed in $(*)_2$
with $\langle\bar a_s:s\in J_\alpha\rangle$ as there. Suppose
$\{M_{J_\alpha}/{\cong}\;:\alpha<2^\lambda\}$ has cardinality
 $<2^\lambda$, then without loss of generality 
$M_{J_\alpha}=M_{J_0}$ for $\alpha<\mu(1)$ and
without loss of generality $M_{J_0}$ has universe $\lambda(3)$. Let
$s(\alpha,\zeta)\in I_{\alpha,\zeta}$, so

\[
\dot{\bold I}_\alpha := \{\bar a_{s(\alpha,\zeta)}:\zeta<\lambda\}
\]

\mn
is a subset of ${}^\partial(\lambda(3))$ of cardinality $\lambda$. By
$(*)_{4,\lambda,\mu(1),\mu(0),\lambda(0)}$ there is $B\subseteq
\lambda(3)$, $|B|<\lambda(0)$ such that

\[
S=:\{\alpha<\mu(1):|\dot{\bold I}_\alpha\cap {}^\partial B| \ge \kappa\}
\]

\mn
has cardinality $\ge \mu(0)$. Choose for each $\alpha\in S$ a set

\[
S_\alpha\subseteq\{\zeta:\bar a_{s(\alpha,\zeta)}\subseteq B\},
\]

\mn
which has order type $\kappa$, and let

\[
\delta_\alpha=:\sup(S_\alpha).
\]

\mn
Clearly $\delta_\alpha\leq\lambda$, hence
$I_{\alpha,\delta_\alpha}$ is well defined. For each 
$\alpha\in S$, as $\langle\bar{a}_s:s\in J_\alpha\rangle$
is ($\kappa,\lambda(0),<\lambda(1),\varphi(\bar{x},\bar{y})$)-skeleton like
and $|B|<\lambda(0)$, there is a subset $J_{\alpha,B}$ of $J_\alpha$ as in
Definition \ref{3.12}. But $I_{\alpha,\delta_\alpha}$ has
cofinality $\lambda(1)>|B|$, hence for all large enough $t\in I_{\alpha,
\delta_{\alpha}}$, the type $\tp_{\{\varphi,\psi\}}(\bar a_t,B,M_{J_0})$ is
the same; choose such $t_\alpha$. Clearly (for $\alpha\in S$)

\[
\tp_{\{\varphi,\psi\}}(\bar a_{t_\alpha},B,M_{J_0})=\Av_{\{\varphi,\psi\}}
(\langle\bar a_{s(\alpha,\zeta)}:\zeta\in S_\alpha\rangle,B,M_{J_0}),
\]

\mn
so by $(*)_{3,\mu(0),\lambda(0)}$ from the assumption of case (d)
\wilog \, for some $\alpha\neq\beta$ we
get the same type. But $I_\alpha,I_\beta$ have different (and well defined)
$\inv^2_\kappa$ (or $\inv^3_\kappa$), contradicting \ref{3.8A}(2).
\bigskip

\noindent
\underline{Case (e)}:

Similar proof (to (d)).
\bigskip

\noindent 
\underline{Case (f)}:

By \ref{3.13} below.
\bigskip

\noindent 
\underline{Case (g)}:

Similar to case (b).
\end{PROOF}

\begin{fact}
\label{3.13}
If $\tau_2=\tau_1\cup\{c_i:i\in I\}$, $c_i$ are individual constants,
$K_\ell$ is a class of $\tau_\ell$-models (for $\ell=1,2$), $M\in K_2\
\Rightarrow\ M \rest \tau_1\in K_1$, and $\mu = \dot{\bbI}(\lambda,K_2)>
\lambda^{|I|}$, \then \, $\dot{\bbI}(\lambda,K_1) \ge \mu$ 
(so if $\mu=2^{\lambda+|\tau_1|}$, equality holds).
\end{fact}

\begin{PROOF}{\ref{3.13}}
Straight (or see \cite[Ch.VIII,1.3]{Sh:a}).
\end{PROOF}

\noindent
In \ref{3.10}-\ref{3.11} above we do not get anything when $\lambda^\partial=
2^\lambda$, however if we assume that $M_J$ has a clearer structure , e.g.,
is an $\EM$-model, we can get better results as done below.

\begin{conclusion}
\label{3.14}
1) Suppose $\psi\in \bbL_{\chi^+,\omega}(\tau_1)$, $\tau\subseteq
\tau_1$, $\varphi(\bar x,\bar y)\in {\Bbb L}_{\chi^+,\omega}(\tau)$,
$\ell g(\bar x) = \ell g(\bar y)=\partial \le \chi$, and $\psi$ 
has the $\varphi (\bar x,\bar y)$-order property that is for
every $\mu$ for some model $M$ of
$\psi$ there are $\bar a_i\in {}^\partial M$ (for $i<\mu$) such that

\[
M\models\varphi[\bar a_i,\bar a_j]\quad\mbox{ iff }\quad i<j.
\]

\mn
\Then \, for every $\lambda$ such that $\lambda>\chi^\partial$ or
$\lambda >\chi$ and $2^\lambda>\lambda^\partial$, $\psi$ has $2^\lambda$ models
of cardinality $\lambda$ with pairwise non-isomorphic $\tau$-reducts.

\noindent
2) Suppose $\psi\in \bbL_{\chi^+,\omega}(\tau)$, $\varphi_\ell(\bar x,
\bar y)\in \bbL_{\chi^+,\omega}(\tau_\ell)$, for $\ell=1,2$, $\ell g(\bar x)=
\ell g(\bar y)=\partial$, $\tau_0=\tau_1\cap\tau_2=\tau_1\cap\tau=\tau_2\cap
\tau$, $\{\psi,\varphi_1(\bar x,\bar y),\varphi_2(\bar x,\bar y)\}$ has no
model and $\psi$ has the $(\varphi_1,\varphi_2)$-order property,
which means that
\mn
\begin{enumerate}
\item[$(*)$]  for every $\alpha$ there is a $\tau_0$-model $M$ and $\bar
a_\beta\in {}^\partial |M|$ for $\beta<\alpha$, such that:
if $\beta<\gamma<\alpha$ \then \, 
\sn
\begin{enumerate}
\item[$(i)$]  for some expansion $M'$ of $M$, $M'\models\varphi_1[\bar a_\beta,
\bar a_\gamma]$,
\sn
\item[$(ii)$] for some expansion $M'$ of $M$, $M'\models\varphi_2[\bar
a_\gamma, \bar a_\beta]$.
\end{enumerate}
\end{enumerate}
\mn
Let $\varphi(\bar x,\bar y)=(\exists\ldots,R,\ldots)_{R\in\tau_1\setminus
\tau_0}\varphi_1(\bar x,\bar y)$; it is a formula in the vocabulary $\tau_0$
(but of second order). 
\Then \,
\mn
\begin{enumerate}
\item[$(a)$]  for $\lambda$ such that $\lambda>\chi^\partial$ or 
$\lambda>\chi$ and $2^\lambda>\lambda^\partial,
\dot{\bbI}_\tau(\lambda,\psi)=2^\lambda$ i.e., there are
$2^\lambda$ non-isomorphic $\tau$-models of $\psi$ of cardinality $\lambda$,
in fact even their $\tau_0$--reducts are not isomorphic;
\sn
\item[$(b)$]  for $\lambda\ge \chi$ there are 
$\langle M_J:J\in (K_{\oor})_\lambda \rangle$, $M_J$ a model 
of $\psi$ of cardinality $\lambda$ with a weakly
$(\partial^+,\varphi)$-skeleton like $\langle\bar a_s:s\in J\rangle$, $\bar
a_s\in {}^\partial |M_J|$, fully represented in $\cM_{\chi,\aleph_0}$
and $\bar a_s=\bar \sigma(s)$ for some sequence $\bar \sigma$ of term of
$\tau_{\chi,\aleph_0}$ see \ref{2.2},  or even 
$\bar a_s=\langle F_{1,i}(s):i<\partial\rangle$.
\end{enumerate}
\end{conclusion}

\begin{PROOF}{\ref{3.14}}
1) Follows from (2), by taking $\varphi(\bar x,\bar y)=
\varphi_1(\bar x,\bar y)=\varphi_2(\bar y,\bar x)$.

\noindent
2) By \ref{1.11}(3), \ref{1.15} there is $\Phi$, proper for the
class of linear orders (see Definition \ref{1.6}) such that for every linear
order $I,\EM_\tau(I,\Phi)$ is a model of $\psi$ of cardinality 
$\chi+|I|$, for $t\in I,\bar a_t$ is a sequence of length 
$\partial$ of members of $\EM_\tau(I,\Phi)$, in fact is 
$\bar \sigma(t)$ for a fixed $\bar \sigma$, such that for $s,t\in I$:

\begin{equation*}
\begin{array}{clcr}
\EM_\tau(I,\Phi) \models \varphi[\bar a_s,\bar a_t] \quad &\text{ iff
} s < t \\
  &\text{ iff } \EM_ tau(I,\Phi)\vDash\neg(\exists\ldots,R,\ldots)_{R\in\tau_2
\setminus\tau_1}\varphi_2[\bar a_t,\bar a_t].
\end{array}
\end{equation*}

\mn
By \ref{3.1A} $\langle\bar a_s:s\in I\rangle$ is weakly $(\partial^+,
\varphi)$-skeleton like (see Definition \ref{3.1}). Clearly
$\EM_ tau(I,\Phi)$ is represented in $\cM_{\chi,\aleph_0}$. So the 
clause (b) of \ref{3.14}(2) holds. To prove clause (a) we can use 
\ref{3c.16}, Case A (as $\theta=\aleph_0)$ below.  
\end{PROOF}
\bigskip

\centerline {$* \qquad * \qquad *$}
\bigskip

We may like in, for example, \ref{3.10} to get not just non-isomorphic
models, but non-isomorphic because of some nice invariant is different. The
following definition serves
\begin{definition}
\label{3.15}
1) Let $\mu$ be a regular uncountable cardinal, $h_0,h_1$ be
functions from some stationary $S\subseteq\mu$ to a set of regular cardinals
$\le \lambda$ satisfying $(\forall \delta\in S)(h_0(\delta)\le h_1(\delta)),
\bar{h}=(h_0,h_1)$.  Let $M$ be a $\tau$-model, $\varphi(\bar x,\bar y)$ a
formula in the vocabulary $\tau$ such that
$\ell g(\bar x)=\ell g(\bar y)=\partial$.

Now, we say that $M$ $\kappa$-obeys $(\bar{h},\varphi)$, or $(h_0,h_1,
\varphi)$, if the following holds:
\mn
\begin{enumerate}
\item[$(*)_0$]  there is a function $\bold H$ from ${}^{\mu>}([M]^{<\mu})$
to $[M]^{<\mu}$ such that:
if $\langle A_i:i<\mu\rangle$ is an increasing continuous sequence of
subsets of $M$, $|A_i|<\mu$, and $\bold H(\langle A_i:i \le j\rangle)\subseteq
A_{j+1}$ for every $j<\mu$, \then \, for some club $\cC \subseteq
\mu$, for  every $\delta\in
\cC \cap S$ of cofinality $\geq\kappa$ the following holds:
\sn
\begin{enumerate}
\item[$\oplus$]  if for each $i<\cf(\delta)$, $\bar a_i\subseteq A_{\alpha_i}$
for some $\alpha_i<\delta$, $\langle\bar a_i:i<\cf(\delta)\rangle$ is weakly
$(\kappa,\varphi(\bar x,\bar y))$-skeleton like inside $M$ 
(so $\ell g(\bar{a}_i)=\partial$), for each $\alpha<\delta$ the sequence
\[
\langle\tp_{\{\varphi,\psi\}}(\bar a_i,A_\alpha):i<\cf(\delta)\rangle
\]
is eventually constant \then:
\sn
\item[$(*)_1=(*)^1_{h_0(\delta),h_1(\delta)}$] if every $B\subseteq |M|$ of
cardinality $<h_0(\delta)$ belongs to $\cP_1$, \then \, every $B\subseteq
|M|$ of cardinality $<h_1(\delta)$ belongs\footnote{so if $h_0(\delta)=
h_1(\delta)$ this is an empty requirement} to $\cP_1$, where
\end{enumerate}
\sn
\item[$(*)_2$]  $\cP_0 = \{B\subseteq M: B\subseteq M$ and
$p^*\rest B$ is realized in $M\;\}$, see on $p^*$ below,
\[
\cP_1 = \{B\in M: B\subseteq M \text{ and } B\cup A_\delta \in
\cP_0\},
\]
\sn
where
\sn
\item[$(*)_3$]  $p^* = p^*_{M,\langle\bar{a}_i:i<\cf(\delta)\rangle}=:
\big\{\vartheta(\bar x,\bar c):\bar c\subseteq M$,
and for every $i < \cf(\delta)$ large enough
$M \models \vartheta[\bar a_i,\bar c]$ \\ 
and
$\vartheta (\bar x,\bar y)\in\{\varphi(\bar x,\bar y),\neg\varphi(\bar x,
\bar y),\varphi(\bar y,\bar x),\neg\varphi(\bar y,\bar x)\}\;\big\}$.
\end{enumerate}
\mn
2) In (1), we say that $M$ obeys $(\bar{h},\varphi(\bar x,\bar y))$
exactly, \when \, in $(*)$, for $\delta\in \cC \cap S$,
the statement $\oplus$ fails for $h_1(\delta)^+$ (i.e., for 
some $\langle\bar a_i:i < \cf(\delta)\rangle$, $p$,
$p^*$ as there, $|p|=h(\delta)$, $p$ is not realized in $M$.)

\noindent
3) We say that $M$ weakly $\kappa$-obeys $(\bar{h},\varphi)$
\when \, the following variant of $(*)$ of part (1) holds:
we replace $(*)^1_{h_0(\delta),h_1(\delta)}$ by
\mn
\begin{enumerate}
\item[$(*)_0=(*)^0_{h_0(\delta),h_1(\delta)}$]  if every $B\subseteq M$ of
cardinality  $<h_0(\delta)$ belongs to $\cP_1$ \then \, every $B\subseteq
M$ of cardinality $<h_1(\delta)$ belongs to $\cP_0$
\end{enumerate}
\mn
4) We say that $M$ weakly obeys $(h_0,h_1,\varphi(\bar{x},\bar{y}))$
exactly \Iff \, in $(*)$ of part (3), for $\delta\in \cC \cap S$, the statement
$(*)^0_{h_0(\delta),h_1(\delta)^+}$ fails.

\noindent
5) We add in the definition above the adjective ``semi'' to
$\kappa$-obeys \Iff \, we change $(*)$ to
\mn
\begin{enumerate}
\item[$(*)'$]  given $\bar b_\alpha\in {}^\partial M$ for $\alpha<\mu$ and
$\langle\bar b_\alpha:\alpha<\mu\rangle$ is weakly 
$(\kappa,\varphi (\bar x,\bar y))$-skeleton like,
there are an unbounded $Y\subseteq \mu$ and
 a function $\bold H$ from $^{}\mu([M]^{<\mu})$ to $[M]^{<\mu}$ such that:
if $\langle A_i:i<\mu\rangle$ is an increasing continuous sequence of
subsets of $M$, $|A_i|<\mu$  and $^{}{i>}\mu \Rightarrow \
\bold H(\langle A_i:i\le j \rangle)\subseteq A_{j+1}$ \then \, 
for some club $\cC$ of $\mu$, for every $\delta \in \cC \cap S$ of 
cofinality $\ge \kappa$, the following holds:
\sn
\begin{enumerate}
\item[$\oplus$]  there are sequences $\langle\alpha_i:i < \cf(\delta)
\rangle$, $\langle\beta_i:i<\cf(\delta)\rangle$ both increasing 
with limit $\delta,\beta_i\in S$, and we let  $\bar a_i=
\bar b_{\beta_i} \subseteq A_{\alpha_i}$ (not necessarily $\langle\bar a_i:i<
\cf(\delta)\rangle$ is weakly $(\kappa,\varphi(\bar x,\bar y))$-skeleton
like inside $M$) and for each $\alpha<\delta$ the sequence $\langle
\tp_{\{\varphi,\psi\}}(\bar a_i, A_\alpha): i<\cf(\delta)\rangle$ is
eventually constant \then\ $(\ast)_0$ or $(\ast)_1$ etc.
\end{enumerate}
\end{enumerate}
\mn
6) We say ``exactly semi $\kappa$-obeys $(h_0,h_1,\varphi)$" \Iff \, $M$
semi $\kappa$-obeys $(h_0,h_1,\varphi)$ and if $\bigwedge\limits_{\delta\in
S} h_1(\delta) \le h^+_1(\delta)$ and $(\exists^{\stat}\delta\in S)(h_1(
\delta)<h^+_1(\delta))$, then $M$ does not semi $\kappa$-obeys $(h_0,h_1^+,
\varphi)$.
We write $(h,\varphi)$ if in $(h_0,h_1,\varphi)$, $h_1=h$ and $h_0$ is
constantly $\kappa$.
\end{definition}

\begin{remark}
\label{3c.15d}
1)  In \ref{3.15}(5), (6) we can avoid $\langle
\alpha_i:i<\cf(\delta)\rangle$ with small changes.

\noindent
2) Note that assuming below $\lambda<\chi^{<\theta}$ is very
reasonable as $\chi^{<\theta}$ is the number of distinct terms, and we
have no information on a representation in $\cM_{\chi,\theta}(I)$ using
every term only once. Also $\lambda<\partial^+$ seems reasonable.
\end{remark}

\begin{theorem}
\label{3c.16}
Assume that $\varphi(\bar x,\bar y)$ is an asymmetric $\tau (K)$-formula,
$\partial = \ell g(\bar x)=\lg(\bar y)$. Suppose that for every 
$I\in K^{\oor}_\lambda$ there is a $\tau$-model $M_I \in K_\lambda$, 
weakly full $\varphi(\bar x,\bar y)$-represented in 
$\cM_{\chi,\theta}(I)$, by the identity
function for notational simplicity (see Definition \ref{2.2}), where
$\lambda>\chi^{<\theta}+\partial^+$ and for $s\in I$, $\bar a_s=\langle
F_{i,1}(s):i<\partial\rangle\in {}^\partial |M_I|$ and
$M_I \models \varphi[\bar a_s,\bar a_t]$ 
\underline{if and only if} $s<t$ \, (for $s,t \in I$)
(where $F_{i,1} \in \tau_{\chi,\theta}$ is  a one place function symbol
for $i < \partial)$.

\Then
\begin{enumerate}
\item[$(a)$]  $\dot{\bbI}(\lambda,K_\lambda)=2^\lambda$ if: $\lambda
\ge \chi^{<\theta}+ \chi^{\partial}$ and
$:\lambda>\chi^{\theta}+\chi^\partial$ or
$\lambda^{\partial}<2^\lambda$ and $\cf(\lambda)>\partial$
or $\lambda^{\partial}<2^\lambda$ and $\theta=\aleph_0$ or there is a linear
order $I$ with $\geq\lambda$ Dedekind cuts of cofinality $\ge \kappa$ with
$2^{|I|}<2^\lambda$,
\sn
\item[$(b)$]  the cardinal invariants from Definition \ref{3.15}(5), suffice to
distinguish $2^\lambda$ models in $K_\lambda$ if $\lambda>\chi^{<\theta}
+\chi^\partial$.
\end{enumerate}
\end{theorem}

\begin{remark}
1) In the cases $M_I=EM_\tau(I,\Phi),|\tau_\Phi| \le \chi,\ell g(\bar a_s)=
\partial$, clearly $M_I$ is weakly full $\varphi(\bar x,\bar y)$-represented
in $\cM_{\chi,\theta}$ by some $f$, $f(\bar a_s)=\langle F_{i,1}(s):
i<\partial\rangle$ for $\theta=\aleph_0$, $\chi=|\tau_\Phi|+\aleph_0$.

\noindent
2)  On ``weakly full $\varphi(\bar{x},\bar{y})$-represented'' see
Definition \ref{2.2} clauses (d)+(f).
\end{remark}

\begin{PROOF}{\ref{3.16}}
Note that, letting $\kappa := \partial^++\theta$, (so it is a 
regular cardinal):
\mn
\begin{enumerate}
\item[$(*)$]  in $M_I$, $\langle\bar a_s:s\in I\rangle$ is weakly $(\kappa,
<\mu,\varphi(\bar x,\bar y))$-skeleton like in $M_I$, see Definition \ref{3.12}
 whenever $\mu \ge \kappa$.
So in particular $(*)$ Definition \ref{3.11} holds.
\end{enumerate}
\mn
[Why? Assume $A\subseteq M_I$ and $|A|<\mu$, so for each $a \in A$ 
let $a=\sigma_a(\bar t_a),\bar t_a \in {}^{\theta>}I$ and let $J =
\cup\{\bar t_a:a\in A\}$ so $J\subseteq I$ is of cardinality $<\mu$
such that $A\subseteq \{\sigma(\bar t):\bar t\in {}^{\theta>}J$ 
and $\sigma$ a $\tau_{\chi,\theta}$-term$\}$.
Clearly $J$ is as required].
\mn
\begin{enumerate}
\item[$(**)$]  $\lambda>\chi^{<\theta}+\partial^+ \ge \kappa=\cf(\kappa)$,
\end{enumerate}
\mn
by the assumption
\mn
\begin{enumerate}
\item[$(***)$]  $\chi \ge \partial$ and of course $\lambda>\chi^{<\theta}
+\partial^+$ hence % 10.06.01 $\lambda\subset\kappa,\lambda>\kappa$.
\end{enumerate}
\mn
We shall use $(*),(**),(***)$, freely.
Let us see why the cases below and \ref{3.11} cover all the possibilities.
\medskip

\noindent 
\underline{Why does clause (a) hold}?

First, if $\lambda>\chi^{<\theta}+\chi^{\partial}$ then clause (b) proved
below suffices, so \wilog \,  $\lambda\leq\chi^{<\theta}+
\chi^{\partial}$, but $\lambda\leq\chi^{<\theta}+\chi^\partial$
so $\lambda=\chi^{<\theta}+ \chi^{\partial}$.

If $\lambda^\partial<2^\lambda$ and $\cf(\lambda)>\partial$ then we can apply
claim \ref{3.11} clause (c); so we have to check the assumptions there.
The general assumption of \ref{3.11}, holds trivially. Now $(*)_0$ there
holds by the general assumption of \ref{3c.16} and $(*)_1$ there holds 
by the case of (a) we are dealing with and $(*)_3$ holds by (*) above.

Second, assume $\lambda^\partial<2^\lambda$ and $\lambda>\chi<\theta=\aleph_0$,
so as \wilog \, the previous case does not holds, we have 
$\cf(\lambda) \le \partial$.

Third, let $\langle\lambda_i:i < \cf(\lambda)\rangle$ be 
strictly increasing with limit $\lambda$, $\lambda_i
=\cf(\lambda_i)> \chi^{<\theta}+\partial^+$, and \wilog \,
$\langle 2^{\lambda_i}: i<\cf(\lambda)\rangle$ is constant (so is constantly
$2^\lambda$) or is strictly increasing (still $2^\lambda = 
\prod\limits_{i<\cf(\lambda)} 2^{\lambda_i}$). In the former case by
Fact \ref{3c.18} below  we can reduce the
problem to any $\lambda_i$, so assume that $\langle 2^{\mu_i}: i<\cf(\lambda)
\rangle$ is strictly increasing. As we are assuming $\chi^{<\theta}<\lambda
\leq \chi^\partial$, clearly $\lambda$ is not strong limit, So without loss of
generality $2^{\lambda_i} \ge \lambda$, and hence $2^{\lambda_i} \ge
\lambda^\partial$, so without loss of generality $2^{\lambda_1} > 
\lambda^\partial$.

Fourth, note that if there is a linear order $I$ with $\ge \lambda$
Dedekind cuts with both cofinalities $\ge \kappa$ and 
$2^{|I|}<2^\lambda$ then we are done as in claim \ref{3.11} clause
(f).  But as $\langle 2^{\mu_i}:i< \cf(\lambda)\rangle$ is 
strictly increasing there is such linear order, see
\cite[3.7=Lc2]{Sh:E62}.
\medskip

\noindent 
\underline{Clause (b)}:

If $\lambda$ is regular $>\kappa^+$, we apply case (C) or case
(F). If $\lambda=\kappa^+$ we apply case (D) (case (G) is empty) and if
$\lambda$ is singular we apply case (E) or (H).
\medskip

\noindent
\underline{Case A}: $\lambda^\partial=\lambda$ or $\lambda^\kappa<2^\lambda$.

As $\kappa=:\partial^++\theta<\lambda$ by (*) above
we can apply \ref{3.11} case (a) or case (b) and get 
$\dot{\bbI}(\lambda,K_\lambda) =2^\lambda$.
\medskip

\noindent
\underline{Case B}: $\lambda^\partial<2^\lambda$
 and $\cf(\lambda)>\partial$ and we get 
$\dot{\bbI}(\lambda,K_\lambda)=2^\lambda$.

By \ref{3.11} case (c) (and $(\ast)$ above).
\medskip

\noindent
\underline{Case C}: $\lambda$ is regular, 
$(\forall \mu < \lambda)[\mu^{<\kappa}<\lambda],\lambda \ge \kappa^{++}$.

Let $S_0=\{\delta<\lambda:\cf(\delta) \ge \kappa\}$ and let $h_0$
be the function with domain $S_0$ and constant value $\chi^{<\theta}$. Let
$J^{[\kappa]}$ be a linear order of cardinality $\kappa$
such that $\alpha<\kappa\Rightarrow J^{[\kappa]} \times (\alpha+1)
\cong J^{[\kappa]} \cong J^{[\kappa]} \times ((\alpha+1)^*)$.
(e.g. let $J$ be a $\kappa$-dense strongly $\kappa$-homogeneous linear order,
hence $\alpha \le \kappa\Rightarrow J \times (\alpha+1)
\cong J= J \times ((\alpha+1)^*)$,
and by the L\"owenheim-Skolem argument there is a dense $J'\subseteq J$ of
cardinality $\kappa$ with this property; alternatively use
\cite[2.21=Lc73]{Sh:E62}).

For a function

\[
h:S_0\longrightarrow\{\mu:\mu\mbox{ is a regular cardinal, }\kappa\leq \mu<
\lambda\}
\]

\mn
let $I_h$ be the linear order with the set of elements

\begin{equation*}
\begin{array}{clcr}
\{(\alpha,\beta,t):&\alpha<\lambda+\kappa,\,t\in J^{[\kappa]} \text{ and}\\
  &\beta < h(\alpha) \text{ if } \alpha\in S_0, \text{ and }
\beta < \kappa \text{ otherwise}\}.
\end{array}
\end{equation*}

\mn
The order is:

\begin{equation*}
\begin{array}{clcr}
(\alpha_1,\beta_1) \le (\alpha_2,\beta_2) \text{ if and only if }&\alpha_1<
\alpha_2, \text{ or}\\
  &\alpha_1 = \alpha_2 \text{ and } \beta_1 \ge \beta_2, \text{ or}\\
  &\alpha_1 = \alpha_2 \text{ and } \beta_1=\beta_2 \text{ and } 
t_1<_{J^*} t_2.
\end{array}
\end{equation*}

\mn
Now
\mn
\begin{enumerate}
\item[$\boxdot$]  $M_{I_h}$ semi $\kappa$-obeys the pair
$(h,(\varphi(\bar x,\bar y))$ exactly (see Definition \ref{3.15}).
\end{enumerate}

First we prove ``obey''.
So (see Definition \ref{3.15}(5) with $\mu=\lambda$) let $\bar b_\alpha\in
{}^\partial(M_I)$ for $\alpha<\lambda$. So for some sequence
$\bar{\sigma}^\alpha$ of $\bar{\sigma}$-terms
$\bar b_\alpha=\bar\sigma^\alpha(\bar t^\alpha)$
with $\bar t^\alpha \in {}^{\kappa>}(I_n) \zeta^*<\kappa,u\subseteq \zeta^*$,
and for some stationary set $Y \subseteq\{\delta<\lambda:\cf(\delta)
=\kappa\}$ and term $\bar \sigma^*$ we have
\mn
\begin{enumerate}
\item[$\circledast_1$]  $\alpha \in Y \Rightarrow 
\bar \sigma^\alpha=\bar \sigma^*,\ell g(\bar t)=\zeta^*$,
order type of  $I \restriction\bar t^\alpha$  is constant, and
$\bar t^\alpha\restriction u=\bar t^*$ and
\sn
\item[$\circledast_2$]  $\epsilon\in\zeta^* \setminus u \Rightarrow$ 
the sequence $\langle t^\alpha_\epsilon:\alpha\in Y\rangle$ is 
$<_I$-increasing
\sn
\item[$\circledast_3$]  the truth value of $t^{\alpha_1}_{\epsilon_1}<_{I_n}
t^{|\alpha_2}_{\epsilon_2}$ for $\alpha_1,\alpha_2 \in Y$ and 
$\epsilon_1,\epsilon_2<\zeta^*$ depend just on the truth values of
$\alpha_1<\alpha_2, \alpha_2<\alpha_1$ and the
values of $\epsilon_1,\epsilon_2$.
\end{enumerate}
\mn
We define a function $\bold H$ from ${}^{\lambda>}([M_I]^{<\lambda})$ to
$[M_I]^{<\lambda}$ by: given $\langle A_j:j<i\rangle$, with $A_j
\subseteq M_{I_h}$ increasing, $|A_j|<\mu$ let

\begin{equation*}
\begin{array}{clcr}
\gamma=\gamma_{A_i} = \Min\{\gamma: &A_j\subseteq \{\sigma^*(\bar t):\bar t\in
{}^{\kappa>}(\gamma\times\mu\times J^{[\kappa]}\cap I_h)\}\text{
and}\\
  &(\forall j \le i)\bar t^j \subseteq \gamma\times \mu\cap I_h\}.
\end{array}
\end{equation*}

\mn
Let $A_i\in [M_{I_h}]^{<\lambda}$ be increasing continuous, 
$\bold H(\langle A_j:
j \le i\rangle)\subseteq A_{j+1}$, and let

\begin{equation*}
\begin{array}{clcr}
\cC = \{\delta<\lambda: &(\forall\alpha,\beta)(\alpha<\delta\cap
(\alpha,\beta)\in I_h \Rightarrow \beta < \delta) \text{ and}\\
  &(\forall i)(\gamma_i<\delta\equiv i<\delta),\text{ and}\\
  &\alpha<\delta \text{ and } i \in Y \setminus \delta \text{ and } j
\in Y \setminus \delta \Rightarrow \\
  &\tp_{\{\varphi,\psi\}}(\bar b_i,A_\alpha,M_{I_h}) =
\tp_{\{\varphi,\psi\}}(\bar b_\delta,A_\alpha,M_{I_n}),\text{ and}\\
  &\delta = \sup(\delta\cap Y) \text{ and}\\
  &\epsilon\in\sigma \setminus u \Rightarrow (\forall i)(t^i_\epsilon
\in \delta\times\mu\times J^{[\kappa]}\equiv i<\delta)\}.
\end{array}
\end{equation*}

\mn
Clearly $\cC$ is a club of $\lambda$. Now let $\delta\in S \cap \cC$.
We can choose $\beta(i)\in Y$ for $i < \cf(\delta)$ increasing with limit
$\delta$. By the definition of representable clearly $\langle\bar
b_{\beta(i)}: i<\cf(\delta)\rangle$ as required from $\langle\bar a_i:i<
\cf(\delta)\rangle$ in Definition \ref{3.15}(1), and so $p^*=p^*_{M_{I_h},
\langle \bar b_{\beta(i)}:i<\cf(\delta)\rangle}$ is well defined.

Now
\mn
\begin{enumerate}
\item[$(*)_0$]  if $B\in [M]^{<h_1(\delta)}$ then $p^*\restriction (A_\delta
\cup B)$ is realized in $N$.
\end{enumerate}
\mn
[Why?  Let $I^*\in [I]^{<h_1(\delta)}$ be such that

\[
B \subseteq\{\sigma(\bar t):\sigma \text{ is a } \tau_{\chi,\theta}
\text{-term and }\bar t\in {}^{\kappa>}(I^*)\}.
\]

\mn
We can find $\beta^*< h_1(\delta)$ such that

\[
(\alpha',\beta',t')\in I' \setminus(\delta\times\delta\times J^{[\kappa]})
\Rightarrow \beta'<\beta^*.
\]

\mn
Now we can choose $\bar t^\otimes\in {}^{\kappa>} I$ such that $\bar
t^\otimes \restriction u=\bar t^*\restriction u$, and

\[
\epsilon\in\partial\setminus u \Rightarrow t^\otimes_\epsilon\in
\{\delta\}\times\{\beta^*\}\times J^{[\kappa]}
\]

\mn
and

\[
epsilon,\zeta<\partial \Rightarrow [t^\otimes_\epsilon < 
t^\otimes_\zeta\equiv t^*_\epsilon< t^*_\zeta],
\]

\mn
possible by the choice of $J^{[\kappa]}$. By ``represented'' and the
definition of $p^*$, clearly $\bar\sigma^*(\bar t^\otimes)$ realizes $p^*
\restriction (A_\delta\cup B)$, so $(*)_0$ holds.]

Now $(*)$ tells us that $M_{I_n}$ semi $\kappa$-obeys $(0,h_1,\varphi(\bar x,
\bar y))$. As for the ``exactly'', it is enough to find $\langle 
\bar b_\alpha:\alpha<\mu\rangle$ exemplifying that, i.e. that for 
every unbounded $S \subseteq \mu,\langle \bar b_\alpha:\alpha \in 
S\rangle$ fulfill the demand there
more then needed it follows by Fact \ref{3c.17} below.
\end{PROOF}

\begin{fact}
\label{3c.17}
Assume
\mn
\begin{enumerate}
\item[$(a)$]  $\mu$ is regular $\leq\lambda$, and $(\forall\alpha<\mu)(\kappa+
\chi+|\alpha|^{<\theta}<\mu)$,
\sn
\item[$(b)$]  $I\in K^{\rm or}_\lambda$,
\sn
\item[$(c)$]  $\langle t_\alpha:\alpha<\mu\rangle$ is $<_I$--increasing,
\sn
\item[$(d)$]  $S =\{\delta<\mu:\cf(\delta)>\kappa\}$ and 
$h$ is the function with domain $S$ defined by $h(\delta)=\cf(I^*
\restriction\{t:(\forall i<\delta) t_i<_I t\})$.
\end{enumerate}
\mn
\Then \, there is a function $\bold H$ from 
${}^{\mu>}([M]^{<\mu})$ to $[M]^{<\mu}$ satisfying
$\bigcup\{\bar a_{t_j}: j<i\}\subseteq \bold H(\langle A_j:j< i\rangle\}$ and
such that: if $A_i\in [M]^{<M}$ is increasing continuous, $H(\langle A_i:j<i
\rangle)\subseteq A_{i+1}$ and

\[
\cC = \{\delta<\mu:\delta \text{ a limit ordinal such that }
(\forall i<\mu)(\bar a_{t_i}\subseteq A_\delta \Leftrightarrow i<\delta\},
\]

\mn
\then
\mn
\begin{enumerate}
\item[$(\alpha)$]  $\cC$ is a club of $\mu$,
\sn
\item[$(\beta)$]  there is an increasing continuous sequence $\langle
I_\alpha:\alpha<\mu\rangle$, $I_\alpha\subseteq I$, $|I_\alpha|<\mu$ such
that
\sn
\begin{enumerate}
\item[$(i)$]  $A_\alpha\subseteq\{\sigma(\bar t):\sigma$ an
$\tau_{\chi,\theta}$-term, $\bar t\in {}^{\theta>}(I_{\alpha+1})\}\subseteq
A_{\alpha+1}$,
\sn
\item[$(ii)$]  $t_\alpha\in I_{\alpha+1}$,
\sn
\item[$(iii)$]  $\cC_1 = \{\delta\in \cC$: if $t_\alpha\in 
I_\delta$ and $(\exists \beta)(t<_I t_\beta) \Rightarrow 
(\exists\beta<\delta)(t<_I t_\beta)\}$ is a club of $\mu$,
\sn
\item[$(iv)$]  $\bar a_{t_\alpha}\in A_{\alpha+1}$
\end{enumerate}
\sn
\item[$(\gamma)$]  if $\delta\in \cC \cap S$ there are $\langle
\alpha_\epsilon:\epsilon<\cf(\delta)\rangle$, $\langle\beta(\epsilon):
\epsilon<\cf(\delta)\rangle$ increasing with limit $\delta$,
such that $\bar a_{t_{\beta(\epsilon)}}\subseteq A_{\alpha_\epsilon}$,
\sn
\item[$(\delta)$]  if $\delta,\langle\alpha_\epsilon,\beta(\epsilon):
\epsilon<\cf(\delta)\rangle$ are as in clause $(\beta)$ \then\ for each
$\alpha<\delta$ the sequence $\langle\tp_{\epsilon,\phi}
(\bar a_{t_{\beta(\epsilon)}},A_\alpha,M):\epsilon<\cf(\delta)\rangle$ is
essentially constant,
\sn
\item[$(\epsilon)$]  if $B\subseteq M$, $|B|<\cf(\delta)+h(\delta)$
\then \, $p^*_{M,\langle\bar a_{t_{\beta(i)}}:i<\cf(\delta)\rangle}
\rest B$ is realized in $M$, see Definition \ref{3.15}(1),$(*)_3$, 
so in Definition \ref{3.15}(1),$(*)_2$'s notation, 
$[M]^{<(\cf(\delta)+ h(\delta))}\subseteq \cP_0$,
\sn
\item[$(\zeta)$]  if $B\subseteq M$, $|B|<h(\delta)$ \then\ $p^*_{M,\langle
\bar a_{t_{\beta(i)}}:i<\cf(\delta)\rangle}\restriction (B\cup A_\delta)$ is
realized in $M$, so in Definition \ref{3.15}(1)$(*)_2$'s notation, $[M]^{<
h(\delta)}\subseteq \cP_1$a
\sn
\item[$(\eta)$]  there are $B^- \subseteq A_\delta$ of cardinality
$\cf(\delta)$ and $B^+\subseteq M$ of cardinality $h(\delta)$ such that
$p^*_{M,\langle \bar a_{t_{\beta(i)}}:i<\cf(\delta)\rangle}\restriction
(B^-\cup B^+)$  is omitted by $M$,
actually $\{\varphi(\bar a_{t_{\beta(i)}},\bar x):i<\cf(\delta)\}\cup
\{\varphi(\bar x, a_t):t\in J\}$ is omitted for some $J\in
[I]^{\cf(\delta)}$.
\end{enumerate}
\end{fact}
\medskip

\begin{PROOF}{\ref{3c.16}}

\noindent 
\underline{Continuation of the proof of Theorem \ref{3c.16}}.
\smallskip

\noindent
\underline{Case D}: $\lambda=\kappa^+>\chi^{<\theta}$.

Similar to Case C, but we have to allow $h(\delta)$ to be $\kappa^+=
\lambda$ in addition to $\kappa$. So $I_h$, defined similarly using $J^{[
\lambda]}$ (not $J^{[\kappa]}$), is no longer $\lambda$-like, $\bar
b_\alpha \in {}^\partial (M_{I_h})$, if the rest is not obvious look at the
proof of Case E.
\smallskip

\noindent
\underline{Case E}:  $0<\gamma^*$, $\chi^{<\kappa}+|\alpha|
<\mu_i<\lambda$, $\mu_i$ ($i<\alpha^*$) strictly increasing, each $\mu_i$
regular, $\mu_{i+1}>\mu^{+++}_i$, $\mu_i>\chi+\partial^++\theta$, $(\forall
\mu<\mu_i)\mu^{<\kappa}<\mu_i$, $\prod\limits_i 2^{\mu_i}=2^\lambda$
(without the last assumption we just get a smaller number of models; note
that if $(\forall \alpha<\lambda)(\chi+|\alpha|^{<\kappa}<\lambda)$, then
there is such $\langle\mu_i:i<\alpha\rangle$).

Let $J^i\cong J^{[\mu^{++}_i]}$ for $i<\alpha^*$ be from Fact \ref{3c.17}
below, and for each $i<\gamma^*$ define $J_h\in K^{\oor}_{\mu^{+3}_i}$ for
$h:\{\delta<\mu^{+++}_i:\cf(\delta)=\mu^{++}_i\}\longrightarrow\{\mu^+_i,
\mu^{++}_i\}$ to be $\sum\limits_{\zeta<(\mu^{+3}_i+\kappa)}(J^i_\zeta)^*$,
where: $\mu ^{+3}_i+\kappa$ is ordinal addition, the $J^i_\zeta$ are
pairwise disjoint, $J^i_\zeta$ is isomorphic to $J^i$ except when $h(\zeta)$
is well defined and equal to $\mu^+_i$, then $J^i_\zeta$ is isomorphic to
$J^i\times(\mu^+_i)^*$. 

Lastly, for every

\[
\bar h\in\prod\limits_i\big\{h:\Dom(h)=S_i=:\{\delta<\mu^{++}_i:\cf(\delta)=
\mu^+_i\},\quad h\mbox{ as above }\big\},
\]

\mn
we let $I_{\bar h}=:\sum\limits_i J_{h_i}+\lambda\times J^{[\kappa]}$.

For each $i<\alpha$ we have to prove that $h_i/{\mathscr D}_{\mu_\gamma^{+++}}$
is an invariant of the isomorphic type of $M_{I_{\bar h}}$. For this it is
enough to prove, for each $\gamma_*<\gamma^*$, that
\mn
\begin{enumerate}
\item[$(*)$]  $M_{I_{\bar h}}$ exactly semi $\kappa$--obeys $(0,h_{\gamma_*}
\varphi)$.
\end{enumerate}
\mn
It is enough to prove ``semi $\kappa$--obeys $(0,h_\gamma,\varphi)$'', as
then the exactness follows by Fact $\alpha$ above. Let $\bar b_\alpha\in
{}^\partial (M_{I_{\bar h}})$ for $\alpha<\mu^{+++}_\gamma$, so $\bar
b_\alpha = \bar\sigma^\alpha(\bar t^\alpha)$, $\bar t^\alpha\in {}^{\kappa>}
(I_{\bar h})$. We can find a stationary set $Y\subseteq\{\delta<\mu^{+++}_{
\gamma_*}:\cf(\delta)=\kappa\}$ such that

\[
\alpha\in Y \Rightarrow \bar\delta^\alpha=\sigma^* \wedge 
\ell g(\bar t^\alpha)=\epsilon^*,
\]

\mn
as $\{(\epsilon,\zeta):t^\alpha_\epsilon<\zeta^\alpha_\zeta\}=v,
u_{i,\gamma}=\{\epsilon<\epsilon^*:t^\alpha_\epsilon\in
J_{h_i}\}=u_\gamma$. By clauses (i)+(h), without loss of generality $\langle
\bar t^\alpha:\alpha\in Y\rangle$ is order indiscernible, as in the proof of
Case C.

So for each $\epsilon<\epsilon^*$, $\langle t^\alpha_\epsilon:
\alpha\in Y\rangle$ is constant, or strictly increasing, or strictly
decreasing, and for some $\gamma<\gamma^*$ they are all in on one 
$I_{h_\gamma}$, moreover if $\langle t^\alpha_\epsilon:\alpha\in 
Y_\epsilon\rangle$ is not constant necessarily $\gamma\geq\gamma_*$. 
So if $\langle
t^\alpha_\epsilon:\alpha\in Y\rangle$ is strictly increasing, $\delta<
\mu^{+++}_{\gamma_*}$, $\cf(\delta)=\mu^+_i$, then

\[
\cf(I^*_{\bar h}\restriction\{t:t<t^*_\kappa\mbox{ for every }\alpha\in Y\})
\]

\mn
is $\mu^+_\gamma$ or $\mu^{++}_\gamma$ when $\epsilon\in u_\gamma$, so is
$\ge \mu^{++}_{\gamma_*}$ except when $\epsilon\in u_{\gamma_*}$ and
$h(\delta)=\mu^+_i$. The situation is similar when $\langle
t^\alpha_\epsilon:\alpha\in Y\rangle$ is strictly decreasing, except that
now $\epsilon \in u_{\gamma_*}$ is impossible.
\smallskip

\noindent
\underline{Case F}: $\lambda$ is regular $>\chi^{<\theta}+
\chi^\partial+\kappa^+$, without loss of generality
$\lambda>(2^\partial)^+$.

(Why the without loss of generality? Otherwise Case C applies.)

\noindent 
First proof: 

Let

\[
S=\{\delta<\lambda:\cf(\delta)=(2^\partial)^+\},
\]

\mn
and for $h:S \longrightarrow \Reg \cap [\kappa,\lambda)$ we define $I_h$ as in
Case C. It suffices to prove
\mn
\begin{enumerate}
\item[$(*)$]  $M_{I_h}$ exactly semi $\kappa$-obeys $(0,h,\varphi)$.
\end{enumerate}
\mn
It suffices to prove $M_{I_h}$ semi $\kappa$-obeys $(0,h,\varphi)$ as the
exactly follows by Fact $\alpha$. Let $\bar b_\alpha\in {}^\partial(M_{I_h})$
for $\alpha<\lambda$ be such that $\langle\bar b_\alpha:\alpha<\lambda
\rangle$ is $(\kappa,\varphi)$-skeleton like and let $\bar b_\alpha=\bar
\sigma^\alpha(\bar t^\alpha)$, and we choose a stationary set $Y_0\subseteq
\{\delta<\lambda:\cf(\delta)=\kappa\}$ such that $\alpha\in Y \Rightarrow
\sigma^\alpha=\sigma^*$ and $\{(\epsilon,\zeta):t^\alpha_\epsilon < 
t^\alpha_\zeta\}=v,\ell g(t^\alpha)=\epsilon^*<\kappa$ (but no
$\Delta$-system!).

Let $\langle A_i: i<\lambda\rangle$, $\langle I_i=(\gamma_i\times\delta
\times J^*)\cap I:i<\lambda\rangle$, $\cC$ be as there.

For $\delta\in S \cap \acc(C)$ let $Y_1\subseteq Y\cap\delta \cap \cC$ be
unbounded of order type $\cf(\delta)$, and $Y_2\subseteq Y_1$ be unbounded
and $\langle t^\alpha:\alpha\in Y_2\rangle$ be indiscernible (for $<_I$)
(exists as $\otp(Y_1)=(2^\partial)^+$).

Let

\[\begin{array}{l}
u_0=\{\epsilon<\epsilon^*:\langle t^\alpha_\epsilon:\alpha\in
Y_2\rangle \mbox{ is constant}\},\\
u_1=\{\epsilon<\epsilon^*:\langle t^\alpha_\epsilon:\alpha\in
Y_2\rangle \mbox{ is increasing and }(\forall\beta<\delta)(\exists\alpha\in
Y_2)(t^\alpha_\epsilon \notin I_\beta)\},\\
u_2=\epsilon^*\setminus u_0\setminus u_1.
\end{array}\]

\mn
Choose $\beta_0<\beta_1<\beta_2$ in $Y_2$ such that $\{t^\alpha_\epsilon:
\alpha\in Y_1,\ \epsilon \in u_2\cup u_0\}\subseteq I_{\beta_0^*}$.

For each $\beta\in Y_2\setminus\beta_2$ define $\bar s^\beta\in
{}^{\epsilon^*}I$, $\bar s^\beta\restriction u_0=\bar t^\alpha
\restriction u_0$ for $\alpha \in Y_2$, $\bar s^\beta\restriction u_1= \bar
t^\beta\restriction u_1$, $\bar s^\beta\restriction u_2=\bar t^{\beta_2}
\restriction u_2$. Now we can continue as in Case C when we note
\mn
\begin{enumerate}
\item[$(\otimes)$]  if $\beta_3<\beta_4$ are from $Y_2\setminus\beta_2$ then
$\bar\sigma^*(\bar t^{\beta_4})$, $\bar\sigma^*(\bar s^{\beta_4})$ realize
the same $\{\varphi,\psi\}$--type over $A_{\beta_3}$.
\end{enumerate}
\mn
[Why? Let $\bar d\in {}^\partial (A_{\beta_3})$ so $\bar d=\bar\sigma'(\bar
t')$, $\bar t\in {}^{\kappa>}(I_{\beta_3})$. If, e.g.,

\[
M_{I_h}\models\vartheta[\bar\sigma^*(\bar t^{\beta_4}),\bar d)\equiv\neg
\vartheta [\bar\sigma^*(\bar s^{\beta_4}),\bar d]
\]

\mn
then

\[
M\models\vartheta[\bar\sigma^*(\bar t^{\beta_4}),\bar\sigma'(\bar t')]
\equiv\neg\vartheta[\bar\sigma^*(\bar t^{\beta_4}),\bar\sigma'(\bar t')].
\]

\mn
So $(\otimes)$ holds.

Now we can find $\bar t^{\prime\prime}\in {}^{\kappa >}(I_{\beta_1})$ such
that $\bar t'',\bar t'$ realizes the same quantifier free
type (in $I$!) over $I_{\beta_0}$, hence over ($\bar t^{\beta_4}\rest
(u_0\cup u_2)) \char 94 \bar t^{\beta_2}\restriction (u_0\cup u_2)$. Hence

\[
M_{I_h}\models\vartheta [\bar\sigma^*(\bar t^{\beta_4}),
\bar\sigma'(\bar t'')] \equiv \neg \vartheta[\bar\sigma^*(\bar s^{\beta_4}),
\bar \beta'(\bar t'')].
\]

\mn
Similarly $\bar s^{\beta_4}$, $\bar t^{\beta_2}$ realize the same quantifier
free type (in $I$) over $I_{\beta_1}$, hence
\[M_{I_h} \models \vartheta[\bar\sigma^*(\bar s^{\beta_4}),
\bar\sigma'(\bar t'')] \equiv\vartheta [\bar\sigma^*(\bar t^{\beta_2}),
\bar\sigma'(\bar t''),
\]

\mn
so together

\[
M_{I_h}\models\vartheta[\bar\sigma^*(\bar
t^{\beta_4}),\bar\sigma'(\bar t'')] \equiv
\neg\vartheta [\bar\sigma^*(\bar t^{\beta_2}),
\bar\sigma'(\bar t'')].
\]

\mn
But this contradicts the choice of $\cC$ (as $Y\subseteq \cC$).
\medskip

\noindent
\underline{Second proof}:

Similar to case $C$ using \cite[3.7=Lc2]{Sh:E62}.
\smallskip

\noindent
\underline{Case G}:$\lambda$ is regular $>\chi^{<\theta} + \chi^\partial$.

If cases  (C) + (F) do not occur then $\lambda=\kappa^+$, so case D applies.
\medskip

\noindent
\underline{Case H}: $\lambda$ is singular $>\chi^{<\theta}+
\chi^\partial$ (hence $>(2^\partial)^+$).

Combine the proof of cases E and F.
\end{PROOF}

\begin{fact}
\label{3c.18}
Assume $\chi\leq\mu=\mu^{<\theta}< \lambda$ and the linear order
$J^{[\lambda]}$ are from \cite[2.21=Lc73]{Sh:E62} with
$(\mu,\mu^+,\mu^+,\aleph_0)$ here standing for
$(\lambda,\mu_1,\mu_2,\theta)$ there and for $I\in K^{\oor}_\mu$ 
we define $M_I$ naturally, as $M_{I+J^{[\lambda]}}\restriction\{\sigma(\bar
t):\sigma$ a $\tau_{\chi,\theta}$-term, $\bar t\in
{}^{\theta>}(I+J^{[\mu]})\}$ (using the fullness of the
representations). 

\Then
\mn
\begin{enumerate}
\item[$\boxtimes_1$] if $I_1$, $I_2\in K^{\rm or}_\mu$, and
$M_{(I_1+J^{[\mu]})}\not \cong M_{(I_2+J^{[\mu]})}$,
then $M_{(I_1+J^{[\lambda]})} \not\cong M_{(I_2 +J^{[\lambda]})}$, so
$M'_I\not\cong M'_J$ 
\sn
\item[$\boxtimes_2$] $|\{M_I/\cong:I\in K^{\rm or}_\lambda\}|\ge |\{M_{I+
J^{[\mu]}}/\cong: I\in K^{\rm or}_\mu\}|=|\{M'_I \not\cong: 
I\in K^{\oor}_M\}|$.
\end{enumerate}
\end{fact}

\begin{proof}
The first clause by clause (j) of \cite[2.21=Lc73(4)]{Sh:E62} below,
the second clause follows.
\end{proof}
\bigskip

\centerline {$* \qquad * \qquad *$}

\begin{remark}
\label{3.1 17}
Note that if we use strongly $\kappa$-homogeneous $J^{[\kappa]}$ and $M_I$
is weakly fully represented in $\cM_{\chi,\theta}(I)$ then this form of $I$
helps to ``eliminate quantifiers" is $\cM_{\chi,\theta} (I)$, i.e.
$\tp(\bar \sigma,\bar t), \emptyset, M_I)$ is determined by 
$\bar \sigma$ and the order of $\bar t$ if 
$\bar t\in {}^{\kappa>}I$. The order $I^{[\kappa]}$ is not really so
homogeneous but it close too, see \cite[\S2]{Sh:E62}.
\end{remark}

\begin{claim}
\label{3.27new}
In the theorems above in the 
assumption we can restrict ourselves to linear order $I$ satisfying
\mn
\begin{enumerate}
\item[$(*)_I$]  $(a) \quad$ for every infinite $J\subseteq I$, the 
number of Dedekind cuts of $J$ realized by elements of $I$ is at 
most $|J|$ (i.e., stable in $\theta$ for every $\theta$),
\sn
\item[$(b)$]  for every infinite $J_0\subseteq I$ there is an $J_1$,
satisfying $J_0 \subseteq J_1 \subseteq I$ such that 
$|J_0|=|J_1|$ and: if $s,t\in I \setminus J_1$ realize the same 
Dedekind cuts of $J_1$ \then \, there is an automorphism $h$ of $I$ over
$J_1$ (i.e. $h \rest J_1 = \id_{J_1})$ mapping $s$ to $t$
(i.e., almost homogeneous for every $\theta$).
See Definition \cite[2.15=Lb56]{Sh:E62} and \cite[2.16=Lb60]{Sh:E62}.
\end{enumerate}
\end{claim}

\begin{PROOF}{\ref{3.27}}
By \ref{3.29new}.
\end{PROOF}

\noindent
We may weaken a little the definition of weakly $\kappa$-skeleton like
(Definition \ref{3.1}(1)).
\begin{definition}
\label{3.28new}
1)  We say $\langle\bar a_s:s\in I\rangle$ is pseudo $\kappa$-skeleton
like for $\Lambda$ \when \,: for every $\varphi(\bar x,\bar a)\in\Lambda$ and a
Dedekind cut $(I_0,I_1)$ of $I$ such that 
$I_1 \ne \emptyset \Rightarrow \cf(I_1) \ge \kappa \qquad \text{ and }
I_2 \ne \emptyset \Rightarrow \cf(I_2^*)\ge \kappa$ there 
are $J_0$, $J_1$ such that 
\mn
\begin{enumerate}
\item[$(*)_1$]  $J_0$ is an end segment of $I_0$ non empty if $I_0\ne
\emptyset$,
\sn
\item[$(*)_2$] $J_1$ is an initial segment of $I_1$, non empty if $I_1\ne
\emptyset$,
\sn
\item[$(*)_3$]  if $s,t\in J_0\cup J_1$ then $M\vDash\varphi[\bar a_s,\bar a]
\equiv \varphi[\bar a_t, \bar a]$; clearly this is a weaker demand than the
``weakly" version.
\end{enumerate}
\mn
2) Similarly we adopt Definition \ref{3.1}(2),(4).
\end{definition}

What is the difference? E.g., for $\kappa=\aleph_0$, $J_{\bar a}$ instead of
being countable it may be a Suslin order or Specker order.
\begin{claim}
\label{3.29new}
We can through all this section ask (a) or (a)+(b) or (a)+(b)$'$, where
\mn
\begin{enumerate}
\item[$(a)$]  replace weakly in ``weakly $\ldots$ skeleton likeq" by pseudo
(including the definitions) and all claims remain true;
\sn
\item[$(b)$]  restricting ourselves to $\lambda \ge 
2^{<\kappa}$, we can replace linear orders by strongly 
$\kappa$-dense linear order (see below);
\sn
\item[$(b)'$]  we can demand that all our linear orders are $\theta$-stable
and almost $\theta$-homogeneous, see Definition \cite[2.21=Lc73]{Sh:E62}.
\end{enumerate}
\end{claim}

\begin{definition}
\label{3.30new}
1) A linear order $I$ is $\kappa$-homogeneous if $\cf(I) \ge \kappa$,
$\cf(I^*) \ge \kappa$ for any subsets $J_0$, $J_1$ of $I$ of cardinality
$<\kappa$ (possibly empty) satisfying
$(\forall s_0\in J_0)(\forall s_1\in J_1)(s_0<_I s_1)$
there is $t\in I$ such that
$(\forall s_0 \in J_0)(s_0<_I t)$ and $(\forall s_1\in J_1)(t<_I
s_1)$.

\noindent
2) A linear order $I$ is strongly $\kappa$-dense if it is
$\kappa$-dense and every partial one-to-one function from $I$ to $I$ of
cardinality $<\kappa$ can be extended to an automorphism.

\noindent
3) A linear order $I$ is $\theta$-stable if for every $J\subseteq I$
of cardinality $\leq\theta$, the number of Dedekind cuts of $J$ induced by
elements of $I$ is at most $\bar \theta$.
\end{definition}

\begin{PROOF}{\ref{3.30new}}
Straightforward, we rely on \cite[2.21=Lc73(5)]{Sh:E62}.
\end{PROOF}

\bibliographystyle{alphacolon}
\bibliography{lista,listb,listx,listf,liste,listz}

\end{document}